\documentclass{article}

\usepackage{pgf,tikz}
\usetikzlibrary{arrows}
\usepackage{array,multirow,makecell}
\setcellgapes{1pt} \makegapedcells
\usepackage{amssymb}
\usepackage{subfig}
\usepackage[normalem]{ ulem }
\usepackage{soul}
\usepackage{graphicx}
\usepackage{here}
\usepackage{epsfig}
\usepackage{color}
\usepackage{pgf,tikz,amsthm}
\usetikzlibrary{arrows}
\usepackage{amsmath}
\usepackage{epstopdf}

\usepackage{hyperref}
\hypersetup{
    colorlinks=true,
    linkcolor=blue,
    filecolor=magenta,      
    urlcolor=cyan,
    citecolor=red
}

\newtheorem{theorem}{\textbf{Theorem}}
\newtheorem{lemma}{\textbf{Lemma}}
\newtheorem{proposition}{\textbf{Proposition}}

\newtheorem{remark}{\textbf{Remark}}

\newcommand{\R}{\mathcal{R}}
\newcommand{\N}{\mathcal{N}}

\newcommand{\Q}{\mathcal{Q}}

\newcommand{\rr}{\mathbb{R}}

\newcommand{\comment}[1]{}

\usepackage{tikz}
\usepackage{anyfontsize}
\usetikzlibrary{arrows.meta}
\tikzset{%
  >={Latex[width=2mm,length=2mm]},
            base/.style = {rectangle, rounded corners, draw=black,
                           minimum width=1.5cm, minimum height=1.5cm,
                           text centered, font=\sffamily},
  activityStarts/.style = {base, fill=blue!30},
  startstop1/.style = {base, fill=blue!30},
  startstop2/.style = {base, fill=green!30},
  startstop3/.style = {base, fill=orange!30},
       startstop/.style = {base, fill=red!30},
    activityRuns/.style = {base, fill=green!30},
         process/.style = {base, minimum width=2.5cm, fill=orange!15,
                           font=\ttfamily},
}

\title{About contamination by sterile females and residual male fertility on the effectiveness of the sterile insect technique. Impact on disease vector control and disease control.}
\author{Y. Dumont$^{1,2,3}$ \footnote{Corresponding author: yves.dumont@cirad.fr}
, I.V. Yatat-Djeumen$^{4,5,6}$ \\
\small $^1$CIRAD, Umr AMAP, P\^ole de Protection des Plantes, F-97410 Saint Pierre, France, \\
\small $^2$AMAP, Univ Montpellier, CIRAD, CNRS, INRA, IRD, Montpellier, France \\
\small $^3$University of Pretoria, Department of Mathematics and Applied Mathematics, Pretoria, South Africa \\
\small $^4$ CIRAD, Umr AMAP, F-34398 Montpellier, France \\
\small $^5$University of Yaound\'e I, National Advanced School of Engineering of Yaound\'e, \\ 
\small Department of Mathematics and Physics, Yaound\'e, Cameroon \\
\small $^6$UMI 209 IRD/UPMC UMMISCO, Bondy, France}
\date{\today}

\begin{document}

\maketitle

\begin{abstract}

The sterile insect technique (SIT) is a technique to control pests and vectors of diseases by releasing mainly sterile males. Several challenges need to be solved before large-scale field application in order to guarantee its success. In this paper we intend to focus on two important issues: residual (sterile) male fertility and contamination by sterile females. Indeed, sterile males are never $100\%$ sterile, that is there is always a small proportion, $\varepsilon$, of fertile males (sperm of) within the sterile males population. Among the sterile insects that are released, a certain proportion, $\epsilon_F$, of them are sterile females due to an imperfect mechanical sex-separation technique. This can be particularly problematic when arthropod viruses are circulating, because mosquito females, even sterile, are vectors of diseases.

Various upper bound values are given in the entomological literature for $\epsilon_F$ and $\varepsilon$ without clear explanations. In this work, we aim to show that these values are related to the biological parameters of the targeted vector, the sterile insects release rate, and the epidemiological parameters of a vector-borne disease, like Dengue. We extend results studied separately in \cite{Aronna2020,DumontYatat2022}.

To study the impact of both issues, we develop and study a SIT-entomological-epidemiological mathematical model, with application to Dengue. Qualitative analysis of the model is carried out to highlight threshold values that shape the overall dynamics of the system.

We show that vector elimination is possible only when $\N \varepsilon<1$, where $\N$ is the basic-offspring number related to the targeted wild population. In that case, we highlight a critical sterile males release rate, $\Lambda_M^{crit}$, above which the control of the wild population is always effective, using a strategy of massive releases, and then small releases, to reach elimination and nuisance reduction. In contrary, when $\N \varepsilon>1$, then SIT-induced vector elimination is unreachable, whatever the size of the releases. 

Moreover, we compute a critical value for the release rate of sterile females, $\Lambda_F^{crit}$, such that if the release rate of the sterilized females is greater than $\Lambda_F^{crit}$, then the epidemiological risk increases. When the sterile females releases rate is low, less than $\Lambda_F^{crit}$, then whatever the value taken by $\varepsilon\N$, the epidemiological risk can be controlled using SIT. However, this is more difficult when $\N \varepsilon>1$. We illustrate our theoretical results with numerical simulations, and we show that early SIT control is better to prevent or mitigate the risk of an epidemic, when residual fertility and contamination by sterile females occur simultaneously. We also highlight  the importance of combining SIT with mechanical control.

In order to guarantee the success of SIT control, we recommend to solve in priority the issue of residual fertility, and, then, to decay the contamination by sterile females as low as possible.

\end{abstract}

\section{Introduction}
Vector-borne diseases have become very strong issues all around the World. After decades of chemical control, the use of biological control methods are more than necessary. Many research programs are ongoing to develop new biocontrol tools. Among them, an old control technique, the Sterile Insect Technique (SIT), is always under study and improvements \cite{SIT,Lee2015}. SIT is an environmentally safe, cost-effective, species-specific, and efficient method of insect control. It is a form of insect population control that relies on the mass-rearing and sterile release of large numbers of male insects to mate with wild female insects. This prevents the production of viable eggs, thus reducing the overall population of the target species. This technique was first developed in the 1950s by entomologists Edward Knipling and Raymond Bushland, who were working for the U.S. Department of Agriculture (USDA) \cite{Knipling1955} (see also \cite{SIT}[chapter 1.1]). The original purpose of SIT was to control the screwworm fly, which was devastating the cattle industry in the southern United States \cite{Knipling1957}. Since then, SIT has been used to control a variety of other insect pests, including the Mediterranean fruit fly, tsetse fly, and also against vectors of diseases, including anopheles and aedes mosquitoes, with more or less success \cite{SIT}. Initially, sterile insects were obtained only by ionization or irradiation, but now new techniques have been developed for mosquitoes control in particular. One of them consists of releasing only males carrying the bacteria \textit{Wolbachia} \cite{sinkin2004}. This is called the Incompatible Insect Technique (IIT) \cite{Lee2015}, where the sperm of Wolbachia-carrying males, W-males, is altered so that it can no longer successfully fertilize uninfected eggs. Thus, IIT can be seen as a classical SIT. A third method exists but it is more controversial since it relies of genetic-modified mosquitoes: this is called the RIDL method, where RIDL stands for "Release of Insects carrying Dominant Lethals" \cite{RIDL_2000}.

However, while conceptually very simple, the conditions and the difficulties to implement SIT in the field are numerous and that is why a drastic control quality is needed. To this end, IAEA (the International Atomic Energy Agency) has published several manuals where several control steps have to be checked in order to ensure/maximize the success of SIT \cite{FAO-IAEA-2020,WHO-IAEA-2020,Clelia2021}.

While several field programs are ongoing, very few have a mathematical modelling component involved. This is a pity because mathematical modelling can bring new insights on several issues that can be detrimental to the efficacy of SIT: see, for instance, \cite{Anguelov2020,Aronna2020,Bliman2022,DumontYatat2022}, and references therein.

Among these controls, it is necessary to evaluate an upper bound for the contamination by sterile females, i.e. the maximal amount of sterile females that can be released during each field release in order to insure that SIT is efficient. Indeed, in order to produce sterile males only, it is necessary to separate the females from the males. Up to now, the sex-separation system is mechanical as male nymphs are (in general) smaller than female nymphs. However, since sex-sorting is highly operator-dependent, a certain number of female nymphs can accidentally fall in the male nymphs bucket and, then, be irradiated to become fully sterilized. Thus, when sterile mosquitoes are released, if the amount of released sterile females is too large, this could maintain or increase the epidemiological risk. Moreover, when the Incompatible Insect technique is considered, releasing Wolbachia-carrying females, even a small amount, can induce a population replacement as showed in \cite{bliman2023}.

For \textit{Aedes albopictus}, estimates of contamination by sterile females, done in Mauritius island \cite{Iyaloo2020b}, were around $4\%$, while in a recent SIT program in R\'eunion island estimates were around $1\%$. Note also carefully that sterilized females are always $100\%$ sterile and thus cannot participate in the wild insect dynamics.
In \cite{DumontYatat2022}, we have showed that when no vector-borne viruses are circulating, then the release of sterile females is not an issue, as long as enough sterile males are released. When a virus is circulating, we showed existence of a contamination threshold for sterile females, such that if the amount of released sterile females per hectare is lower than this threshold, then it is possible to control the wild mosquitoes population. Otherwise, whatever the size of the releases, the basic reproduction number 
will always be greater than $1$ and thus it will be impossible to control the epidemiological risk even if the wild population has been reduced using massive sterile insects releases.

Another control test to take care is the (sterile) male residual fertility, when sterilized males (sperm of) are not necessarily $100\%$ sterile, even if an optimal dose of radiation is used.  Indeed, males are sterilized in boxes such that full sterility is not guaranteed: There are always irradiated males with a small amount of sperm that remains fertile. This is called residual fertility. For \textit{Aedes albopictus}, some estimates done in Mauritius \cite{Iyaloo2020a} lead to a residual fertility between $3.8\%$ and $4.1\%$, while in the SIT-program in R\'eunion island, an average value of $1\%$ was obtained. In Italy, in \cite{Bellini2013}, the authors found a residual fertility between $0.82 \pm 0.14 \%$ and $4.93 \pm 4.72 \%$ thanks to the age of the males, for an irradiation at $40$ Gy.

In \cite{Aronna2020}, using a very simple model, the authors showed that the proportion of fertile sperms, $\varepsilon$, has to be lower than $1/\mathcal{N}$, where $\mathcal{N}$ is the basic offspring number related to the targeted wild population. If, for any reason, $\varepsilon>1/\mathcal{N}$, then, whatever the amount of sterile males released, the wild population will always be above a threshold, that can be estimated, numerically at least.

Up to know we have studied these two issues separately in \cite{Aronna2020,DumontYatat2022}, while, in fact, they do occur simultaneously. Thus, it would be useful to know how the combination of both issues could be problematic in the implementation of SIT program either for nuisance reduction or to reduce the epidemiological risk.


The paper is organized as follows. In section \ref{section2}, we present the full SIT-entomological-epidemiological model and we recall theoretical results without SIT obtained in \cite{Anguelov2020,DumontYatat2022} and we derive theoretical results for the SIT-entomological model. The full SIT model is studied in section \ref{section3}. Finally, in section \ref{section4}, we derive some numerical simulations to illustrate our theoretical results and to discuss the impact of low/high residual fertility as well as low/high contamination by sterile females. The paper ends with a conclusion and perspectives in section \ref{conclusion}.

\section{The SIT-entomological-epidemiological Model \label{section2}}

Based on \cite{DumontYatat2022}, we briefly describe the full SIT model, taking into account residual male fertility and contamination by sterile females.

From the entomological point of view, we split the mosquito population into immature stage (larvae and pupae), $A$, male adults, $M$, and mature females, $F_W$. 

We consider $\Lambda_{tot}$ the release rate of all sterile insects, i.e. sterile males and sterile females, such that $\Lambda_{tot}=\Lambda_M+\Lambda_F$, where   $\Lambda_M=(1-\epsilon_F)\Lambda_{tot}$, $\Lambda_F=\epsilon_F \Lambda_{tot}$, and  $\epsilon_F$ is the proportion of sterile females released. 

Male residual sterility is modeled by considering that a proportion, $\varepsilon M_S$, of sterile males is fertile, such that emerging immature females will become fertile with a probability of $\dfrac{M+\varepsilon M_S}{M+M_S}$ or they will become sterile with a probability of $\dfrac{\varepsilon M_S}{M+M_S}$.

\begin{figure}[h!]
\centering
\begin{tikzpicture}[node distance=2.0cm,
    every node/.style={rectangle,rounded corners, font=\sffamily},roundnode/.style={circle, draw=black!60, fill=green!5, very thick, minimum size=7mm},
squarednode/.style={rectangle, draw=black!60, fill=blue!5, very thick, minimum size=5mm},
]
  \node (A)             [activityStarts]              {$A$};
  \node (M)      [startstop1, below of=A, xshift=-2.0cm, yshift=-2cm]
                                                        {$M$};
  \node (FS)      [startstop1, below of=A, xshift=3cm, yshift=-2cm]{$F_{W,S}$};
  
    \node (FI)      [startstop1, below of=FS, xshift=0cm, yshift=-1cm]{$F_{W,E}$};
     \node (FEI)      [startstop1, below of=FI, xshift=0cm, yshift=-1cm]{$F_{W,I}$};
    
        \node (SS)      [startstop3, right of=FS, xshift=1cm, yshift=0cm]{$S_S$};
        \node (SI)      [startstop3, below of=SS, xshift=0cm, yshift=-1cm]{$S_E$};
         \node (SEI)      [startstop3, below of=SI, xshift=0cm, yshift=-1cm]{$S_I$};

\node (MT) [startstop3, right of=A, xshift=4cm, yshift=-0.5cm] {$M_S$};

\node (S) [startstop2, left of=A, xshift=-4cm, yshift=-0.5cm] {$S_h$};
\node (I) [startstop2, below of=S, xshift=0cm, yshift=-2.0cm] {$I_h$};
\node (R) [startstop2, below of=I, xshift=0cm, yshift=-1.5cm] {$R_h$};

  \draw[->] (A) -- (-3.5,0) node [midway, above] {$\mu_{A,1}+\mu_{A,2}A$};
 \draw[->] (FS) -- (5.0,-4.0) node [midway, below] {$\mu_{S}$};
 \draw[->] (FI) -- (5.0,-7.0) node [midway, below] {$\mu_{S}$};
 \draw[->] (FEI) -- (5.0,-10) node [midway, above] {$\mu_{I}$};
 \draw[->] (M) -- (-2,-2.6) node [midway, left] {$\mu_{M}$};
 \draw[->] (MT) -- (8.0,-.5) node [midway, above] {$\mu_{M_S}$};
 \draw[->] (3.0,-.5) -- (MT) node [midway, above] {$\left(1-\epsilon_{F}\right)\Lambda_{tot}$};
 \draw[->] (SS) -- (8.0,-4.0) node [midway, above] {$\mu_{S}$};
 \draw[->] (6.0,-2.0) -- (SS)  node [midway, right] {$\epsilon_F \Lambda_{tot}$};
 \draw[->] (SI) -- (8.0,-7) node [midway, above] {$\mu_{S}$};
 \draw[->] (SEI) -- (8.0,-10.0) node [midway, above] {$\mu_{I}$};
    \draw[->] (A) -- (M)  node [midway, left] {$(1-r)\gamma$};
   \draw[->] (A.east) -- (FS.north)  node [near end, left] {$r\gamma \frac{M+\varepsilon M_S}{M+M_S}$};  
    \draw[->,dashed] (MT.west) -- (2.0,-2.0);
    \draw[->] (A.east) -- (6.00,-3.25)  node [midway, right] {$r\gamma \frac{(1-\varepsilon) M_S}{M+M_S}$};  
    \draw[->,dashed] (MT) -- (5.25,-2.75);
    \draw[->,dashed] (M)[anchor=east] -- (4.75,-2.5);
    \draw[->,dashed] (M)[anchor=east] -- (1.75,-1.25);
    \draw[->] (FS.west)-- (0.,-4)node [midway,above]{$\phi$} -- (A);
    \draw[->] (FI.west)-- (0.,-7)node [midway,above]{$\phi$} -- (A);
    \draw[->] (FEI.west)-- (0.,-10)node [midway,above]{$\phi$} -- (A);
    \draw[->] (FS) -- (FI)node [midway,right] {$\frac{B\beta_{hm} I_h}{N_h}$};
     \draw[->] (FI) -- (FEI)node [midway,right] {$\nu_m$};
    \draw[->,dashed] (I) -- (FI.north) ;
 \draw[->] (SS) -- (SI)node [midway,right] {$\frac{B\beta_{hm} I_h}{N_h}$};
 \draw[->] (SI) -- (SEI)node [midway,right] {$\nu_m$};
    \draw[->,dashed] (I) -- (SI.north) ;   
     \draw[->] (S) -- (I) node [midway,right] {$\frac{B\beta_{mh} \left(F_I+S_I\right)}{N_h}$};
     \draw [->,dashed] (FEI.south) -- (-8.5,-10.75) -- (-8.5,-2.75) -- (-6.0,-2.75);
      \draw [->,dashed] (SEI.south) -- (6.0,-11) -- (-8.75,-11) -- (-8.75,-2.25) -- (-6.0,-2.25);
    \draw[->] (I) -- (R) node [midway, left] {$\nu_h$};
    \draw[->] (R) -- (-8.0,-8.0) node [midway, above] {$\mu_h$};
    \draw[->] (I) -- (-8.0,-4.5) node [midway, above] {$\mu_h$};
    \draw[->] (S) -- (-8.0,-0.5) node [midway, above] {$\mu_h$};
    \draw[->] (-4.0,-0.5) -- (S) node [midway, above] {$\mu_h N_h$};
    
  \end{tikzpicture}
\caption{Flow diagram of model \eqref{Human-EDO}-\eqref{Wild-Insect-EDO}.}
\label{flow_diagramme}
\end{figure}
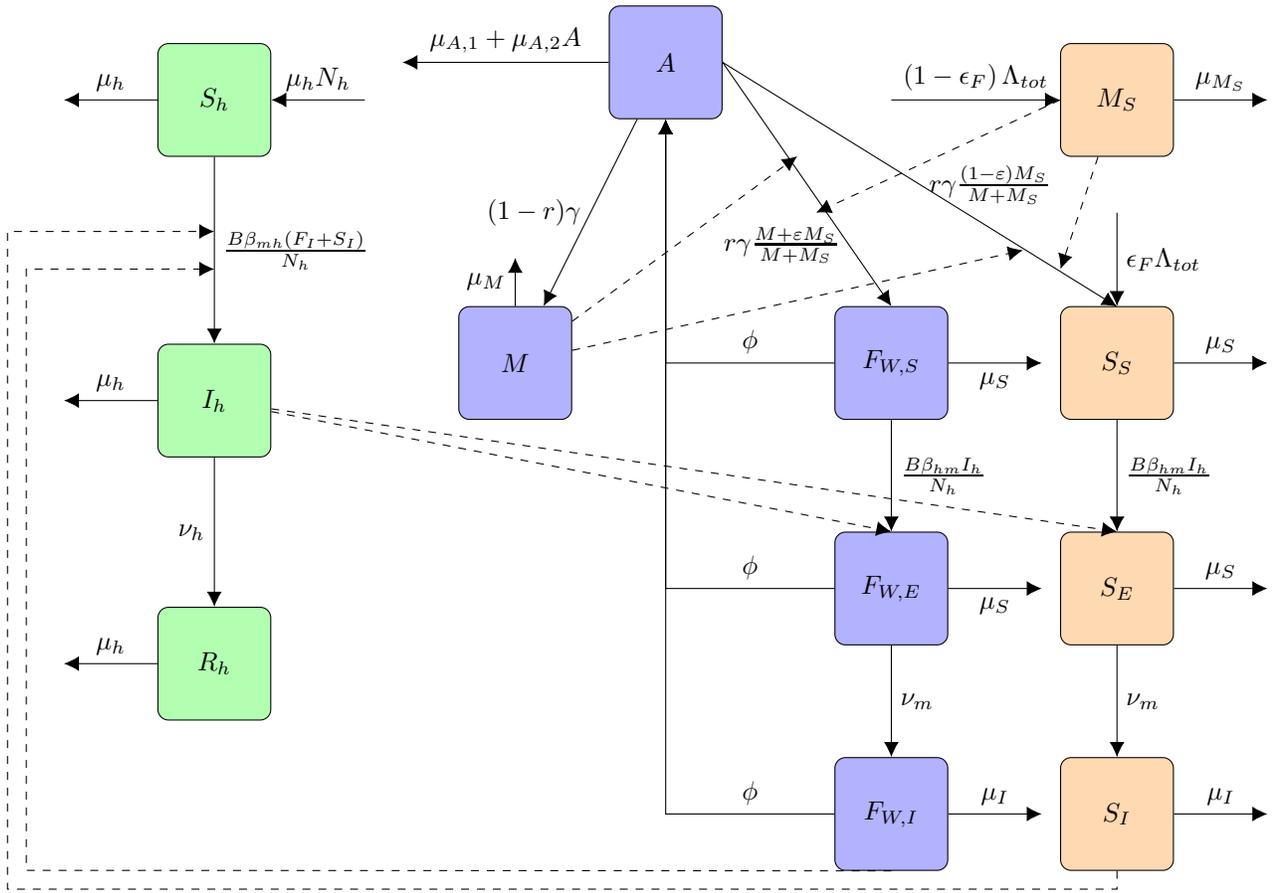

Thus, in order to take into account the release of sterile females and the effect of residual fertility, we have to consider a sub-populations of sterile females, $S$. Moreover, to take into account the circulation of a vector-borne virus, with an extrinsic incubation period of the virus within the vector population, we consider three epidemiological states, i.e. the susceptible, exposed and infected states, for the sterile  and the wild females, $S_S$, $S_E$, $S_I$, $F_{W,S}$, $F_{W,E}$, and $F_{W,I}$. We assume that the total population of humans, $N_h$, is positive and constant. It is also divided in three epidemiological states, i.e. $N_h=S_h+I_h+R_h$. When (wild and sterile) female mosquitoes are infected, we assume that their mortality rate can be impacted. 
Thus following \cite{DumontYatat2022}, and the flow diagram given in Fig. \ref{flow_diagramme}, page \pageref{flow_diagramme}, we derive the following SIT-entomological-epidemiological model
\begin{equation}\label{Human-EDO}
    \left\{ %
\begin{array}{lcl}
\dfrac{dS_{h}}{dt} & = & \mu_{h}N_{h}-B\beta_{mh}{\displaystyle \frac{F_{W,I}+S_{I}}{N_{h}}S_{h}-\mu_{h}S_{h},}\\
{\displaystyle \frac{dI_{h}}{dt}} & = & B\beta_{mh}{\displaystyle \frac{F_{W,I}+S_{I}}{N_{h}}S_{h}-\nu_hI_{h}-\mu_{h}I_{h},}\\
\dfrac{dR_{h}}{dt} & = & \nu_hI_{h}-\mu_{h}R_{h},
\end{array}\right.
\end{equation}
\begin{equation}\label{Wild-Insect-EDO}
\left\{ %
\begin{array}{lcl}
{\displaystyle \frac{dA}{dt}} & = & \phi (F_{W,S}+F_{W,E}+F_{W,I})-(\gamma+\mu_{A,1}+\mu_{A,2}A)A,\\
{\displaystyle \frac{dM}{dt}} & = & (1-r)\gamma A-\mu_{M}M,\\
{\displaystyle \frac{dF_{W,S}}{dt}} & = & \dfrac{M+\varepsilon M_S}{M+M_{S}}r\gamma A-B\beta_{hm}{\displaystyle \frac{I_{h}}{N_{h}}F_{W,S}-\mu_{S}F_{W,S},}\\
{\displaystyle \frac{dF_{W,E}}{dt}} & = & B\beta_{hm}{\displaystyle \frac{I_{h}}{N_{h}}F_{W,S}-(\nu_{m}+\mu_{S})F_{W,E},}\\
{\displaystyle \frac{dF_{W,I}}{dt}} & = & \nu_{m}F_{W,E}-\mu_{I}F_{W,I},\\
{\displaystyle \frac{dS_{S}}{dt}} & = & \epsilon_F \Lambda_{tot}+\dfrac{(1-\varepsilon)M_{S}}{M+M_{S}}r\gamma A-B\beta_{hm}{\displaystyle \frac{I_{h}}{N_{h}}S_{S}-\mu_{S}S_{S},}\\
{\displaystyle \frac{dS_{E}}{dt}} & = & B\beta_{hm}{\displaystyle \frac{I_{h}}{N_{h}}S_{S}-(\nu_{m}+\mu_{S})S_{E},}\\
{\displaystyle \frac{dS_{I}}{dt}} & = & \nu_{m}S_{E}-\mu_{I}S_{I},\\
\dfrac{dM_{S}}{dt} & = & (1-\epsilon_F) \Lambda_{tot}-\mu_{M_{S}}M_{S},
\end{array}\right.
\end{equation}
with appropriate non-negative initial conditions.

\comment{
\begin{equation}\label{Released-Insect-EDO}
\left\{ \begin{array}{l}
\dfrac{dF_{S}}{dt}=\epsilon_F \Lambda_{tot}-B\beta_{hm}\dfrac{I_{h}}{N_{h}}F_{S}-\mu_{F_S}F_{S}\\
{\displaystyle \frac{dF_{E}}{dt}}=B\beta_{hm}\dfrac{I_{h}}{N_{h}}F_{S}{\displaystyle -(\nu_{m}+\mu_{F_S})F_{E},}\\
{\displaystyle \frac{dF_{I}}{dt}}=\nu_{m}F_{E}-\mu_{F_I}F_{I},\\
\dfrac{dM_{S}}{dt}=\Lambda_{M}-\mu_{M_{S}}M_{S}.
\end{array}\right.
\end{equation}
}
We summarize all the model parameters in Table \ref{table:parameters}, page \pageref{table:parameters}. In \cite{Duprez2022}, the authors have considered varying parameters to take into account variations of temperature and raining along the year in R\'eunion island and their impact on SIT strategies to reduce the nuisance or the epidemiological risk. Thus, in Table \ref{table:parameters}, page \pageref{table:parameters}, we derive the variations for each parameters from a daily average temperature varying between $15^\circ$ and $30^\circ$. These interval values will be used for a global sensitivity analysis done in section \ref{section4}. In the simulations part, we will consider parameter values related to an average temperature of $25^\circ$, that is (close to) the most favorable temperature for \textit{Aedes albopictus} mosquito dynamics.

\begin{table}[h]
  \centering
 \begin{tabular}{|p{1.7cm}|p{5.0cm}|p{1.7cm}|p{2.2cm}|p{2.5cm}|p{1.8cm}|}
   \hline
   Parameters & Description & Unit & Range & Baseline for simulation ($T=25^\circ$) & Reference\\ \hline \hline
   $1/\mu_{h}$& Average human lifespan & Day  & $[60,80]\times365$ & $78\times365$ &  \\
   $1/\nu_h$& Average DENV viremic period & Day  & $[1,7]$ & $7$ & \cite{Vaughn2000} \\
   $B$& daily number of mosquito bites on human & - &  $[0.1,1]$ & $0.25$ &  \\
   $\beta_{mh}$ & Rate of transmission of DENV from Infected mosquito to Susceptible human& Day$^{-1}$ &  $[0.12;0.57]$ & $0.3427$ & \cite{Duprez2022}\\
   $\beta_{hm}$ &  Rate of transmission of DENV from Infected human to Susceptible mosquito &Day$^{-1}$  & $[0.4;0.96]$ & $0.872$  & \cite{Duprez2022}\\
   $\mu_{A,1}$ & Natural death rate for larvae and pupae stage. &Day$^{-1}$   & $[0.019;0.299]$ &  $0.0262$  & \cite{Duprez2022} \\
   $\mu_{A,2}$ & Density-induced death rate for larvae and pupae stage. &Day$^{-1}$Ind$^{-1}$  & $[2\times 10^{-5};0.02]$ & $1.76 \times 10^{-4}$ & \cite{DumontYatat2022,Duprez2022} \\
   $\phi$ &  Daily hatching eggs deposit  & Day$^{-1}$  &$[0,11]$ &  $10$  & \cite{Duprez2022}\\
   $\gamma$ &  Transition rate from non-adult stage to adult-stage. & Day$^{-1}$  & $[0.028,0.12]$& $0.0962$ & \cite{Duprez2022} \\
   $r$ &  Sex-ratio & -  & $[0.4,0.6]$ & $0.5$& \\
   $\mu_{S}$& Female mosquito death rate & Day$^{-1}$  & $[0.035,0.07]$ & $0.0453$ &  \cite{Duprez2022} \\
   $\mu_I$ & Infected female mosquito death rate & Day$^{-1}$  &  $[0.035,0.07]$ & $0.0453$ & \cite{Duprez2022}\\
    $\mu_M$ & Male mosquito death rate & Day$^{-1}$  &  $[0.05,0.082]$ & $0.0722$ & \cite{Duprez2022}\\
    $\mu_{M_S}$ & Sterile Male mosquito death rate & Day$^{-1}$  &  $[0.1,0.2]$ & $0.1$ & \cite{Duprez2022}\\
    $\nu_m$ & Extrinsic incubation rate & Day$^{-1}$ & $[0.015,0.25]$ & $0.184$ & \cite{Duprez2022} \\
   $\Lambda_{tot}$ & Sterile insect release rate & Ind Day$^{-1}$& $[0;18 000]$ & varying &  \\
   $\varepsilon$ & Residual fertility & - & $[0;0.05]$ & varying &  \\
   $\epsilon_F$ & Sterile female contamination & - & $[0;0.05]$ & varying &  \\
   \hline
 \end{tabular}
  \caption{Parameters description and parameters values for the entomological-epidemiological model related to Dengue circulation, for an average temperature of $T=25^\circ$C and $N_h=20000$. }
 \label{table:parameters}
 \end{table}

\subsection{The wild insect model without SIT}
We deduce from system \eqref{Human-EDO}-\eqref{Wild-Insect-EDO} that dynamics of wild insects, without SIT, is modelled by system \eqref{Wild-Insect-ode-withoutSIT}:
\begin{equation}\label{Wild-Insect-ode-withoutSIT}
\left\{%
\begin{array}{lcl}
\displaystyle\frac{dA}{d t} &=& \phi F_{W,S}-(\gamma+\mu_{A,1}+\mu_{A,2}A)A,\\
\displaystyle\frac{dM}{d t} &=& (1-r)\gamma A-\mu_MM,\\
\displaystyle\frac{dF_{W,S}}{d t} &=& r\gamma A-\mu_SF_{W,S}.\\
\end{array}
\right.
\end{equation}
System \eqref{Wild-Insect-ode-withoutSIT} is quite simple and assumes implicitly that there are always adults of both sex (male and female), such that emerging females will always mate with a male and thus become fertile. In addition, system \eqref{Wild-Insect-ode-withoutSIT} has been considered and studied in previous works, see e.g. \cite{Anguelov2020,DumontYatat2022}. Hence, below we recall its main qualitative results without any proofs.

The basic offspring number related to model (\ref{Wild-Insect-ode-withoutSIT}) is 
\begin{equation}
\N=\displaystyle\frac{r\gamma\phi}{\mu_S(\gamma+\mu_{A,1})}.
\label{BO}
\end{equation}

Setting the right-hand side of system (\ref{Wild-Insect-ode-withoutSIT}) to zero we obtain the extinction equilibrium $\textbf{0}_{\mathbb{R}^3}=(0,0,0)^T$ and the equilibrium $E^*=(A^*,M^*,F_{W,S}^*)^T$ given by 
\begin{equation}\label{Wild-equilibria-definition}
 \left\{%
\begin{array}{rcl}
A^*&=& \displaystyle\frac{(\gamma+\mu_{A,1})}{\mu_{A,2}}(\N-1),\\
 M^*&=& \displaystyle\frac{(1-r)\gamma A^*}{\mu_M}=\displaystyle\frac{(1-r)\gamma}{\mu_M}\displaystyle\frac{(\gamma+\mu_{A,1})}{\mu_{A,2}}(\N-1),\\
 F_{W,S}^*&=& \displaystyle\frac{r\gamma A^*}{\mu_S}=\displaystyle\frac{r\gamma}{\mu_S}\displaystyle\frac{(\gamma+\mu_{A,1})}{\mu_{A,2}}(\N-1).\\
\end{array}
\right.
\end{equation}
The inequalities between vectors are considered here in their usual coordinate-wise sense. Clearly, $E^*>\bf{0}_{\mathbb{R}^3}$ if and only if $\N>1$. 
We summarize these results with some more details related to basins of attraction of equilibria in the following theorem. 

\begin{theorem}[\cite{Anguelov2020,DumontYatat2022}]
\label{Mosquitoes-ode-theorem}
Model (\ref{Wild-Insect-ode-withoutSIT}) defines a forward dynamical system on 
$\mathcal{D}=\{x\in\mathbb{R}^3:x\geq\bf{0}_{\mathbb{R}^3}\}$. Furthermore, 
\begin{enumerate}
    \item[1)] If $\N\leq1$ then $\bf{0}_{\mathbb{R}^3}$ is globally asymptotically stable on
    $\mathcal{D}$.
    \item[2)] If $\N>1$ then $E^*$ is stable with
    basin of attraction $$\mathcal{D}\setminus\{x=(A,M,F_{W,S})^T\in\mathbb{R}^3_+:A=F_{W,S}=0\},$$  and $\bf{0}_{\mathbb{R}^3}$ is unstable with the non negative $M-$axis being a stable manifold.
\end{enumerate}
\end{theorem}
\begin{proof} See \cite[Theorem 1]{Anguelov2020,DumontYatat2022}.
\end{proof}

\begin{remark}
    Mechanical control, that is the removing of mosquito breeding sites, has an impact on $\mu_{A,2}$ because it depends on $K$, the larvae-carrying capacity that is defined by $K=3 \times N_h$ \cite{DumontYatat2022}[section 7]
    \begin{equation}
        \mu_{A,2}=\dfrac{\gamma+\mu_{A,1}}{K} \N.
    \end{equation}
Thus reducing $K$ by a certain percentage, say $p_{mc}$, will increase $\mu_{A,2}$ by a factor $\dfrac{1}{1-p_{mc}}$.
\end{remark}

\subsection{The wild insect model with SIT}
\label{WimwSIT}
We now consider the following SIT-entomological model that occurs
when no virus is circulating. Its study is helpful to derive the Disease
Free Equilibrium, DFE, thanks to several release sizes. Thanks to
the fact that $t$ is sufficiently large or that the initial releases
are such that $M_{S}(0)=M_{S}^{*}=(1-\epsilon_F)\dfrac{\Lambda_{tot}}{\mu_{M_{S}}}$. The entomological model assumes the form

\begin{equation}\label{ODE-entomo-old}
    \left\{ \begin{array}{l}
{\displaystyle \frac{dA}{dt}}=\phi F_{W,S}-(\gamma+\mu_{A,1}+\mu_{A,2}A)A,\\
{\displaystyle \frac{dM}{dt}}=(1-r)\gamma A-\mu_{M}M,\\
{\displaystyle \frac{dF_{W,S}}{dt}}={\displaystyle \frac{M+\varepsilon M_S^*}{M+M_{S}^{*}}r\gamma A-\mu_{S}F_{W,S},}\\
{\displaystyle \frac{dS_{S}}{dt}}={\displaystyle \epsilon_F \Lambda_{tot}+ \frac{(1-\varepsilon)M_{S}^{*}}{M+M_{S}^{*}}r\gamma A-\mu_{S}S_{S}.}
\end{array}\right.
\end{equation}

Since the released sterile females do not play a role in the wild mosquito dynamics, we derive the following reduced SIT-entomological model

\begin{equation}\label{ODE-entomo}
    \left\{ \begin{array}{l}
{\displaystyle \frac{dA}{dt}}=\phi F_{W,S}-(\gamma+\mu_{A,1}+\mu_{A,2}A)A,\\
{\displaystyle \frac{dM}{dt}}=(1-r)\gamma A-\mu_{M}M,\\
{\displaystyle \frac{dF_{W,S}}{dt}}={\displaystyle \frac{M+\varepsilon M_S^*}{M+M_{S}^{*}}r\gamma A-\mu_{S}F_{W,S}.}\\
\end{array}\right.
\end{equation}
We now deal with equilibria of model \eqref{ODE-entomo}. Of course, given an equilibrium $\Bar{E}=(\Bar{A}, \Bar{M}, \Bar{F}_{W,S})^T$ of system \eqref{ODE-entomo}, we can recover the $S_S$-component of the corresponding equilibrium of system \eqref{ODE-entomo-old}, by setting $$\Bar{S}_S=\dfrac{1}{\mu_S}\left(\epsilon_F \Lambda_{tot}+ \frac{(1-\varepsilon)M_{S}^{*}}{\Bar{M}+M_{S}^{*}}r\gamma \Bar{A}\right).$$
We follow the methodology developed in \cite{Anguelov2020}. When $A=0$, we obtain the elimination equilibrium $E_0=(0,0,0)^T$.
Assuming $A\neq0$, then from the first equation, we
derive
\begin{equation}\label{equation-M}
    \dfrac{\phi r\gamma}{\mu_{S}}\frac{M+\varepsilon M_{S}^{*}}{M+M_{S}^{*}}=(\gamma+\mu_{A,1}+\mu_{A,2}A).
\end{equation}

Then, using the fact that
\[
A=\dfrac{\mu_{M}}{(1-r)\gamma}M,
\]
 setting
\[
\mathcal{Q}=\dfrac{\mu_{A,2}\mu_{M}}{(\gamma+\mu_{A,1})(1-r)\gamma},
\]
and
$$\alpha=\dfrac{M_S^{*}}{M},$$
we derive
\begin{equation}\label{aux1}
\dfrac{1+\varepsilon\alpha}{1+\alpha}\mathcal{N}=1+\dfrac{\mathcal{Q} M_{S}^{*}}{\alpha}.
\end{equation}

Setting $\Q_S=  M_{S}^{*}\Q>0$, equation \eqref{aux1} becomes
\begin{equation}\label{aux2}
	\begin{array}{lcl}
		\left(1-\mathcal{N}\varepsilon\right)\alpha^2+\left(1+\mathcal{Q}_S-\mathcal{N}\right)\alpha+\mathcal{Q}_S=0.
	\end{array}
\end{equation}
The discriminant of \eqref{aux2} is 
\begin{equation}\label{aux3}
	\begin{array}{lcl}
		\Delta(\mathcal{Q}_S)=(\mathcal{Q}_S)^2 + \mathcal{Q}_S\left(4\mathcal{N}\varepsilon-2\left(\mathcal{N}+1\right)\right)+\left(\mathcal{N}-1\right)^2.
	\end{array}
\end{equation}
To study the sign of $\Delta(\mathcal{Q}_S)$, we consider the sub-determinant of $\Delta$
\begin{equation}\label{aux4}
	\begin{array}{lcl}
		\Delta'=16\left(1-\mathcal{N}\varepsilon\right)\left(1-\varepsilon\right)\mathcal{N}.
	\end{array}
\end{equation}
Since $1-\varepsilon \geq 0$, $\Delta'$ has the same sign as $1-\mathcal{N}\varepsilon$.
\begin{enumerate}
	\item Assume that $\mathcal{N}\varepsilon <1$. Then, $\Delta' > 0$ and $\Delta$ has two real roots $\mathcal{Q}_{S_1}$ and $\mathcal{Q}_{S_2}$ such that:
	\begin{equation}\label{aux5}
		\left\{
		\begin{array}{lcl}
			\mathcal{Q}_{S_1}\mathcal{Q}_{S_2} = \left(1-\mathcal{N}\right)^2 > 0,\\
			\mathcal{Q}_{S_1} + \mathcal{Q}_{S_2} =  2\left(1-\mathcal{N}\varepsilon+\mathcal{N}\left(1-\varepsilon\right)\right) > 0,\\
			\mathcal{Q}_{S_1} = \left(\sqrt{\mathcal{N}\left(1-\varepsilon\right)}-\sqrt{1-\mathcal{N}\varepsilon}\right)^{2}>0,\\
			\mathcal{Q}_{S_2} = \left(\sqrt{\mathcal{N}\left(1-\varepsilon\right)}+\sqrt{1-\mathcal{N}\varepsilon}\right)^{2} > \mathcal{Q}_{S_1}.
		\end{array}
		\right.
	\end{equation}
	It therefore follows that $\Delta(\mathcal{Q}_S)\geq0$ when $\mathcal{Q}_S \in \left(0,\mathcal{Q}_{S_1}\right] \cup \left[\mathcal{Q}_{S_2}, +\infty\right)$ and $\Delta(\mathcal{Q}_S) < 0$ when $\mathcal{Q}_S \in \left(\mathcal{Q}_{S_1}, \mathcal{Q}_{S_2}\right)$. The following discussion is valid:
	\begin{itemize}
		\item Assume that $\mathcal{Q}_S \in \left(0,\mathcal{Q}_{S_1}\right).$ Then, \eqref{aux2} admits two real roots $\alpha_{-}$, $\alpha_{+}$ where 
		\begin{equation}\label{alpha-1}
			\begin{array}{lcl}
				\alpha_{\pm} = \dfrac{\left(\mathcal{N}-\mathcal{Q}_S-1\right) \pm \sqrt{\Delta(\mathcal{Q}_S)}}{2\left(1-\mathcal{N}\varepsilon\right)}.
			\end{array}
		\end{equation}
		Note that 
		\begin{equation*}
			\mathcal{N}-1-\mathcal{Q}_S > \mathcal{N}-1-\mathcal{Q}_{S_1} = 2\left(\sqrt{\left(1-\mathcal{N}\varepsilon\right)\left(1-\varepsilon\right)\mathcal{N}}-\left(1-\mathcal{N}\varepsilon\right)\right) > 0.
		\end{equation*}
		Since\\
		$\alpha_{-}\alpha_{+} = \dfrac{\mathcal{Q}_S}{1-\mathcal{N}\varepsilon} > 0$,  $\mathcal{N}-1-\mathcal{Q}_S > 0$ and $\alpha_{+}+\alpha_{-} = \dfrac{\mathcal{N}-1-\mathcal{Q}_S}{1-\mathcal{N}\varepsilon} > 0$, we deduce that $0 < \alpha_{-} < \alpha_{+}.$
		\item Assume that $\mathcal{Q}_S \in \left(\mathcal{Q}_{S_2}, +\infty\right).$ Then, \eqref{aux2} admits two real roots $\alpha_{-}, \alpha_{+}$ where
		\begin{equation}\label{alpha-2}
			\begin{array}{lcl}
				\alpha_{\pm} = \dfrac{\left(\mathcal{N}-\mathcal{Q}_S-1\right) \pm \sqrt{\Delta(\mathcal{Q}_S)}}{2\left(1-\mathcal{N}\varepsilon\right)}.
			\end{array}
		\end{equation}
		Note that
		\begin{equation*}
			\mathcal{N}-1-\mathcal{Q}_S < \mathcal{N}-1-\mathcal{Q}_{S_2} = -2\left(\sqrt{\left(1-\mathcal{N}\varepsilon\right)\left(1-\varepsilon\right)\mathcal{N}}+\left(1-\mathcal{N}\varepsilon\right)\right) < 0.
		\end{equation*}
		Since\\
		$\alpha_{-}\alpha_{+} = \dfrac{\mathcal{Q}_S}{1-\mathcal{N}\varepsilon} > 0$,  $\mathcal{N}-1-\mathcal{Q}_S < 0$ and $\alpha_{+}+\alpha_{-} = \dfrac{\mathcal{N}-1-\mathcal{Q}_S}{1-\mathcal{N}\varepsilon} < 0$, we deduce that $\alpha_{-} < \alpha_{+} < 0.$
		\item Assume that $\mathcal{Q}_S \in \left(\mathcal{Q}_{S_1}, \mathcal{Q}_{S_2}\right).$ Then, \eqref{aux2} does not admit real roots.
		\item Assume that $\mathcal{Q}_S = \mathcal{Q}_{S_1}.$ Then, \eqref{aux2} has only one real solution
		\begin{equation}\label{alpha-3}
			\begin{array}{lcl}
				\alpha_{\diamond} = \dfrac{\mathcal{N}-1-\mathcal{Q}_{S_1}}{2\left(1-\mathcal{N}\varepsilon\right)} > 0.
			\end{array}
		\end{equation}
		\item Assume that $\mathcal{Q}_S = \mathcal{Q}_{S_2}.$ Then, \eqref{aux2} has only one real solution $$\alpha_{-} = \alpha_{+} = \dfrac{\mathcal{N}-1-\mathcal{Q}_{S_2}}{2\left(1-\mathcal{N}\varepsilon\right)} < 0.$$
	\end{itemize}
	\item Assume that $\mathcal{N}\varepsilon >1.$ Then $\Delta' < 0$ and $\Delta(\mathcal{Q}_S) > 0$. Therefore, \eqref{aux2} admits two real roots $\alpha_{-}$, $\alpha_{+}$.
	Since $\alpha_{-}\alpha_{+} = \dfrac{\mathcal{Q}_S}{1-\mathcal{N}\varepsilon} < 0$. It follows that
	\begin{equation}\label{alpha-4}
		\begin{array}{lcl}
			\alpha_{-} < 0 < \alpha_{+} = \dfrac{-\left(\mathcal{N}-1-\mathcal{Q}_S\right) + \sqrt{\Delta(\mathcal{Q}_S)}}{2\left(\mathcal{N}\varepsilon-1\right)}.
		\end{array}
	\end{equation}
			
  
	\item Assume that $\mathcal{N}\varepsilon = 1$. Then,  
	\eqref{aux2} admits a unique solution
	\begin{equation}\label{alpha-5}
		\begin{array}{lcl}
			\alpha_{\sharp}=
			\dfrac{\mathcal{Q}_S}{\mathcal{N}-1-\mathcal{Q}_S}.
		\end{array}
	\end{equation}
	$\alpha_{\sharp} > 0$ whenever $\mathcal{Q}_S < \mathcal{N}-1.$
\end{enumerate}
From the previous discussion, we deduce, for $M_S^*=\dfrac{\Lambda_M}{\mu_{M_S}}$, the following:

\begin{theorem}\label{theo-equilibre-sit-only}  System \eqref{ODE-entomo-old} always admits the trivial equilibrium $E_{0} = \left(0, 0, 0,\dfrac{\epsilon_F \Lambda_{tot}}{\mu_S}\right)^T$. In addition:
	\begin{enumerate}
		\item Assume that $ \mathcal{N}\varepsilon < 1.$ Consider the threshold 
		\begin{equation}\label{seuil-ms}
			\begin{array}{lcl}
				\Lambda_{M}^{crit} = \dfrac{\mu_{M_S}}{\mathcal{Q}}\left(\sqrt{\mathcal{N}\left(1-\varepsilon\right)}-\sqrt{1-\mathcal{N}\varepsilon}\right)^{2}.
			\end{array}
		\end{equation}
		\begin{enumerate}
			\item If $(1-\epsilon_F) \Lambda_{tot} \in \left(0, \Lambda_{M}^{crit}\right)$, then system \eqref{ODE-entomo-old} admits two positive equilibria $E_{1} = \left(A_{1}, M_{1}, F_{W,S_1},S_{S_1}\right)^T$ and $E_{2} = \left(A_{2}, M_{2}, F_{W,S_2},S_{S_2}\right)^T$, such that $(A_{1}, M_{1}, F_{W,S_1})^T < (A_{2}, M_{2}, F_{W,S_2})^T$ and
			\begin{equation*}
				\left\{
				\begin{array}{lcl}
					M_{1} = \dfrac{(1-\epsilon_F) \Lambda_{tot}}{\mu_{M_S}\alpha_{+}},\quad\mbox{where $\alpha_{+}$ is computed from \eqref{alpha-1}},\\
					M_{2} = \dfrac{(1-\epsilon_F) \Lambda_{tot}}{\mu_{M_S}\alpha_{-}},\quad\mbox{where $\alpha_{-}$ is computed from \eqref{alpha-1}},\\
					A_{1,2} = \dfrac{\mu_M}{\left(1-r\right)\gamma}M_{1,2},\\
					F_{W,S_{1,2}} = \dfrac{\left(\gamma+\mu_{1}+\mu_{2}A_{1,2}\right)A_{1,2}}{\phi}, \\
					S_{S_{1,2}} = \displaystyle\frac{1}{\mu_S}\left(\epsilon_F \Lambda_{tot}+\dfrac{(1-\varepsilon)M_S^*}{M_{1,2}+M_S^*} r\gamma A_{1,2}\right).
				\end{array}
				\right.
			\end{equation*}
			\item If $(1-\epsilon_F) \Lambda_{tot}=\Lambda_{M}^{crit}$, then system \eqref{ODE-entomo-old} admits a unique equilibrium $E_{\diamond} = \left(A_{\diamond}, M_{\diamond}, F_{W,S_\diamond},S_{S_{\diamond}}\right)^T$ where
			\begin{equation*}
				\left\{
				\begin{array}{lcl}
					M_{\diamond} = \dfrac{\Lambda_M}{\mu_{M_S}\alpha_{\diamond}}, \quad\mbox{where $\alpha_{\diamond}$ is computed from \eqref{alpha-3}},\\
					A_{\diamond} = \dfrac{\mu_M}{\left(1-r\right)\gamma}M_{\diamond},\\
					F_{\diamond} = \dfrac{\left(\gamma+\mu_{1}+\mu_{2}A_{\diamond}\right)A_{\diamond}}{\phi}, \\
					S_{S_{\diamond}} = \displaystyle\frac{1}{\mu_S}\left(\epsilon_F \Lambda_{tot}+\dfrac{(1-\varepsilon)M_S^*}{M_{\diamond}+M_S^*} r\gamma A_{\diamond}\right).
				\end{array}
				\right.
			\end{equation*}
		\end{enumerate}
		\item Assume that $\mathcal{N}\varepsilon > 1.$ Then, for any $(1-\epsilon_F) \Lambda_{tot} > 0$, system \eqref{ODE-entomo-old} admits a unique positive equilibrium $E_{\dag} = \left(A_{\dag}, M_{\dag}, F_{W,S_\dag},S_{S_\dag}\right)^T$ where
		\begin{equation*}
			\left\{
			\begin{array}{lcl}
				M_{\dag} = \dfrac{(1-\epsilon_F) \Lambda_{tot}}{\mu_{M_S}\alpha_{+}}, \quad\mbox{ where $\alpha_{+}$ is computed from \eqref{alpha-4}},\\
				A_{\dag} = \dfrac{\mu_M}{\left(1-r\right)\gamma}M_{\dag},\\
				F_{W,S_\dag} = \dfrac{\left(\gamma+\mu_{1}+\mu_{2}A_{\dag}\right)A_{\dag}}{\phi}, \\
				S_{S_{\dag}} = \displaystyle\frac{1}{\mu_S}\left(\epsilon_F \Lambda_{tot}+\dfrac{(1-\varepsilon)M_S^*}{M_{\dag}+M_S^*} r\gamma A_{\dag}\right).
			\end{array}
			\right.
		\end{equation*}
		\item Assume that $\mathcal{N}\varepsilon = 1.$ Consider the threshold
		\begin{equation*}
			\Lambda_{M, \sharp}^{crit} = \Lambda_{M}^{crit}\vert_{\mathcal{N}\varepsilon=1} = \dfrac{\mu_{M_S}}{ \mathcal{Q}}(\mathcal{N}-1) > 0.
		\end{equation*}
		If $(1-\epsilon_F) \Lambda_{tot} \in \left(0, \Lambda_{M, \sharp}^{crit}\right)$, then system \eqref{ODE-entomo-old} admits a unique positive equilibrium $E_{\sharp} = \left(A_{\sharp}, M_{\sharp}, F_{W,S_\sharp},S_{S_\sharp}\right)^T$ where
		\begin{equation*}
			\left\{
			\begin{array}{lcl}
				M_{\sharp} = \dfrac{(1-\epsilon_F) \Lambda_{tot}}{\mu_{M_S}\alpha_{\sharp}}, \quad\mbox{ where $\alpha_{\sharp}$ is computed from \eqref{alpha-5}},\\
				A_{\sharp} = \dfrac{\mu_M}{\left(1-r\right)\gamma}M_{\sharp},\\
				F_{W,S_\sharp} = \dfrac{\left(\gamma+\mu_{1}+\mu_{2}A_{\sharp}\right)A_{\sharp}}{\phi}, \\
				S_{S_{\sharp}} = \displaystyle\frac{1}{\mu_S}\left(\epsilon_F \Lambda_{tot}+\dfrac{(1-\varepsilon)M_S^*}{M_{\sharp}+M_S^*} r\gamma A_{\sharp}\right).
			\end{array}
			\right.
		\end{equation*}
	\end{enumerate}
\end{theorem}	

\begin{remark}
    When $\varepsilon=0$, we recover the critical rate $\Lambda_{M}^{crit}$ defined in \cite{Anguelov2020,DumontYatat2022}. 
\end{remark}

Taking into account the fact that system \eqref{ODE-entomo} is cooperative, we are able to study stability properties of its equilibria and then to deduce the stability properties for system \eqref{ODE-entomo-old}. Thus, following \cite{Anguelov2012TIS,Anguelov2020,DumontYatat2022}, we obtain Theorem \ref{theo-stabilite-equilibre-sit} where $x=(A,M,F_{W,S},S)^T$.
\begin{theorem}\label{theo-stabilite-equilibre-sit}
	The following results are valid for system \eqref{ODE-entomo-old}:
	\begin{enumerate}
		\item Assume that $ \mathcal{N}\varepsilon <1$.
		\begin{enumerate}
			\item If $(1-\epsilon_F) \Lambda_{tot} > \Lambda_M^{crit},$ then $E_{0}$ is globally asymptotically stable.
			\item If $(1-\epsilon_F) \Lambda_{tot} \in \left(0, \Lambda_M^{crit}\right)$, then $E_1$ is unstable, and the set $\{x \in \mathbb{R}^4 : (0,0,0)^T \leq (A,M,F_{W,S})^T < (A_1,M_1,F_{W,S_1})^T\}$ is in the basin of attraction of $E_{0}$ and the set $\{x \in \mathbb{R}^4 : (A_1,M_1,F_{W,S_1})^T < (A,M,F_{W,S})^T \}$ is in the basin of attraction of $E_{2}.$
			\item If $(1-\epsilon_F) \Lambda_{tot} = \Lambda_M^{crit},$ then the set $\{x \in \mathbb{R}^4 : (0,0,0)^T \leq (A,M,F_{W,S})^T < (A_\diamond,M_\diamond,F_{W,S_\diamond})^T\}$ is in the basin of attraction of $E_{0}$, while the set $\{x \in \mathbb{R}^4 : (A_\diamond,M_\diamond,F_{W,S_\diamond})^T \leq (A,M,F_{W,S})^T \}$ is in the basin of attraction of $E_{\diamond}.$
		\end{enumerate}
		\item Assume that $\mathcal{N}\varepsilon > 1.$ Then, the elimination equilibrium $E_{0}$ is unstable and the coexistence equilibrium $E_{\dag}$ is globally asymptotically stable for any $(1-\epsilon_F) \Lambda_{tot}>0.$
		\item Assume that $\mathcal{N}\varepsilon = 1.$
		\begin{enumerate}
			\item If $(1-\epsilon_F) \Lambda_{tot} \geq \Lambda_{M, \sharp}^{crit},$ then $E_{0}$ is globally asymptotically stable.
			\item If $(1-\epsilon_F) \Lambda_{tot} \in \left(0, \Lambda_{M, \sharp}^{crit}\right)$, then the elimination equilibrium $E_{0}$ is unstable and the coexistence equilibrium $E_{\sharp}$ is globally asymptotically stable.
		\end{enumerate}
	\end{enumerate}
\end{theorem}

\begin{proof}
See Appendix \ref{AppendixB}.    
\end{proof}

\section{Qualitative analysis of the full SIT epidemiological model \label{section3}}
Now we turn to the more complex model described in the introduction. In the sequel, we assume that $\N>1$. Indeed, in the case where $\N\leq1$, by a comparison argument, the system will always converge toward the trivial disease-free equilibrium. 

Without SIT, this model has been studied in \cite{DumontYatat2022} where we derived the Basic Reproduction Number defined as follows
\begin{equation}\label{R0}
\R_0^2=\dfrac{\nu_m}{\nu_m+\mu_S}\dfrac{B\beta_{mh}}{\mu_I}\dfrac{B\beta_{hm}}{\nu_h +\mu_h}\dfrac{F_{W,S}^*}{N_h}.
\end{equation}
We assume that, without any control, $$\R_0^2>1.$$ 
From \cite{DumontYatat2022}, there exists a unique endemic equilibrium $$EE=(S_h^\sharp,I_h^\sharp,R_h^\sharp,A^\sharp,M^\sharp,F_{W,S}^\sharp,F_E^\sharp,F_I^\sharp)^T$$ when $\R_0^2>1$.

We will now proceed like in \cite[section 5]{DumontYatat2022}. 
In this section, we consider that constant and permanent SIT releases are done as a control tool. Hence, following \eqref{ODE-entomo}, the dynamics of human and mosquito populations  are described by system \eqref{Human-ode-good}-(\ref{Mosquitoes-ode-SIT-constant}):

\begin{equation}\label{Human-ode-good}
\left\{ %
\begin{array}{lcl}
\dfrac{dS_{h}}{dt} & = & \mu_{h}N_{h}-B\beta_{mh}{\displaystyle \frac{F_{W,I}+S_{I}}{N_{h}}S_{h}-\mu_{h}S_{h},}\\
{\displaystyle \frac{dI_{h}}{dt}} & = & B\beta_{mh}{\displaystyle \frac{F_{W,I}+S_{I}}{N_{h}}S_{h}-\nu_hI_{h}-\mu_{h}I_{h},}\\
\dfrac{dR_{h}}{dt} & = & \nu_hI_{h}-\mu_{h}R_{h},
\end{array}\right.
\end{equation}

\begin{equation}\label{Mosquitoes-ode-SIT-constant}
\left\{ %
\begin{array}{lcl}
{\displaystyle \frac{dA}{dt}} & = & \phi(F_{W,S}+F_{W,E}+F_{W,I})-(\gamma+\mu_{A,1}+\mu_{A,2}A)A,\\
{\displaystyle \frac{dM}{dt}} & = & (1-r)\gamma A-\mu_{M}M,\\
{\displaystyle \frac{dF_{W,S}}{dt}} & = & \dfrac{M+\varepsilon M_S^*}{M+M_{S}^*}r\gamma A-B\beta_{hm}{\displaystyle \frac{I_{h}}{N_{h}}F_{W,S}-\mu_{S}F_{W,S},}\\
{\displaystyle \frac{dF_{W,E}}{dt}} & = & B\beta_{hm}{\displaystyle \frac{I_{h}}{N_{h}}F_{W,S}-(\nu_{m}+\mu_{S})F_{W,E},}\\
{\displaystyle \frac{dF_{W,I}}{dt}} & = & \nu_{m}F_{W,E}-\mu_{I}F_{W,I},\\
{\displaystyle \frac{dS_{S}}{dt}} & = &\epsilon_F \Lambda_{tot}+\dfrac{(1-\varepsilon)M_{S}^*}{M+M_{S}^*}r\gamma A-B\beta_{hm}{\displaystyle \frac{I_{h}}{N_{h}}S_{S}-\mu_{S}S_{S},}\\
{\displaystyle \frac{dS_{E}}{dt}} & = & B\beta_{hm}{\displaystyle \frac{I_{h}}{N_{h}}S_{S}-(\nu_{m}+\mu_{S})S_{E},}\\
{\displaystyle \frac{dS_{I}}{dt}} & = & \nu_{m}S_{E}-\mu_{I}S_{I}.
\end{array}
\right.
\end{equation}

In the sequel, we provide qualitative results of system (\ref{Human-ode-good})-(\ref{Mosquitoes-ode-SIT-constant}). Let us set $$x(t)=(S_h(t), I_h(t),R_h(t),A(t),M(t),F_{W,S}(t),F_{W,E}(t),F_{W,I}(t),S_{S}(t),S_E(t),S_I(t))^T.$$
\subsection{Boundedness of solutions and existence of equilibria}

Using similar arguments as in \cite[Lemmas 1 \& 2]{DumontYatat2022}, it is straightforward to obtain the following Lemma 
\begin{lemma}[Boundedness of solutions]\label{boundeness-lemma-bis} The set 		
$$
\begin{array}{cr}
\Gamma=&\left\{x\in\rr^{11}_+:S_h+I_h+R_h=N_h; (A,M)^T\leq \left(A^*,M^*\right)^T; F_{W,S}+F_{W,E}+F_{W,I}\leq  F_{W,S}^*;\right.  \\
     & \left. S_S+S_E+S_I\leq \dfrac{\epsilon_F \Lambda_{tot}+r\gamma A^*}{\mu_{S}} \right\}
\end{array}
$$
is positively invariant for system \eqref{Human-ode-good}-\eqref{Mosquitoes-ode-SIT-constant} where $(A^*,M^*,F_{W,S}^*)^T$ is given by \eqref{Wild-equilibria-definition}.
	
\end{lemma}

Using Theorem \ref{theo-equilibre-sit-only}, page \pageref{theo-equilibre-sit-only}, we deduce: 

\begin{proposition}[Trivial and non-trivial disease-free equilibria]
\label{prop3}
Whatever $\N \varepsilon\geq 0$, system (\ref{Human-ode-good})-(\ref{Mosquitoes-ode-SIT-constant}) always has a trivial disease-free equilibrium, $TDFE$, such that 
\begin{equation}
    TDFE=\left(N_h,0_{\rr^{7}},\dfrac{\epsilon_F \Lambda_{tot}}{\mu_{S}},0_{\rr^{2}}\right)^T.
\label{TDFE}
\end{equation}
\begin{enumerate}
    \item Assume $\N\varepsilon<1$.  Let $\Lambda_M^{crit}$ defined by \eqref{seuil-ms}, page \pageref{seuil-ms}. 
    \begin{itemize}
    \item If $(1-\epsilon_F) \Lambda_{tot}\in(0, \Lambda_M^{crit})$, then system (\ref{Human-ode-good})-(\ref{Mosquitoes-ode-SIT-constant}) has two non-trivial disease-free equilibria $DFE_{1,2}=\left(N_h,0_{\rr^{2}},A_{1,2},M_{1,2},F_{W,S_{1,2}},0_{\rr^{2}},S_{S_{1,2}},0_{\rr^{2}}\right)^T$
		with $(A_{1},M_{1},F_{W,S_{1}})^T<(A_{2},M_{2},F_{W,S_{2}})^T$ and $A_{1,2}$, $M_{1,2}$, $F_{W,S_{1,2}}$, and  $S_{S_{1,2}}$ given in Theorem  \ref{theo-equilibre-sit-only}.
		
		\item If $(1-\epsilon_F) \Lambda_{tot}=\Lambda_M^{crit}$, then system (\ref{Human-ode-good})-(\ref{Mosquitoes-ode-SIT-constant}) has one  non-trivial disease-free equilibrium
$$
DFE_{\diamond}=\left(N_h,0_{\rr^{2}},A_{\diamond},M_{\diamond},F_{W,S_{\diamond}},0_{\rr^{2}},S_{S_{\diamond}},0_{\rr^{2}}\right)^T,
$$
		with $A_{\diamond}$, $M_{\diamond}$, $F_{W,S_{\diamond}}$, and  $S_{S_{\diamond}}$ given in Theorem  \ref{theo-equilibre-sit-only}.
		\end{itemize}
	\item Assume $\N\varepsilon>1$. System (\ref{Human-ode-good})-(\ref{Mosquitoes-ode-SIT-constant}) admits one non-trivial disease-free equilibrium
$$
DFE_{\dag}=\left(N_h,0_{\rr^{2}},A_{\dag},M_{\dag},F_{W,S_{\dag}},0_{\rr^{2}},S_{S_{\dag}},0_{\rr^{2}}\right)^T,    
$$
where $A_{\dag}$, $M_{\dag}$, $F_{W,S_{\dag}}$, and $S_{S_{\dag}}$ are given in Theorem  \ref{theo-equilibre-sit-only}.
		\item Assume that $\N\varepsilon=1$. If $(1-\epsilon_F) \Lambda_{tot}\in(0,\Lambda_{M, \sharp}^{crit})$, where $\Lambda_{M, \sharp}^{crit}=\mu_{M_S}\dfrac{\N-1}{\Q}$, then system (\ref{Human-ode-good})-(\ref{Mosquitoes-ode-SIT-constant}) has the following  non-trivial disease-free equilibrium  $$DFE_{\sharp}=\left(N_h,0_{\rr^{2}},A_{\sharp},M_{\sharp},F_{W,S_{\sharp}},0_{\rr^{2}},S_{S_{\sharp}},0_{\rr^{2}}\right)^T$$
		where $A_{\sharp}$, $M_{\sharp}$, $F_{W,S_{\sharp}}$, and $S_{S_{\sharp}}$ are given in Theorem  \ref{theo-equilibre-sit-only}.
\end{enumerate}

\end{proposition}
Note that using the relation $\N\varepsilon=1$ in the expression of $\Lambda_M^{crit}$, we recover $\Lambda_{M, \sharp}^{crit}$. Thus, in order to simplify the reading of the paper, we will not consider the particular case $\N\varepsilon=1$ in the rest of the paper because most of the forthcoming results are similar to those obtained when $\N\varepsilon<1$.

Following point 1.b) of Theorem \ref{theo-stabilite-equilibre-sit}, page \pageref{theo-stabilite-equilibre-sit}, in the disease-free case, equilibrium $DFE_1$ is unreachable because it is always unstable. Therefore, in addition to $TDFE$, the meaningful disease-free equilibrium of system (\ref{Human-ode-good})-(\ref{Mosquitoes-ode-SIT-constant}) is
\begin{equation}
DFE_{SIT_c}=\left\{\begin{array}{l}
DFE_\dag, \quad\mbox{when}\quad \N\varepsilon>1,\\
\\
DFE_2, \quad\mbox{when}\quad \N\varepsilon<1\quad\mbox{and}\quad(1-\epsilon_F) \Lambda_{tot}\in(0, \Lambda_M^{crit}),\\
\\
DFE_\diamond, \quad\mbox{when}\quad \N\varepsilon<1 \quad\mbox{and}\quad (1-\epsilon_F) \Lambda_{tot}=\Lambda_M^{crit},\\
\\
TDFE, \quad\mbox{when}\quad \N\varepsilon<1\quad\mbox{and}\quad(1-\epsilon_F) \Lambda_{tot}>\Lambda_M^{crit}.\\
\end{array} 
\right.
\end{equation}
\begin{remark}
Note that in the last case, only $TDFE$ exists, while in the other cases $DFE_{SIT_c}$ and $TDFE$ co-exist.
\end{remark}
Using the next generation matrix approach, see e.g. \cite{VanDENDRIESSCHE2002}, the basic reproduction number of  system (\ref{Human-ode-good})-(\ref{Mosquitoes-ode-SIT-constant}) is

\begin{equation}\label{R0-SIT-c}
 \R^2_{0,SIT_c}=\left\{\begin{array}{l}
 \dfrac{\nu_m}{\nu_m+\mu_S}\dfrac{B\beta_{mh}}{\mu_I}\dfrac{B\beta_{hm}}{\nu_h+\mu_h}\dfrac{(F_{W,S_\dag}+S_{S_\dag})}{N_h}, \quad\mbox{when}\quad\N\varepsilon>1,\\
\\
\\
 \dfrac{\nu_m}{\nu_m+\mu_S}\dfrac{B\beta_{mh}}{\mu_I}\dfrac{B\beta_{hm}}{\nu_h+\mu_h}\dfrac{(F_{W,S_2}+S_{S_2})}{N_h}, \quad\mbox{when}\quad\N\varepsilon<1\quad\mbox{and}\quad (1-\epsilon_F) \Lambda_{tot}\in(0, \Lambda_M^{crit}),\\
 \\
 \dfrac{\nu_m}{\nu_m+\mu_S}\dfrac{B\beta_{mh}}{\mu_I}\dfrac{B\beta_{hm}}{\nu_h+\mu_h}\dfrac{(F_{W,S_\diamond}+S_{S_\diamond})}{N_h}, \quad\mbox{when}\quad\N\varepsilon<1\quad\mbox{and}\quad (1-\epsilon_F) \Lambda_{tot} = \Lambda_M^{crit},\\
\\
\dfrac{\nu_m}{\nu_m+\mu_S}\dfrac{B\beta_{mh}}{\mu_I}\dfrac{B\beta_{hm}}{\nu_h+\mu_h}\dfrac{\epsilon_F \Lambda_{tot}}{\mu_S N_h},
\quad\mbox{when}\quad\N\varepsilon<1\quad\mbox{and}\quad (1-\epsilon_F) \Lambda_{tot}>\Lambda_M^{crit}.\\
 \end{array} 
 \right.
\end{equation}
\begin{remark}
In some cases, as expected, $\R^2_{0,SIT_c}$ has two parts: the first part $$\R^2_{0,SIT_c,W}= \dfrac{\nu_m}{\nu_m+\mu_S}\dfrac{B\beta_{mh}}{\mu_I}\dfrac{B\beta_{hm}}{\nu_h+\mu_h}\dfrac{F_{W,S_{\dag,2,\diamond}}}{N_h},$$
is related to the wild susceptible females that are still fertile while the second part, $$\R^2_{0,SIT_c,S}= \dfrac{\nu_m}{\nu_m+\mu_S}\dfrac{B\beta_{mh}}{\mu_I}\dfrac{B\beta_{hm}}{\nu_h+\mu_h}\dfrac{S_{S_{\dag,2,\diamond}}}{N_h},$$
is related to susceptible females, wild and released ones, that are sterile.

The main question is: when $\R^2_{0,SIT_c,W}<1$, is it possible that the releases of sterile females together with the releases of males which are assumed not to be fully sterile imply $\R^2_{0,SIT_c}>1$?
\end{remark}
\begin{remark}\label{remark4}
Since $F_{W,S_{2,\dag,\diamond}}+S_{S_{2,\dag,\diamond}}=\dfrac{r\gamma A_{2,\dag,\diamond}+\epsilon_F \Lambda_{tot}}{\mu_S}$ and $F_{W,S}^*=\dfrac{r\gamma A^*}{\mu_S}$, and using \eqref{R0}, it is interesting to observe that 
\comment{
\begin{equation}\label{R0-SIT-d}
\R^2_{0,SIT_c} = \R_0 ^2
\left\{\begin{array}{l}
\dfrac{r\gamma A_{2}+\epsilon_F \Lambda_{tot}}{r\gamma A^*}, \quad\mbox{when}\quad (1-\epsilon_F) \Lambda_{tot}\in(0, \Lambda_M^{crit}),\\
\\
\dfrac{r\gamma A_{\dag}+\epsilon_F \Lambda_{tot}}{r\gamma A^*}, \quad\mbox{when}\quad (1-\epsilon_F) \Lambda_{tot}=\Lambda_M^{crit}, \\
\\
\dfrac{\epsilon_F \Lambda_{tot}}{r\gamma A^*}, \quad\mbox{when}\quad (1-\epsilon_F) \Lambda_{tot}>\Lambda_M^{crit}, 
 \end{array} 
 \right.
\end{equation}
}

\begin{equation}\label{R0-SIT-d}
 \R^2_{0,SIT_c}=\R_0 ^2\left\{\begin{array}{l}
 \dfrac{\epsilon_F \Lambda_{tot}}{r\gamma A^*}+\dfrac{A_\dag}{A^*}, \quad\mbox{when}\quad\N\varepsilon>1,\\
\\
\dfrac{\epsilon_F \Lambda_{tot}}{r\gamma A^*}+\dfrac{A_2}{A^*}, \quad\mbox{when}\quad\N\varepsilon<1\quad\mbox{and}\quad (1-\epsilon_F) \Lambda_{tot}\in(0, \Lambda_M^{crit}),\\
\\
\dfrac{\epsilon_F \Lambda_{tot}}{r\gamma A^*}+\dfrac{A_\diamond}{A^*}, \quad\mbox{when}\quad\N\varepsilon<1\quad\mbox{and}\quad (1-\epsilon_F) \Lambda_{tot}=\Lambda_M^{crit},\\
\\
\dfrac{\epsilon_F \Lambda_{tot}}{r\gamma A^*}, \quad\mbox{when}\quad\N\varepsilon<1\quad\mbox{and}\quad (1-\epsilon_F) \Lambda_{tot}>\Lambda_M^{crit},\\

 \end{array} 
 \right.
\end{equation}

where $A^*$ is defined in \eqref{Wild-equilibria-definition}, page \pageref{Wild-equilibria-definition}. Thus, clearly, when  $\epsilon_F \Lambda_{tot}$ is too large, i.e. $\epsilon_F \Lambda_{tot}>r\gamma A^*$, we always have $\R^2_{0,SIT_c}>\R_0 ^2$. In this case, if we already have $\R_0^2>1$, then $\R^2_{0,SIT_c}>1$ such that the SIT will fail to lower the epidemiological risk. Conversely, since $A^*>A_{2,\dag,\sharp}$, then $\R^2_{0,SIT_c}<\R_0 ^2$ whenever 
 $\epsilon_F \Lambda_{tot}$ is sufficiently low, i.e.
\begin{equation}
\epsilon_F \Lambda_{tot} <r\gamma (A^*-A_{2,\dag,\diamond}).
\end{equation}
We recover the same result like in \cite{DumontYatat2022} when $\N \varepsilon \leq 1$.
\end{remark}

\begin{remark}
Since $A_{2,\dag,\diamond}$ is an increasing function of $\epsilon_F \Lambda_{tot}$, it is straightforward to deduce that $\R_{0,SIT_c}^2$ increases with respect to $\epsilon_F \Lambda_{tot}$.
\end{remark}

\begin{remark}\label{remark5}
According to \eqref{R0-SIT-d}, when $\N \varepsilon \leq 1$ and $(1-\epsilon_F) \Lambda_{tot}>\Lambda_M^{crit}$, then $\R^2_{0,SIT_c}<1$ iff
\begin{equation}
\label{Lambda_F_critique}
\epsilon_F \Lambda_{tot} < \dfrac{r\gamma A^*}{\R_0^2}=\dfrac{r\gamma (\gamma+\mu_{A,1})(\N-1)}{\mu_{A,2}\R_0 ^2}:=\Lambda^{crit}_F.
\end{equation}
Also, it follows from \eqref{R0-SIT-d} that $$\epsilon_F \Lambda_{tot}>\Lambda_F^{crit}\Rightarrow \R_{0,SIT_c}^2>1.$$
\end{remark}
\begin{remark}
Clearly, $\epsilon_{F}$ has to be chosen such that
\begin{equation}
\epsilon_F < \dfrac{\Lambda^{crit}_F}{\Lambda_{tot}}.
\end{equation}
This result is in complete contradiction with the constant maximal percentage given by IAEA for contamination by sterile females: we can clearly see that the percentage of contamination may depend on the total amount of sterile insects per release.
\end{remark}
Thanks to the case of sterile female contamination, straightforward computations lead to
\begin{proposition}
\label{existence_WIFEE}
    When $\epsilon_F\Lambda_{tot}>\Lambda_F^{crit}$, then there exists a wild insects-free boundary equilibrium, $WIFE$, such that $A^\#=M^\#=F^\#_S=F^\#_E=F^\#_I=0$, $S^\#_S>0$, $S^\#_E>0$, $S^\#_I>0$ and
\begin{equation}\label{QEE}
    \begin{array}{ccl}
    S_{h}^{\#}&=&\dfrac{\mu_{S}+B\beta_{hm}\dfrac{\mu_{h}}{\mu_{h}+\nu_h}}{B\beta_{hm}\dfrac{\mu_{h}}{\mu_{h}+\nu_h}+\mu_{S}\dfrac{\epsilon_{F}\Lambda_{tot}}{\Lambda_{F}^{crit}}} N_{h}, \\
    S_{I}^{\#}&=&\dfrac{\nu_{m}}{\mu_{I}\left(\nu_{m}+\mu_{S}\right)}\left(1-\dfrac{{\displaystyle \mu_{S}}}{\mu_{S}+B\beta_{hm}\dfrac{\mu_{h}}{\mu_{h}+\nu_h}\left(1-\dfrac{S_{h}^{\#}}{N_{h}}\right)}\right)\epsilon_{F}\Lambda_{tot}.
    \end{array}
\end{equation}
\end{proposition}
\begin{proof}
See Appendix \ref{AppendixC}.
\end{proof}

We now have a look at the existence of non-trivial endemic equilibria.

\begin{proposition}
\label{EE_existence}
Assume $\mu_I=\mu_S$.
\begin{itemize}
\comment{
{\color{blue}
\item If $\R_{0,TDFE}^2>1$, or equivalently, $\epsilon_F \Lambda_{tot}>\Lambda_F^{crit}$, then system \eqref{Human-ode-good}-\eqref{Mosquitoes-ode-SIT-constant} admits a unique wild insects-free endemic equilibrium $WIFEE=(S_h^{\sharp},I_h^{\sharp},R_h^{\sharp},0,0,0,0,0,S_S^{\sharp},S_E^{\sharp},S_I^{\sharp})$ where
\begin{equation}\label{QEE}
    \begin{array}{ccl}
     S_I^\sharp  &= & \dfrac{\mu_hN_h}{B\beta_{mh}}\dfrac{1}{1+\dfrac{B\beta_{hm}\mu_h}{\mu_S(\nu_h+\mu_h)}}(\R_{0,TDFE}^2-1), \\
     S_S^\sharp &=& \dfrac{\epsilon_F \Lambda_{tot}}{\mu_S}\dfrac{1+\dfrac{\mu_S(\nu_m+\mu_S)N_h\mu_h}{\epsilon_F \Lambda_{tot}\nu_mB\beta_{mh}}}{1+\dfrac{B\beta_{hm}}{\mu_S}\dfrac{\mu_h}{\nu_h+\mu_h}},\\
     S_E^\sharp &=& \dfrac{\mu_S}{\nu_m}S_I^\sharp,\\
         I_h^\sharp &=& \dfrac{B\beta_{mh}}{\nu_h+\mu_h}\dfrac{S_I^\sharp}{N_h}\dfrac{\mu_hN_h}{\mu_h+B\beta_{mh}\dfrac{S_I^\sharp}{N_h}},\\
         S_h^\sharp &=& \dfrac{\mu_hN_h}{\mu_h+B\beta_{mh}\dfrac{S_I^\sharp}{N_h}},\\
         R_h^\sharp &=& \dfrac{\nu_h}{\mu_h}I_h^\sharp.\\
    \end{array}
\end{equation}
}}
\item Let $\N \varepsilon \leq 1$, and set
\begin{equation}
\label{seuil_EE_SIT}
    \Lambda_{M,EE}^{crit}=\dfrac{\mu_{M_S}}{\Q}\left(\sqrt{\N+\left(1-\N\varepsilon\right)}-\sqrt{1-\N\varepsilon}\right)^{2}.
\end{equation}
Assume $0<(1-\epsilon_F) \Lambda_{tot}<\Lambda_{M,EE}^{crit}$, and $\epsilon_F \Lambda_{tot} \geq 0$ is chosen such that 
\begin{equation}
    \epsilon_F \Lambda_{tot}+ r\gamma A_{1}^{EE}>\dfrac{F_{W,S}^{*}}{\R_{0}^{2}},
    \label{cond_existence_EE_SIT}
\end{equation}
where
\[
A_{1}^{EE}=\dfrac{1}{2\Q\dfrac{(1-r)\gamma}{\mu_{M}}}\left(\N-\Q M_{S}^{*}- \sqrt{\left(\left(\Q M_{S}^{*}-\N\right)^{2}-4\Q\left(1-\N\varepsilon\right)M_{S}^{*}\right)}\right).
\]
Then there exists two endemic equilibria, $EE_{SIT,1}$  and $EE_{SIT,2}$. In addition $EE_{SIT,1}=EE_{SIT,2}$ when $\N \varepsilon=1$.
\item Let $\N \varepsilon > 1$. For all $(1-\epsilon_F) \Lambda_{tot}>0$, assume that $\epsilon_F \Lambda_{tot} \geq 0$ is chosen such that 
$$
\epsilon_F \Lambda_{tot}+ r\gamma A_{*}^{EE}>\dfrac{F_{W,S}^{*}}{\R_{0}^{2}},
$$
where
\[
A_{1}^{EE}=\dfrac{1}{2\Q\dfrac{(1-r)\gamma}{\mu_{M}}}\left(\N-\Q M_{S}^{*}- \sqrt{\left(\left(\Q M_{S}^{*}-\N\right)^{2}+4\Q\left(\N\varepsilon -1\right)M_{S}^{*}\right)}\right).
\]
Then, there exists one positive equilibrium $EE_{SIT,*}$.
\end{itemize}
\end{proposition}
\begin{proof}
See Appendix \ref{AppendixC}.
\end{proof}

We consider the case where $\mu_S<\mu_I$. We first set the following thresholds
$$
\begin{array}{ccl}
    \alpha &=& \dfrac{\nu_m}{\mu_I}\dfrac{B\beta_{mh}}{\nu_h+\mu_h}\dfrac{B\beta_{hm}}{\nu_m+\mu_S}\dfrac{1}{N_h^2},
    \end{array}
$$
\begin{equation}
\label{Lambda_EE_1}
    \Lambda_{F,EE}^{crit,1}=\dfrac{ \mu_S}{\epsilon_F\alpha}\dfrac{\N\varepsilon\left(1-\dfrac{1+\dfrac{\nu_m}{\mu_I}}{1+\dfrac{\nu_m}{\mu_S}}\right)}{1-\N \varepsilon \left(\dfrac{1+\dfrac{\nu_m}{\mu_I}}{1+\dfrac{\nu_m}{\mu_S}} \right)},
\end{equation}
\begin{equation}
\label{Lambda_EE_2}
    \Lambda_{tot}^{crit,2}=\dfrac{r\gamma(\gamma+\mu_{A,1})\left(\N\dfrac{1+\dfrac{\nu_{m}}{\mu_{I}}}{1+\dfrac{\nu_{m}}{\mu_{S}}}-1\right)}{\mu_{A,2}\left((1-\epsilon_F)\dfrac{r}{1-r}\dfrac{\mu_{M}}{\mu_{M_{S}}}+\epsilon_F\right)},
    \end{equation}
\begin{equation}
\label{Lambda_EE_3}
\Lambda_{tot}^{crit,3}=\dfrac{1}{2\dfrac{\Q(1-\epsilon_F)}{\mu_{M_{S}}}\alpha\epsilon_F}\left[\sqrt{\Delta}+\left(\dfrac{\alpha\mu_{M}(1-\epsilon_F)r}{\left(1-r\right)\mu_{M_{S}}}\left(1-\N\varepsilon\dfrac{1+\dfrac{\nu_{m}}{\mu_{I}}}{1+\dfrac{\nu_{m}}{\mu_{S}}}\right)+\alpha\epsilon_F\left(1-\N\dfrac{1+\dfrac{\nu_{m}}{\mu_{I}}}{1+\dfrac{\nu_{m}}{\mu_{S}}}\right)\right)\right],
\end{equation}
    where 
    $$
    \Delta=\left(\left(\dfrac{\alpha\mu_{M}(1-\epsilon_F)r}{\left(1-r\right)\mu_{M_{S}}}\left(1-\N\varepsilon\dfrac{1+\dfrac{\nu_{m}}{\mu_{I}}}{1+\dfrac{\nu_{m}}{\mu_{S}}}\right)+\alpha\epsilon_F\left(1-\N\dfrac{1+\dfrac{\nu_{m}}{\mu_{I}}}{1+\dfrac{\nu_{m}}{\mu_{S}}}\right)\right)\right)^{2}+4\dfrac{\Q(1-\epsilon_F)}{\mu_{M_{S}}}\alpha\epsilon_F\N\mu_{S}\left(1-\dfrac{1+\dfrac{\nu_{m}}{\mu_{I}}}{1+\dfrac{\nu_{m}}{\mu_{S}}}\right)>0.
    $$
Then, we derive
\begin{proposition}\label{EE-existence_general}
    \begin{itemize}
    Assume $\mu_S<\mu_I$.
        \item Let $\N \varepsilon \leq \dfrac{1+\dfrac{\nu_m}{\mu_S}}{1+\dfrac{\nu_m}{\mu_I}}$.
        \begin{itemize}
        \item If $\Lambda_{tot}^{crit,1}<\Lambda_{tot}< \Lambda_{tot}^{crit,3}$, or
        \item If $\Lambda_{tot}>\max \{\Lambda_{tot}^{crit,3},\Lambda_{tot}^{crit,1} \}$, and $\N \geq \dfrac{1+\dfrac{\nu_m}{\mu_S}}{1+\dfrac{\nu_m}{\mu_I}}$ and $\Lambda_{tot} < \Lambda_{tot}^{crit,2}$,
    \end{itemize}
     then, there exist no or $2$ endemic equilibria.
     \begin{itemize}
         \item If $\Lambda_{tot}>\max \{\Lambda_{tot}^{crit,3},\Lambda_{tot}^{crit,1} \}$, and 
        \begin{itemize}
        \item  $\N > \dfrac{1+\dfrac{\nu_m}{\mu_S}}{1+\dfrac{\nu_m}{\mu_I}}$ and $\Lambda_{tot} > \Lambda_{tot}^{crit,2}$, or
        \item $\N <\dfrac{1+\dfrac{\nu_m}{\mu_S}}{1+\dfrac{\nu_m}{\mu_I}}$,
        \end{itemize}
     then, no endemic equilibrium exists.
     \item If $\Lambda_{tot}<\min \{\Lambda_{tot}^{crit,3},\Lambda_{tot}^{crit,1} \}$, then exists one endemic equilibrium.
         \end{itemize}
        \item Let $\N \varepsilon \geq \dfrac{1+\dfrac{\nu_m}{\mu_S}}{1+\dfrac{\nu_m}{\mu_I}}$.
        \begin{itemize}
        \item If $\Lambda_{tot}< \Lambda_{tot}^{crit,3}$, or
        \item If $\Lambda_{tot}> \max \{\Lambda_{tot}^{crit,2}, \Lambda_{tot}^{crit,3} \}$, 
    \end{itemize}
    then, only one endemic equilibrium exists.
    \begin{itemize}
        \item If $\Lambda_{tot}^{crit,3} <\Lambda_{tot}< \Lambda_{tot}^{crit,2}$ 
    \end{itemize}
    then, one endemic equilibrium exists or three endemic equilibria.
    \end{itemize}
\end{proposition}

\begin{proof}
See Appendix \ref{AppendixC}.
\end{proof}

\subsection{Stability analysis of the disease-free equilibria and uniform persistence}
Let us set
\begin{equation}
    \R_{0,TDFE}^{2}=\dfrac{B\beta_{mh}}{\nu_h+\mu_{h}}\dfrac{\nu_{m}}{\left(\nu_{m}+\mu_{S}\right)\mu_{S}}\dfrac{B\beta_{hm}}{N_{h}}\dfrac{\epsilon_F \Lambda_{tot}}{\mu_{S}}=\R_0^2\dfrac{\epsilon_F \Lambda_{tot}}{r\gamma A^*}.
\end{equation}
A straightforward computation of the Jacobian related to system (\ref{Human-ode-good})-(\ref{Mosquitoes-ode-SIT-constant}) at equilibrium $TDFE$ leads to
\begin{theorem}
    Assume $\N\varepsilon<1$ and $\Lambda_{tot}>0$. Let $\epsilon_F\geq 0$ such that $\R_{0,TDFE}^{2}<1$, then, the Trivial Disease-Free Equilibrium, $TDFE$, is locally asymptotically stable, and unstable when $\R_{0,TDFE}^{2}>1$.
    \label{LAS_TDFE}
\end{theorem}

The previous theorem shows that, when $\N \varepsilon<1$, nuisance reduction with SIT is always possible with low  contamination by sterile females, as long as $\Lambda_{tot}>0$, and the wild population is small or not yet established. When the wild population is large or established we need further results.

Using \cite[Theorem 2]{VanDENDRIESSCHE2002}, the stability properties of the biological disease-free equilibrium $DFE_{SIT_c}\in \left\{DFE_{\dag,2,\diamond}, TDFE\right\}$ is summarized as follows.
\begin{theorem}\label{theorem-LAS-SITcont} The following results hold true for system (\ref{Human-ode-good})-(\ref{Mosquitoes-ode-SIT-constant}).

\noindent Assume $\N\varepsilon<1$.
	\begin{enumerate}
	\item Let $(1-\epsilon_F) \Lambda_{tot}\in(0, \Lambda_M^{crit})$
	\begin{enumerate}
		\item If $\R^2_{0,SIT_c}<1$, then $DFE_{2}$, defined in Proposition \ref{prop3}-(1), is locally asymptotically stable.
		\item If $\R^2_{0,SIT_c}>1$, then $DFE_{2}$ is unstable.
	\end{enumerate}
	\item Let $(1-\epsilon_F) \Lambda_{tot} = \Lambda_M^{crit}$
	\begin{enumerate}
		\item If $\R^2_{0,SIT_c}<1$, then $DFE_{\diamond}$, defined in Proposition \ref{prop3}-(1), is locally asymptotically stable.
		\item If $\R^2_{0,SIT_c}>1$, then $DFE_{\diamond}$ is unstable.
	\end{enumerate}
    \item Let $(1-\epsilon_F) \Lambda_{tot}> \Lambda_M^{crit}$.
	\begin{enumerate}
		\item If $\R^2_{0,SIT_c}=\R_{0,TDFE}^{2}<1$, then $TDFE$, defined in Proposition \ref{prop3}, is globally asymptotically stable.
		\item If $\R^2_{0,SIT_c}=\R_{0,TDFE}^{2}>1$, then $TDFE$ is unstable.
	\end{enumerate}
	\end{enumerate}	
Assume $\N\varepsilon>1$.
	\begin{enumerate}
		\item If $\R^2_{0,SIT_c}<1$, then $DFE_{\dag}$, defined in Proposition \ref{prop3}-(2), is locally asymptotically stable.
		\item If $\R^2_{0,SIT_c}>1$, then $DFE_{\dag}$ is unstable.
	\end{enumerate}
\end{theorem}
In fact, when the residual fertility level is low, i.e. $\varepsilon<\dfrac{1}{\N}$, system (\ref{Human-ode-good})-(\ref{Mosquitoes-ode-SIT-constant}) may exhibit a bistable dynamics in the disease-free context. Indeed, based on Theorem \ref{theo-stabilite-equilibre-sit} together with Theorems \ref{LAS_TDFE} and \ref{theorem-LAS-SITcont}, it is straightforward to establish: 

\begin{theorem}\label{Stability-NTDFE-Mosquito-Human-SIT} 	 Assume $\N\varepsilon<1$ and $(1-\epsilon_F) \Lambda_{tot}\in(0, \Lambda_{M}^{crit})$. If $\R^2_{0,SIT_c}<1$, then equilibria $DFE_{2}$ and $TDFE$ are locally asymptotically stable (LAS). 
\end{theorem}

Clearly, from the two previous theorems, when contamination by sterile females is low, such that $\R_{0,TDFE}^{2}<1$, we derive that:
\begin{itemize}
    \item nuisance reduction is only possible when $\N \varepsilon<1$. In particular, for established wild population, massive sterile insects releases can drive the wild population close to $TDFE$.
\item reducing the epidemiological risk is possible whatever the values taken by $\N \varepsilon$.
\end{itemize}

\begin{remark}
\label{Claim-1}
Based on a comparison argument and a limit system argument we observe the following:
    \begin{itemize}
        \item System (\ref{Human-ode-good})-(\ref{Mosquitoes-ode-SIT-constant}) may undergo a bistability involving the wild insects-free boundary equilibrium, WIFE and the `full' endemic equilibrium $EE$ when $\N\varepsilon\leq1$, $\R^2_{0,TDFE}>1$ and $(1-\epsilon_F) \Lambda_{tot}\in(0, \Lambda_{M}^{crit})$.
        \item The wild insects-free boundary equilibrium, WIFE is GAS when $\N\varepsilon\leq1$, $\R^2_{0,TDFE}>1$ and $(1-\epsilon_F) \Lambda_{tot}> \Lambda_{M}^{crit}$.
    \end{itemize}
\end{remark}

In order to deal with the uniform persistent of system (\ref{Human-ode-good})-(\ref{Mosquitoes-ode-SIT-constant}), we prove the following result: 
\begin{theorem}\label{Persistence}
    If $\N\varepsilon>1$ and $\R^2_{0,TDFE}>1$, then the system is uniformly persistent.
\end{theorem}

\begin{proof}
See Appendix \ref{AppendixE}.
\end{proof}
However, the previous result does not give information on how SIT can impact $\R^2_{0,SIT_c}$.
\subsection{Impact of insect releases on the SIT basic reproduction number}
Now, we want to find $\Lambda_{tot}$ and $\epsilon_F$, such that the epidemiological risk is low, i.e. lead $\R^2_{0,SIT_c}<1$. 

As stated in Remark \ref{remark5}, page \pageref{remark5}, if $\epsilon_F \Lambda_{tot}$ is large, that is $\epsilon_F \Lambda_{tot}>\Lambda_F^{crit}$, then whatever the release rate of sterile males $(1-\epsilon_F) \Lambda_{tot}$ is, we will always have $\R^2_{0,SIT_c}>1$. Hence, in the sequel, we first assume that  $$\epsilon_F \Lambda_{tot}<\Lambda_F^{crit}.$$ Moreover, following Remark \ref{remark4}, page \pageref{remark4},
$\R^2_{0,SIT_c}\leq \R^2_{0}$ 
iff $\epsilon_F \Lambda_{tot}$ is sufficiently low. However, this does not necessarily imply that there exists $(1-\epsilon_F) \Lambda_{tot}>0$ such that $\R^2_{0,SIT_c}<1$.
Straightforward computations lead:
\begin{equation}\label{observations-preliminaires}
\left\{
\begin{array}{ccl}
    A_2&=&\dfrac{1}{2}A^*\left(1-\dfrac{\Q_S}{\N-1}\right)\left(1+\sqrt{1-\dfrac{4\Q_S(1-\N\varepsilon)}{\left(\N-1-\Q_S\right)^2}}\right)>0,\, \mbox{when }\N\varepsilon \leq 1\quad\mbox{and}\quad (1-\epsilon_F) \Lambda_{tot}\in(0, \Lambda_M^{crit})\\
    A_\diamond&=& \dfrac{A^{*}}{\mathcal{N}-1}\left(\sqrt{\dfrac{\left(1-\varepsilon\right)\mathcal{N}}{\left(1-\mathcal{N}\varepsilon\right)}}-1\right), \qquad \qquad \qquad \qquad \qquad \qquad \mbox{when }\N\varepsilon < 1\quad\mbox{and}\quad (1-\epsilon_F) \Lambda_{tot}=\Lambda_M^{crit} \\
 A_\dag &=& \dfrac{1}{2}A^*\left(1-\dfrac{\Q_S}{\N-1}+\sqrt{\left(1-\dfrac{\Q_S}{\N-1}\right)^2+\dfrac{4\Q_S(\N\varepsilon-1)}{\left(\N-1\right)^2}}\right)>0,\, \mbox{when}\quad\N\varepsilon>1.\\
\end{array}
\right.
\end{equation}

Using \eqref{R0-SIT-d}, \eqref{Lambda_F_critique} and \eqref{observations-preliminaires}, we deduce that

\begin{equation}\label{link-R0-LambdaF}
\R^2_{0,SIT_c} =
\left\{
	\begin{array}{l}
\dfrac{\epsilon_F \Lambda_{tot}}{\Lambda_F^{crit}} + \dfrac{\R^2_{0}}{2}\left(1-\dfrac{\Q_S}{\N-1}\right)\left(1+\sqrt{1-\dfrac{4\Q_S(1-\N\varepsilon)}{\left(\N-1-\Q_S\right)^2}}\right),\,\mbox{when }\N\varepsilon\leq 1\quad\mbox{and}\quad (1-\epsilon_F) \Lambda_{tot}\in(0, \Lambda_M^{crit}),\\
\\
\dfrac{\epsilon_F \Lambda_{tot}}{\Lambda_F^{crit}} + \dfrac{\R^2_{0}}{\mathcal{N}-1}\left(\sqrt{\dfrac{\left(1-\varepsilon\right)\mathcal{N}}{\left(1-\mathcal{N}\varepsilon\right)}}-1\right), \qquad \qquad \qquad \qquad \qquad \mbox{when }\N\varepsilon < 1\quad\mbox{and}\quad (1-\epsilon_F) \Lambda_{tot}=\Lambda_M^{crit} \\
\\
\dfrac{\epsilon_F \Lambda_{tot}}{\Lambda_F^{crit}}, \quad\quad\quad\quad\quad\quad \qquad \quad\quad\quad\quad\quad\quad\quad \quad\quad\quad\quad\quad\quad\quad \quad\mbox{when}\quad\N\varepsilon \leq 1\quad\mbox{and}\quad (1-\epsilon_F) \Lambda_{tot}>\Lambda_M^{crit},\\

\dfrac{\epsilon_F \Lambda_{tot}}{\Lambda_F^{crit}} + \dfrac{\R^2_{0}}{2}
\left(1-\dfrac{\Q_S}{\N-1}+\sqrt{\left(1-\dfrac{\Q_S}{\N-1}\right)^2+\dfrac{4\Q_S(\N\varepsilon-1)}{\left(\N-1\right)^2}}\right), \quad\mbox{when}\quad\N\varepsilon>1.\\
\\
 \end{array} 
	\right.
\end{equation}

It is straightforward to obtain the following result.
\begin{lemma}\label{Theo-observation-lambdaF}
    \begin{enumerate}
        \item If $\epsilon_F \Lambda_{tot}>\Lambda_F^{crit}$, then $\R^2_{0,SIT_c}>1$.
        \item Assume that $\N\varepsilon \leq 1$ and $(1-\epsilon_F) \Lambda_{tot}>\Lambda_M^{crit}$. Then $\R^2_{0,SIT_c}<1$ iff $0<\epsilon_F \Lambda_{tot}<\Lambda_F^{crit}$.
    \end{enumerate}
\end{lemma}

Lemma \ref{Theo-observation-lambdaF} depicts the fact that when the epidemiological risk is high, that is, when $\R_0^2>1$, and if in addition the release rate of sterile females is large, that is $\epsilon_F \Lambda_{tot}>\Lambda_F^{crit}$, then whatever the amount of released sterile males, the SIT will fail since we will always have $\R^2_{0,SIT_c}>1$. However, massive releases of sterile males ($(1-\epsilon_F) \Lambda_{tot}>\Lambda_M^{crit}$) could be successful provided that $\epsilon_F \Lambda_{tot}<\Lambda_F^{crit}$.

The next question to investigate deals with the possibility to lower the epidemiological risk using small sterile males releases when $\epsilon_F \Lambda_{tot}<\Lambda_F^{crit}$ and also to investigate if there exist necessary conditions to ensure that $\R^2_{0,SIT_c}<1$ when $\N \varepsilon >1$.

\comment{
We set 
\begin{equation}
\label{seuil_R0}
    \R_{0,*}^2=\dfrac{1}{\dfrac{\epsilon_F \Lambda_{tot}}{\Lambda_F^{crit}}+1-\dfrac{Q\Lambda_M}{\mu_{M_S}(\R-1)}}.
    \label{seuilMTC}
\end{equation}
}

\subsection{When $\N \varepsilon<1$}
Using \eqref{link-R0-LambdaF}$_2$, we define the following threshold
\begin{equation}
    \R_{0, \N\varepsilon<1}^2=\dfrac{\mathcal{N}-1}{\dfrac{\epsilon_F \Lambda_{tot}\mu_{A,2}}{r\gamma(\gamma+\mu_{A,1})}+\sqrt{\dfrac{\left(1-\varepsilon\right)\mathcal{N}}{\left(1-\mathcal{N}\varepsilon\right)}}-1}.
    \label{seuil_R0}
\end{equation}
We derive the following result
\begin{theorem}\label{theorem-LAS-SITcont-discussion-Cgt0}Assume $0 \leq \epsilon_F \Lambda_{tot}<\Lambda_F^{crit}$. Consider system (\ref{Human-ode-good})-(\ref{Mosquitoes-ode-SIT-constant}) and set
\begin{equation}
\Lambda_{M,\R_0^2,\varepsilon}^*=\dfrac{\mu_{M_{S}}\left(\N-1\right)}{Q}\left(1-\dfrac{\R_{0}^{4}(1-\N\varepsilon)+\left(\N-1\right)\left(1-\dfrac{\epsilon_F \Lambda_{tot}}{\Lambda_{F}^{crit}}\right)^{2}}{\R_{0}^{4}(1-\N\varepsilon)+\R_{0}^{2}\left(\N-1\right)\left(1-\dfrac{\epsilon_F \Lambda_{tot}}{\Lambda_{F}^{crit}}\right)}\right).
\label{seuilMTC}
\end{equation}
	\begin{enumerate}
		\item If $\R^2_{0}\geq \R_{0, \N\varepsilon<1}^2$, then the following results hold true: 
  \begin{itemize}
      \item  When $(1-\epsilon_F) \Lambda_{tot}>\Lambda_M^{crit}$, the equilibrium $TDFE$ is globally asymptotically stable.
      \item  When $(1-\epsilon_F) \Lambda_{tot}\leq \Lambda_M^{crit}$, then $\R^2_{0,SIT_c}>1$ and SIT fails.
  \end{itemize}

  
		\item If $1<\R^2_{0}<\R_{0, \N\varepsilon<1}^2$, then the following results hold true: 
		\begin{itemize}
		\item When $(1-\epsilon_F) \Lambda_{tot}>\Lambda_M^{crit}$, the equilibrium $TDFE$ is globally asymptotically stable. 
		
		\item When $(1-\epsilon_F) \Lambda_{tot}=\Lambda_M^{crit}$, then $\R^2_{0,SIT_c}<1$, $DFE_{\diamond}$
		and TDFE are locally asymptotically stable. The set $$\{(S,I,R,A,M,F_{W,S},F_{W,E},F_{W,I},S_S,S_E,S_I)^T\in\rr^{11}_+:(A,M,F_{W,S})^T<(A_\diamond,M_\diamond,F_{W,S_\diamond})^T\}$$ belongs to the basin of attraction of TDFE while the set 
		$$\{(S,I,R,A,M,F_{W,S},F_{W,E},F_{W,I},S_S,S_E,S_I)^T\in\rr^{11}_+:(A,M,F_{W,S})^T\geq(A_\diamond,M_\diamond,F_{W,S_\diamond})^T\}$$ belongs to the basin of attraction of $DFE_{\diamond}$.
		 \item when $(1-\epsilon_F) \Lambda_{tot}>\Lambda_{M,\R_0^2,\varepsilon}^*$, then $\R^2_{0,SIT_c}<1$, and the equilibria $DFE_{2}$ and $TDFE$ are locally asymptotically stable. Moreover, the set $$\{(S,I,R,A,M,F_{W,S},F_{W,E},F_{W,I},S_S,S_E,S_I)^T\in\rr^{11}_+:(A,M,F_{W,S})^T<(A_1,M_1,F_{W,S_1})^T\}$$ belongs to the basin of attraction of TDFE while the set 
		$$\{(S,I,R,A,M,F_{W,S},F_{W,E},F_{W,I},S_S,S_E,S_I)^T\in\rr^{11}_+:(A,M,F_{W,S})^T>(A_1,M_1,F_{W,S_1})^T\}$$ belongs to the basin of attraction of $DFE_{2}$.
		 \end{itemize}
	\end{enumerate}
\end{theorem}
\begin{proof}
We follow the same methodology used in \cite[Theorem 6]{DumontYatat2022} to derive \eqref{seuilMTC}. Then, the results follow from Theorem \ref{Stability-NTDFE-Mosquito-Human-SIT}, page \pageref{Stability-NTDFE-Mosquito-Human-SIT}.
\end{proof}
\begin{remark}
Of course, when $\varepsilon=0$, we recover the results obtained in \cite{DumontYatat2022}.
\end{remark}
Clearly the constraint on the releases size given by \eqref{seuilMTC} can be strong, i.e. close to $\Lambda_M^{crit}$, such that it seems to be preferable to use massive releases, i.e. $(1-\epsilon_F) \Lambda_{tot}>\Lambda_M^{crit}$.

In that case, the strategy developed in \cite{Anguelov2020b,Anguelov2020}, using massive and then small releases can be adequate to reduce the epidemiological risk and maintain this risk at a lower level.

\vspace{1cm}
Thus, in terms of vector control: when $\R_0^2\leq 1$, vector control is not necessary; 
when $\R_0^2>1$ and $0\leq \epsilon_F \Lambda_{tot}<\Lambda_F^{crit}$,  then two cases should be considered:
\begin{itemize}
\item when $\R_0^2\geq \R_{0, \N\varepsilon<1}^2$, then massive releases of sterile insect,  i.e.  $(1-\epsilon_F) \Lambda_{tot}>\Lambda_M^{crit}$, should be advocated. 
\item When $\R_0^2<\R_{0, \N\varepsilon<1}^2$, then small, but large enough ($\Lambda_{M,\R_0^2,\varepsilon}^*<(1-\epsilon_F) \Lambda_{tot}\leq \Lambda_M^{crit}$), releases of sterile insects could be useful to control the disease. However, since $\Lambda_{M,\R_0^2,\varepsilon}^*$ is close to $\Lambda_M^{crit}$, from a practical point of view, it is preferable to consider massive releases of sterile insects too.
\end{itemize}
When $\N \varepsilon\leq 1$, we summarize all qualitative results of system (\ref{Human-ode-good})-(\ref{Mosquitoes-ode-SIT-constant}) related to the disease-free equilibria in Table \ref{table:recapitulatif}, page \pageref{table:recapitulatif}.

\begin{table}[H]
 \centering
\begin{tabular}{|l|l|c|c|c|l|}
  \hline
  $\N$ & $\mathcal{R}_0^2$ & $\epsilon_F \Lambda_{tot}$ & $\mathcal{R}^2_{0}$ & $(1-\epsilon_F) \Lambda_{tot}$ & Observations \\
  \hline
   $\leq1$ & & & & & $TDFE$ is GAS\\
  \hline
     & $\leq1$ & & & & Releases of sterile insects are useless  \\
    &  & & & &  because the $DFE$ is already GAS\\
   \cline{2-6}
      $>1$  &  & $\geq \Lambda_F^{crit}$ & &  & Even massive releases could not be efficient   \\
    &  & & & &  to reduce the epidemiological risk: $\R_{0,SIT_c}^2>1$.\\
    &  & & & &   WIFE and/or EE are/is LAS\\
   \cline{3-6}
                     & $>1$  &  & $\geq\R_{0, \N\varepsilon<1}^2$ & $>\Lambda_M^{crit}$ & $TDFE$ is GAS\\
    \cline{5-6}
                     & &  &  & $\leq \Lambda_M^{crit}$ & SIT failed since $\R_{0,SIT_c}^2>1$\\
\cline{4-6}
                     &  &  &  & $>\Lambda_M^{crit}$ & $TDFE$ is GAS  \\
  \cline{5-6}
                     &  & $<\Lambda_F^{crit}$ & $<\R_{0, \N\varepsilon<1}^2$ & $=\Lambda_M^{crit}$ & $\R_{0,SIT_c}^2<1$: $TDFE$ and $DFE_\diamond$ are both stable \\
  \cline{5-6}
                     &  &  &  & $>\Lambda_{M,\R_0^2,\varepsilon}^*$ & $\R_{0,SIT_c}^2<1$: $TDFE$ and $DFE_2$ are both stable \\
  \hline
\end{tabular}
\caption{Summary table of the qualitative analysis of system
(\ref{Human-ode-good})-(\ref{Mosquitoes-ode-SIT-constant}) when $\N \varepsilon\leq 1$.}\label{table:recapitulatif}
\end{table}

\subsection{The case where $\N \varepsilon>1$}
We want to derive if, for a given $\epsilon_F \Lambda_{tot}<\Lambda_F^{crit}$, there exists $\Lambda_{M,\N \varepsilon>1}^{crit}$ such that for all $(1-\epsilon_F) \Lambda_{tot}>\Lambda_{M,\N \varepsilon>1}^{crit}$, we always have $\R^2_{0,SIT_c}>1$. Conversely, for a given $\Lambda_{tot}$ it is possible to find a rate $\epsilon_F$ such that $\R^2_{0,SIT_c}>1$?

Assuming $\R_0^2>1$, $0\leq\epsilon_F \Lambda_{tot}<\Lambda_{F}^{crit}$, and using \eqref{link-R0-LambdaF}$_4$, we have the following:

$\bullet$ Assume that $\dfrac{\epsilon_F \Lambda_{tot}}{\Lambda_F^{crit}} + \dfrac{\R^2_{0}}{2}
\left(1-\dfrac{\Q_S}{\N-1}\right)\geq1$ or equivalently $\left(1-\dfrac{\epsilon_{F}\Lambda_{tot}}{\Lambda_{F}^{crit}}\right)\dfrac{2}{\R_{0}^{2}}-\left(1-\dfrac{\Q_S}{\N-1}\right)\leq0$. Then it holds $$\R^2_{0,SIT_c}>1.$$ Note also that 
$$\left(1-\dfrac{\epsilon_{F}\Lambda_{tot}}{\Lambda_{F}^{crit}}\right)\dfrac{2}{\R_{0}^{2}}-\left(1-\dfrac{\Q_S}{\N-1}\right)\leq 0\Longleftrightarrow (1-\epsilon_F)\Lambda_{tot}\leq \dfrac{\mu_{M_S}}{\Q}(\N-1)\left(1-\dfrac{2\left(1-\dfrac{\epsilon_F \Lambda_{tot}}{\Lambda_F^{crit}}\right)}{\R_0^2}\right):=\Lambda_M^{crit,\sharp}.$$

$\bullet$ Assume that $\dfrac{\epsilon_F \Lambda_{tot}}{\Lambda_F^{crit}} + \dfrac{\R^2_{0}}{2}
\left(1-\dfrac{\Q_S}{\N-1}\right)<1$ or equivalently $\left(1-\dfrac{\epsilon_{F}\Lambda_{tot}}{\Lambda_{F}^{crit}}\right)\dfrac{2}{\R_{0}^{2}}-\left(1-\dfrac{\Q_S}{\N-1}\right)>0$ or equivalently 

$$(1-\epsilon_F)\Lambda_{tot}> (\N-1)\dfrac{\mu_{M_S}}{\Q}\left(1-\dfrac{2\left(1-\dfrac{\epsilon_F \Lambda_{tot}}{\Lambda_F^{crit}}\right)}{\R_0^2}\right):=\Lambda_M^{crit,\sharp}.$$

Let us set
$$\R_{0,\N\varepsilon>1}^2=\dfrac{(\N-1)}{(\N\varepsilon-1)}\left(1-\dfrac{\epsilon_{F}\Lambda_{tot}}{\Lambda_{F}^{crit}}\right).$$
Then we have
$$
\begin{array}{ccl}
  \R^2_{0,SIT_c}>1   & \Longleftrightarrow & \dfrac{\epsilon_F \Lambda_{tot}}{\Lambda_F^{crit}} + \dfrac{\R^2_{0}}{2}
\left(1-\dfrac{\Q_S}{\N-1}+\sqrt{\left(1-\dfrac{\Q_S}{\N-1}\right)^2+\dfrac{4\Q_S(\N\varepsilon-1)}{\left(\N-1\right)^2}}\right)>1,\\
    & \Longleftrightarrow &  1-\dfrac{\Q_S}{\N-1}+\sqrt{\left(1-\dfrac{\Q_S}{\N-1}\right)^2+\dfrac{4\Q_S(\N\varepsilon-1)}{\left(\N-1\right)^2}}>\left(1-\dfrac{\epsilon_{F}\Lambda_{tot}}{\Lambda_{F}^{crit}}\right)\dfrac{2}{\R_{0}^{2}},\\
    & \Longleftrightarrow &  \sqrt{\left(1-\dfrac{\Q_S}{\N-1}\right)^2+\dfrac{4\Q_S(\N\varepsilon-1)}{\left(\N-1\right)^2}}>\left(1-\dfrac{\epsilon_{F}\Lambda_{tot}}{\Lambda_{F}^{crit}}\right)\dfrac{2}{\R_{0}^{2}}-\left(1-\dfrac{\Q_S}{\N-1}\right).\\
    & \Longleftrightarrow & \Q_S\left(\dfrac{\N\varepsilon-1}{(\N-1)^2}-\dfrac{1-\dfrac{\epsilon_{F}\Lambda_{tot}}{\Lambda_{F}^{crit}}}{\R_0^2(\N-1)}\right)>\dfrac{1-\dfrac{\epsilon_{F}\Lambda_{tot}}{\Lambda_{F}^{crit}}}{\R_0^2}\left(\dfrac{1-\dfrac{\epsilon_{F}\Lambda_{tot}}{\Lambda_{F}^{crit}}}{\R_0^2}-1\right),\\
    & \Longleftrightarrow & 
    \Q_S\left(\dfrac{\N\varepsilon-1}{(\N-1)}\dfrac{\R_0^2}{1-\dfrac{\epsilon_{F}\Lambda_{tot}}{\Lambda_{F}^{crit}}}-1\right)>(\N-1)\left(\dfrac{1-\dfrac{\epsilon_{F}\Lambda_{tot}}{\Lambda_{F}^{crit}}}{\R_0^2}-1\right),\\
     & \Longleftrightarrow &
    \Q_S\left(1-\dfrac{\R_0^2}{\R_{0,\N\varepsilon>1}^2}\right)<(\N-1)\left(1-\dfrac{1-\dfrac{\epsilon_{F}\Lambda_{tot}}{\Lambda_{F}^{crit}}}{\R_0^2}\right).\\
\end{array}
$$
Thus, we deduce the two following cases:
\begin{itemize}
    \item[$(i)$] If $\R_0^2>\R_{0,\N\varepsilon>1}^2$, then $\R^2_{0,SIT_c}>1$ for all $(1-\epsilon_F) \Lambda_{tot}>0$.
\item[$(ii)$] If $\R_0^2<\R_{0,\N\varepsilon>1}^2$, then we set  $$\Lambda_{M,\N \varepsilon>1}^{crit}=\dfrac{\mu_{M_S}}{\Q}\dfrac{(\N-1)\left(1-\dfrac{1-\dfrac{\epsilon_{F}\Lambda_{tot}}{\Lambda_{F}^{crit}}}{\R_0^2}\right)}{\left(1-\dfrac{\R_0^2}{\R_{0,\N\varepsilon>1}^2}\right)}$$ and 
we have
    $$
\left\{
\begin{array}{l}
\R^2_{0,SIT_c}>1    \Longleftrightarrow   (1-\epsilon_F) \Lambda_{tot} < \Lambda_{M,\N \varepsilon>1}^{crit} ,  \\
     \R^2_{0,SIT_c}<1    \Longleftrightarrow   (1-\epsilon_F) \Lambda_{tot} > \Lambda_{M,\N \varepsilon>1}^{crit}.  \\
\end{array}
\right.
$$
\end{itemize}
To summarize the previous discussion, when $\N\varepsilon>1$, we have three configurations
\begin{enumerate}
    \item When $(1-\epsilon_F)\Lambda_{tot}\leq\Lambda_{M}^{crit,\sharp}$ or $((1-\epsilon_F)\Lambda_{tot}>\Lambda_{M}^{crit,\sharp}$ and  
  $\R_{0}^{2}>\max(1,\R_{0,\N\varepsilon>1}^2))$, then $\R^2_{0,SIT_c}>1$.
    \item When $1<\R_{0}^{2}<\R_{0,\N\varepsilon>1}^2$ and $(1-\epsilon_F)\Lambda_{tot}>\max(\Lambda_{M,\N \varepsilon>1}^{crit},\Lambda_{M}^{crit,\sharp})$ then $\R^2_{0,SIT_c}<1$.
    \item When $1<\R_{0}^{2}<\R_{0,\N\varepsilon>1}^2$ and $\Lambda_{M}^{crit,\sharp}<(1-\epsilon_F)\Lambda_{tot}<\Lambda_{M,\N \varepsilon>1}^{crit}$ then $\R^2_{0,SIT_c}>1$.
\end{enumerate}
We therefore summarize all qualitative results of system (\ref{Human-ode-good})-(\ref{Mosquitoes-ode-SIT-constant}) related to the disease free equilibria in Table \ref{table: recapitulatif2}, page \pageref{table: recapitulatif2}.

\begin{table}[H]
 \centering
\begin{tabular}{|l|l|c|c|c|l|}
  \hline
  $\N$ & $\mathcal{R}_0^2$ & $\epsilon_F \Lambda_{tot}$ & $\mathcal{R}^2_{0}$ & $(1-\epsilon_F) \Lambda_{tot}$ & Observations \\
  \hline
     & $\leq1$ & & & & Releases of sterile insects are useless  \\
    &  & & & &  because the $DFE$ is already GAS\\
   \cline{2-6}
      $>1$  &  & $\geq \Lambda_F^{crit}$ & &  & Even massive releases could not be efficient  \\
    &  & & & &  to reduce the epidemiological risk\\
   \cline{3-6}
                     &   &  & $\geq\R_{0,\N\varepsilon>1}^2$ &  & SIT fails since $\R_{0,SIT_c}^2>1$ \\
\cline{4-6}
                     &  & $<\Lambda_{F}^{crit}$ & $<\R_{0,\N\varepsilon>1}^2$ & $>\max\{\Lambda_{M,\N \varepsilon>1}^{crit},\Lambda_{M}^{crit,\sharp}\}$ & $\R_{0,SIT_c}^2<1$, $DFE_\dag$ is LAS \\
  \cline{5-6}
                     &  &  &  & $ <\max\{\Lambda_{M}^{crit,\sharp},\Lambda_{M,\N \varepsilon>1}^{crit}\}$ & SIT fails since $\R_{0,SIT_c}^2>1$ \\
  \cline{5-6}
  \hline
\end{tabular}
\caption{Summary table of the qualitative analysis of system
(\ref{Human-ode-good})-(\ref{Mosquitoes-ode-SIT-constant}) when $\N \varepsilon> 1$.}\label{table: recapitulatif2}
\end{table}

\section{Numerical simulations \label{section4}}

\subsection{Sensitivity analysis}

\begin{figure}[h!]
\includegraphics[width=1.0\linewidth]{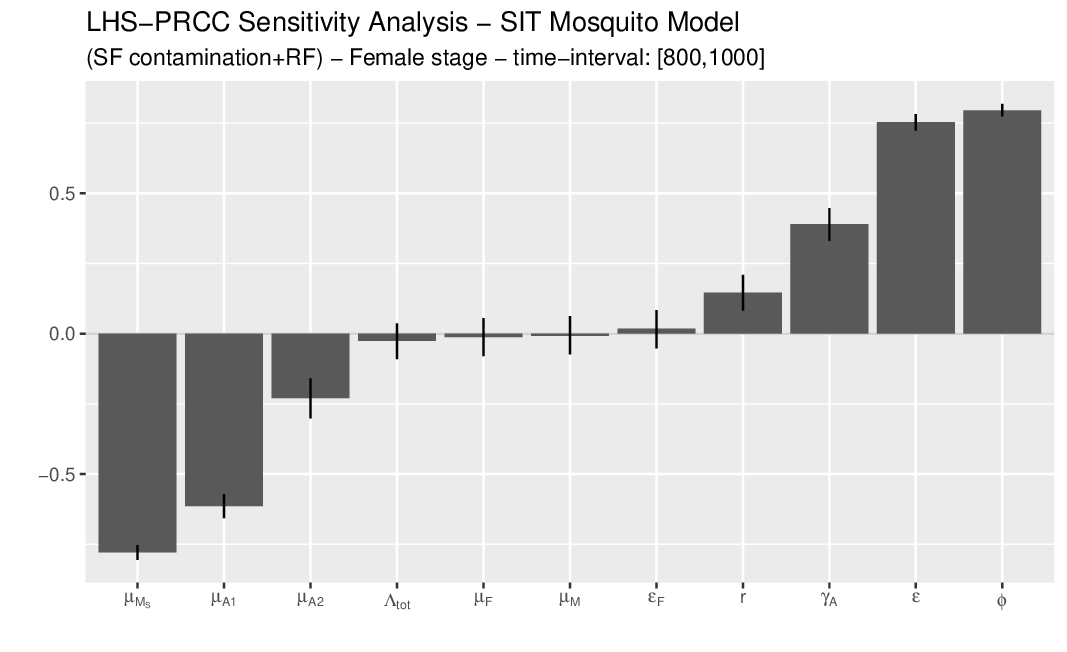} 
    \caption{LHS-PRCC Sensitivity analysis of the Entomological model - Wild Females}
    \label{fig:sensi1}
\end{figure}

\begin{figure}[h!]
    \begin{center}
\includegraphics[width=0.9\linewidth]{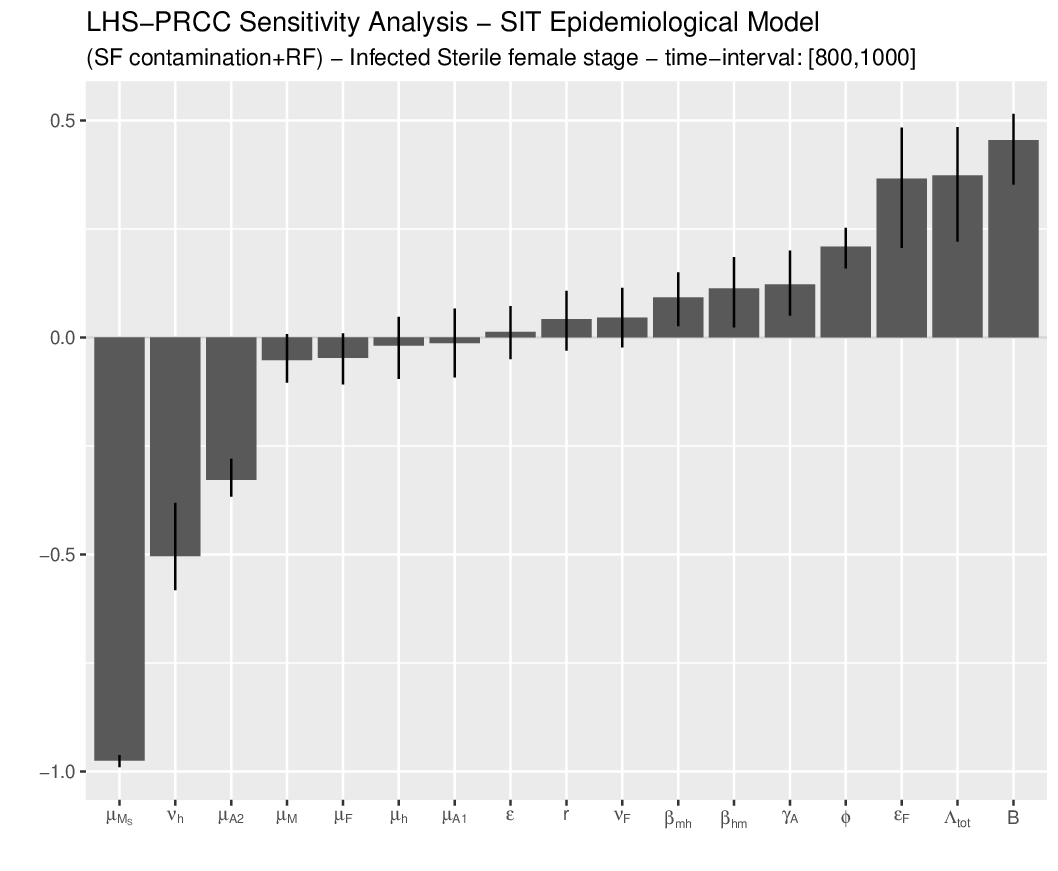} 
\end{center}
    \caption{LHS-PRCC Sensitivity analysis of the Epidemiological model - Infected Sterile Females}
    \label{fig:sensi2}
\end{figure}

\begin{figure}[h!]
    \begin{center}
\includegraphics[width=0.9\linewidth]{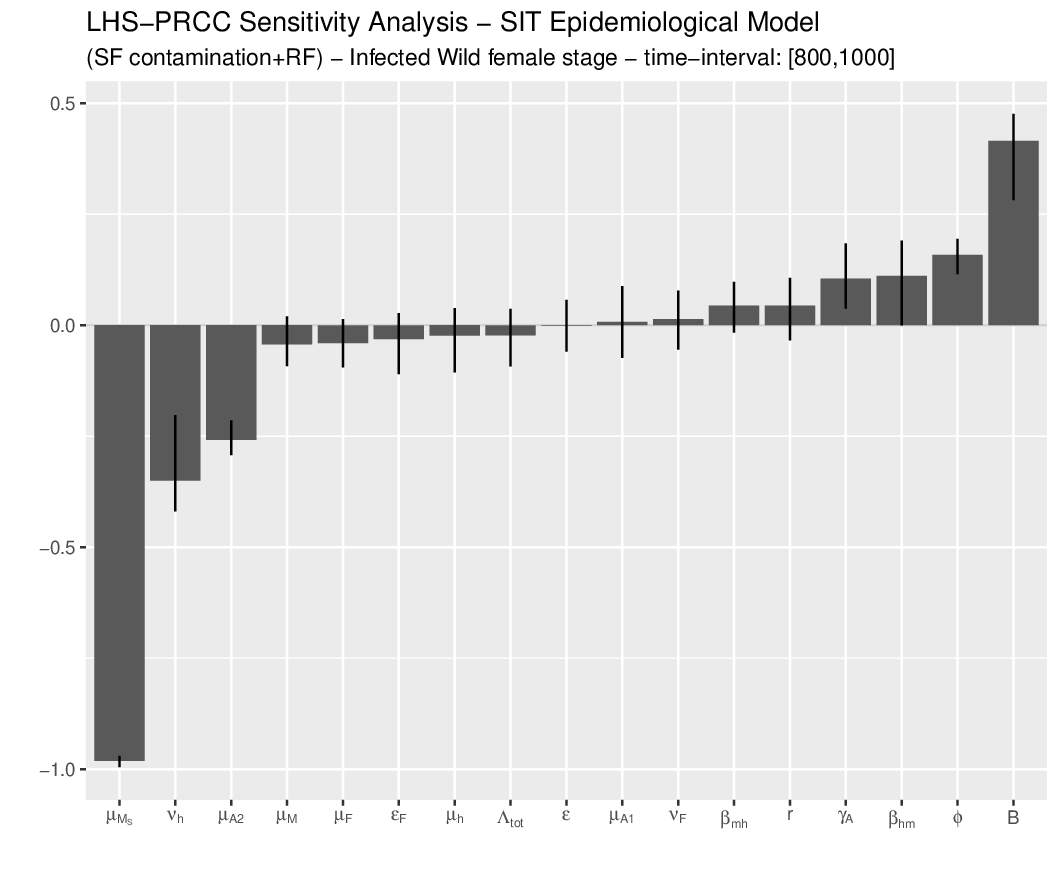} 
\end{center}
    \caption{LHS-PRCC Sensitivity analysis of the Epidemiological model - Infected Wild Females}
    \label{fig:sensi3}
\end{figure}

\begin{figure}[h!]
    \begin{center}
\includegraphics[width=0.9\linewidth]{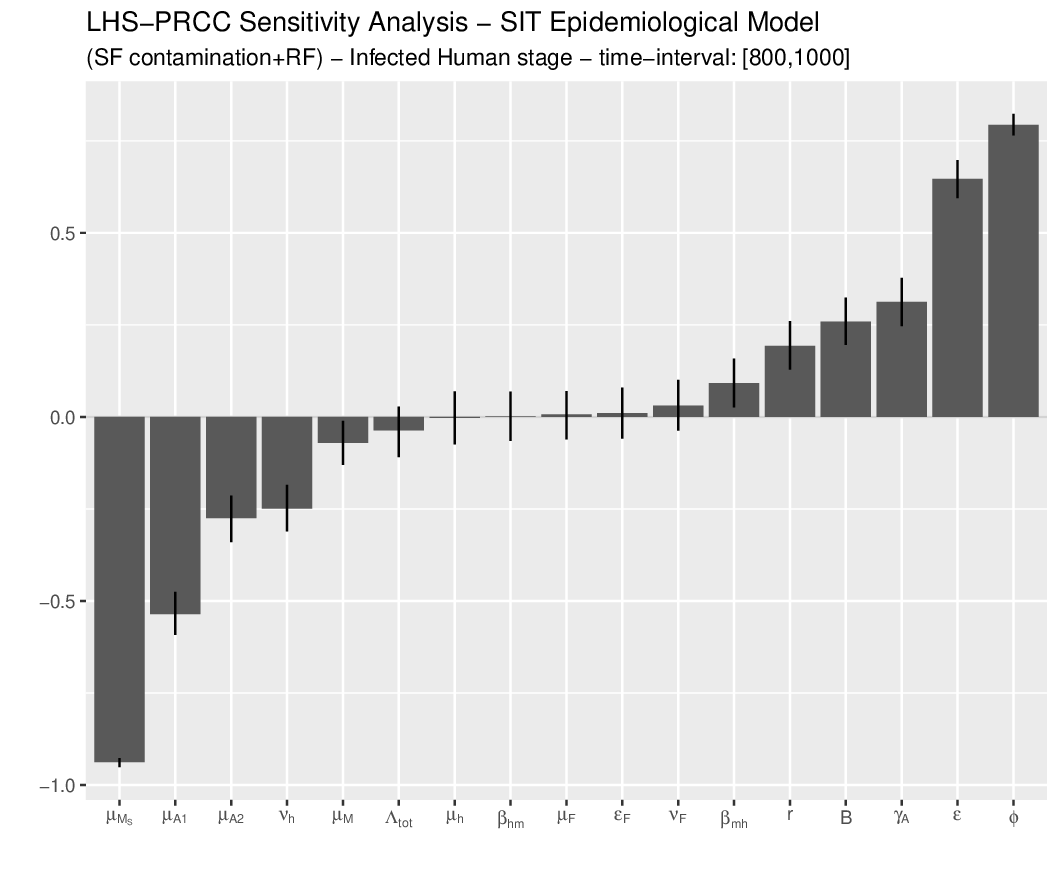} 
\end{center}
    \caption{LHS-PRCC Sensitivity analysis of the Epidemiological model - Infected Humans}
    \label{fig:sensi4}
\end{figure}

It is interesting to study the impact of parameter changes on the dynamics of our systems, and to find which parameters are the most sensitive on the variable outputs.
In Figs \ref{fig:sensi1}, \ref{fig:sensi2}, \ref{fig:sensi3} and \ref{fig:sensi4}, we provide a LHS-PRCC sensitivity analysis, where LHS stands for Latin Hypercube Sampling and PRCC for Partial Rank Correlation Coefficient. The LHS-PRCC method provides mainly information about how the outputs are impacted if we increase (or decrease) the inputs of a specific parameter. The analysis is done on the time interval [800,1000]. The results are ordered from the most negative to the most positive ones.
We derive a LHS-PRCC analysis for the variable $F$ from the entomological model, and the variables $S_I$, $F_{W,I}$ and $I_h$ from the epidemiological model. It is very interesting to compare the impact of the parameters thanks to the considered variables. In Fig. \ref{fig:sensi1}, the parameters $\phi$, $\varepsilon$, $\mu_{M_S}$ and $\mu_{A,1}$ are the parameters for which the Female variable, related to the entomological model \eqref{ODE-entomo-old}, is the more sensitive to. Then, the infected sterile female variable, $S_I$, is mostly sensitive to $\mu_{M_S}$, $\nu_h$, $\mu_{A,2}$, $\epsilon_F$, $\Lambda_{tot}$, and $B$. A similar trend is observed in Fig. \ref{fig:sensi3}, when dealing with wild infected female variable, $F_{W,I}$, except that now $\beta_{hm}$ and $\phi$ are now the main parameters, while $\epsilon_F$ and $\Lambda_{tot}$ not. The residual fertility parameter, $\varepsilon$ has also almost no effect. Finally, considering the infected human variable, $I_h$, it is mostly sensitive to parameters $\mu_{M_S}$, $\mu_{A,1}$ $\varepsilon$ and $\phi$ (see also Fig. \ref{fig:sensi4}). 

We can notice that the two parameters of interest throughout this work $\varepsilon$ and $\epsilon_F$ have a strong impact on $F$, $I_h$, and $S_I$.

For all PRCC analysis, we used the PCC function (R software \cite{R2022}) and $1000$ bootstrap replicates, with a probability level of $0.95$ for (the bootstrap) confidence intervals.

\subsection{Simulations}
All forthcoming numerical simulations are done using the ode23 solver of Matlab \cite{MATLAB2012}. Results are obtained in a couple of seconds.

Like in \cite{DumontYatat2022}, we will consider the effective reproduction number, $\R_{eff}(t)$ for all time $t>0$. Indeed, SIT control is a long term strategy and the starting time of SIT treatment is important thanks to the starting time of the risky period from the epidemiological point of view, that is when Dengue virus starts to circulate, $t_{DENV}$. That is why, it is important to consider the effective reproduction number, $\R_{eff}(t)$, that is defined as follow
\begin{equation}
    \R_{eff}(t)=\dfrac{\nu_m}{\nu_m+\mu_S}\dfrac{B^2\beta_{mh}\beta_{hm}}{\mu_I\left(\nu_h+\mu_h\right)}\dfrac{F_{W,S}(t)+S_S(t)}{N_h}.
\end{equation}
In particular, we will estimate $\R_{eff}$ at time $t_{DENV}$. Clearly, if $\R_{eff}(t_{DENV})<1$ and $\R_{0,SIT_c}^2<1$, then no epidemics will occur. In contrary, even if $\R_{0,SIT_c}^2<1$ but $\R_{eff}(t_{DENV})>1$ then an outbreak will occur.

We consider the parameter values defined in Table \ref{table:parameters}, page \pageref{table:parameters}. For these values we derive $\N \approx 86.75$. This is a high value but meaningful since we have considered the ``best" case for the mosquito dynamics, i.e. the most difficult case in terms of control. For the epidemiological parameters, at a mean temperature of $T=25^\circ C$, we find that $\R_0^2\approx 7.298$, which is quite  large value. 

Then, according to formula \eqref{Lambda_F_critique} and the parameters values, the critical sterile females release rate, $\Lambda_F^{crit}$, is around $391$.

\begin{figure}[h!]
\includegraphics[width=1.0\linewidth]{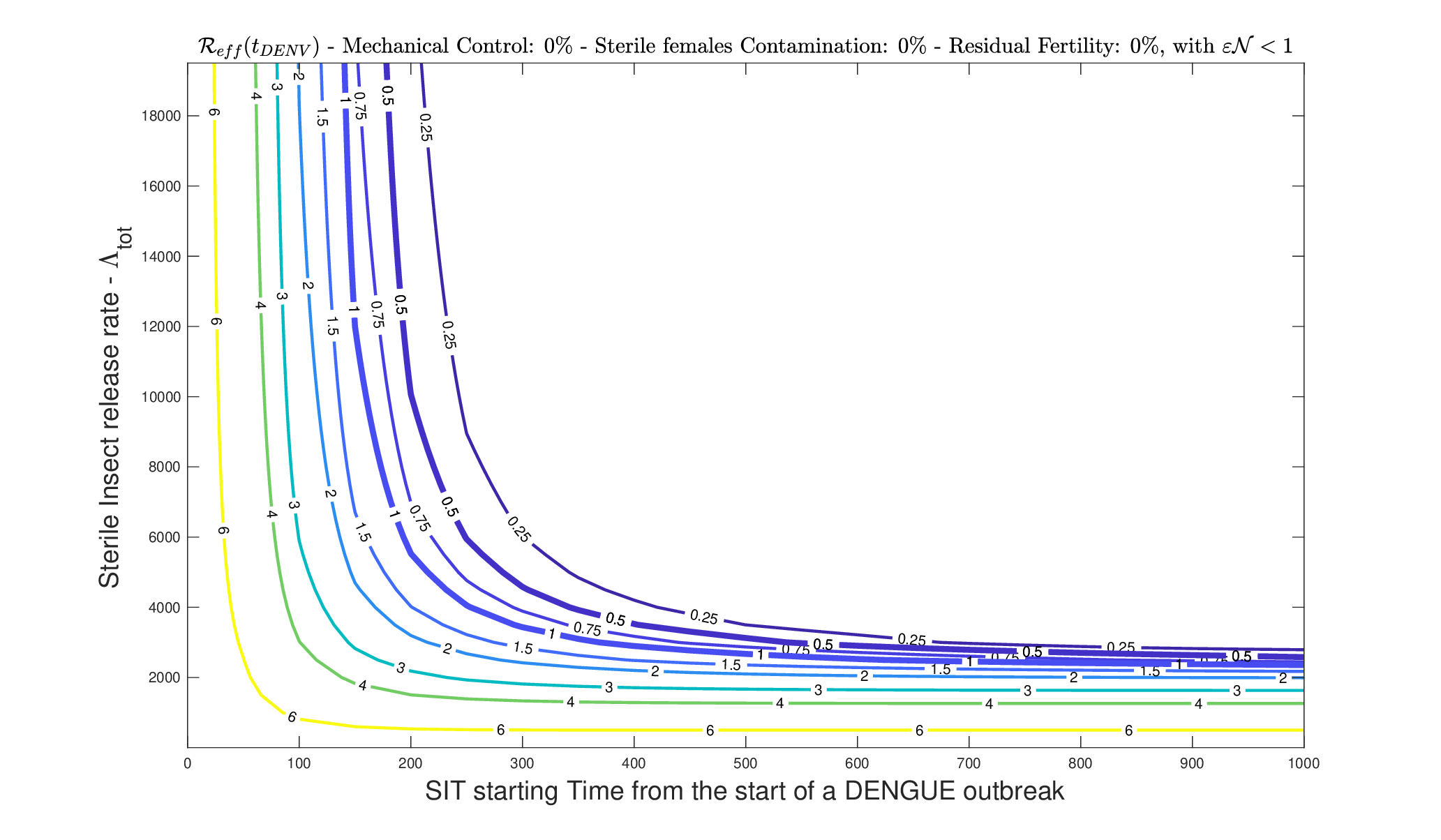} 
    \caption{$\R_{eff}(t_I)$ vs the starting time without sterile female contamination, without residual fertility}
    \label{fig:2a}
\end{figure}
\begin{figure}[h!]
\includegraphics[width=1.0\linewidth]{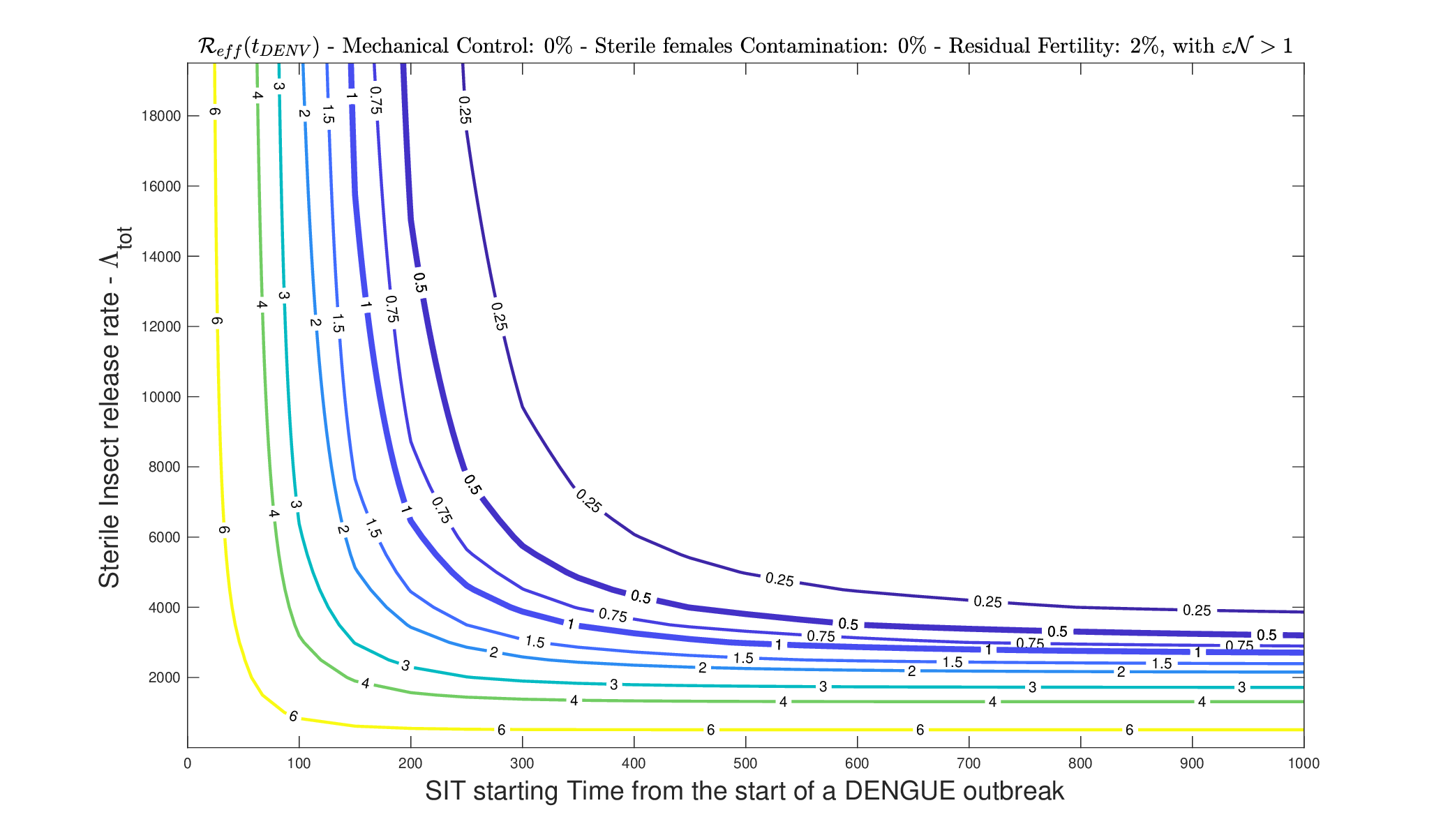} 
    \caption{$\R_{eff}(t_I)$ vs the starting time and the level of the control, without contamination by sterile females, and with $2\%$ of residual fertility, without Mechanical control}
    \label{fig:2b}
\end{figure}
\begin{figure}[h!]
\includegraphics[width=0.9\linewidth]{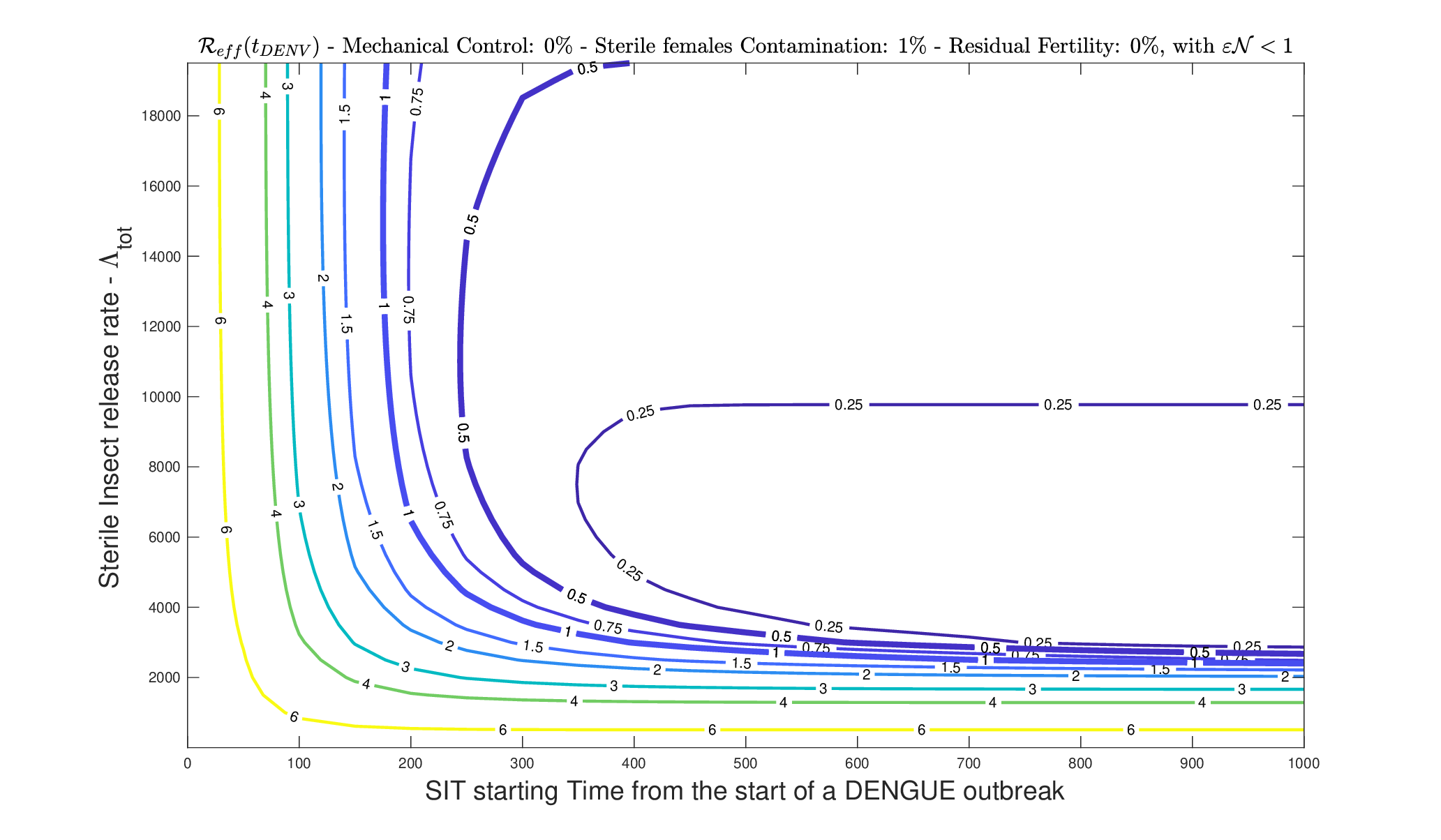} 
   \caption{$\R_{eff}(t_I)$ vs the starting time and the level of the control with $1\%$ of contamination by sterile females, $0\%$ of residual fertility, and without Mechanical control}
    \label{fig:2c}
\end{figure}
\begin{figure}[h!]
\includegraphics[width=0.9\linewidth]{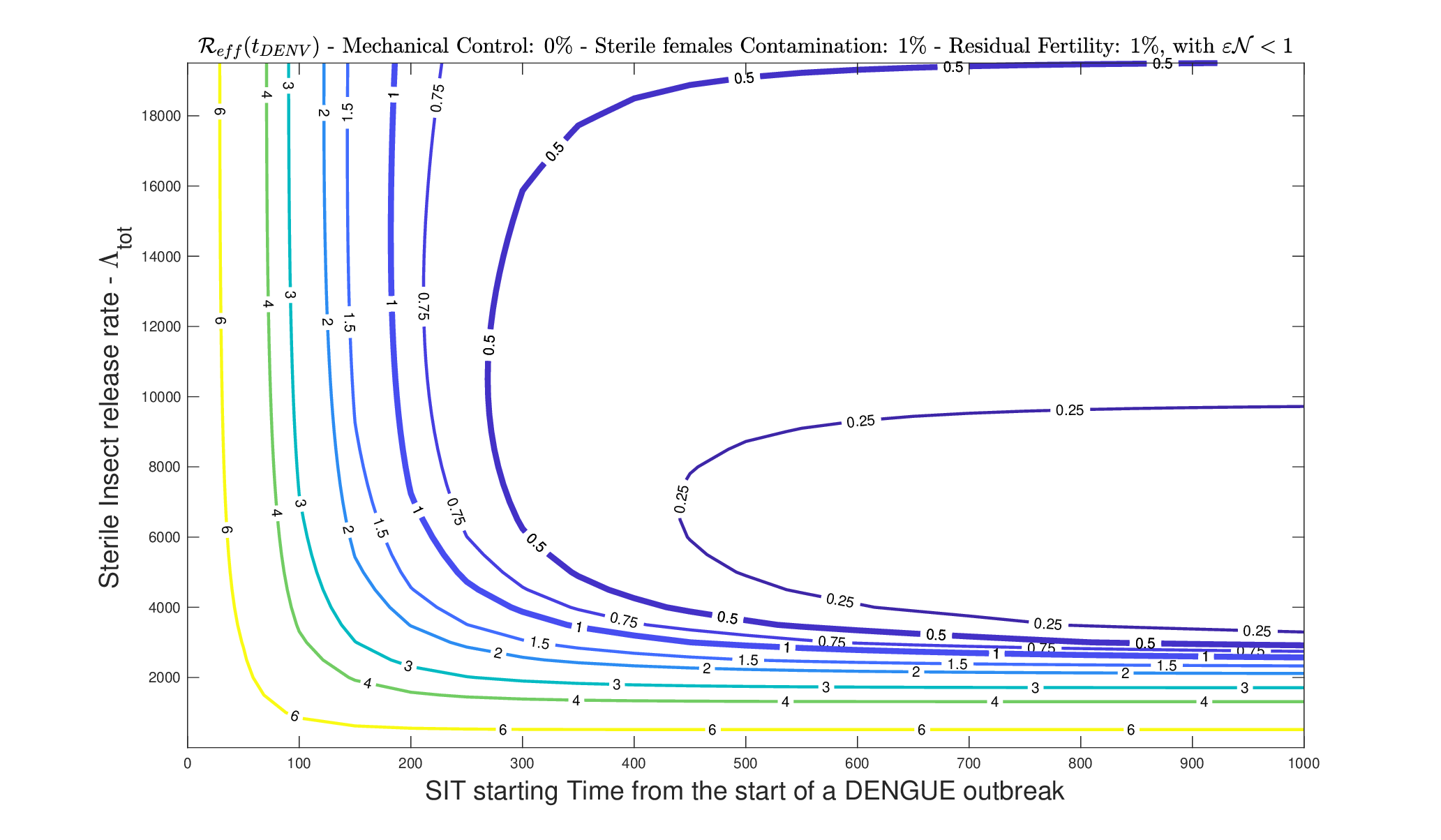} 
   \caption{$\R_{eff}(t_I)$ vs the starting time and the level of the control with $1\%$ of contamination by sterile females, $1\%$ of residual fertility, and without Mechanical control}
    \label{fig:2cd}
\end{figure}
\begin{figure}[h!]
\includegraphics[width=0.9\linewidth]{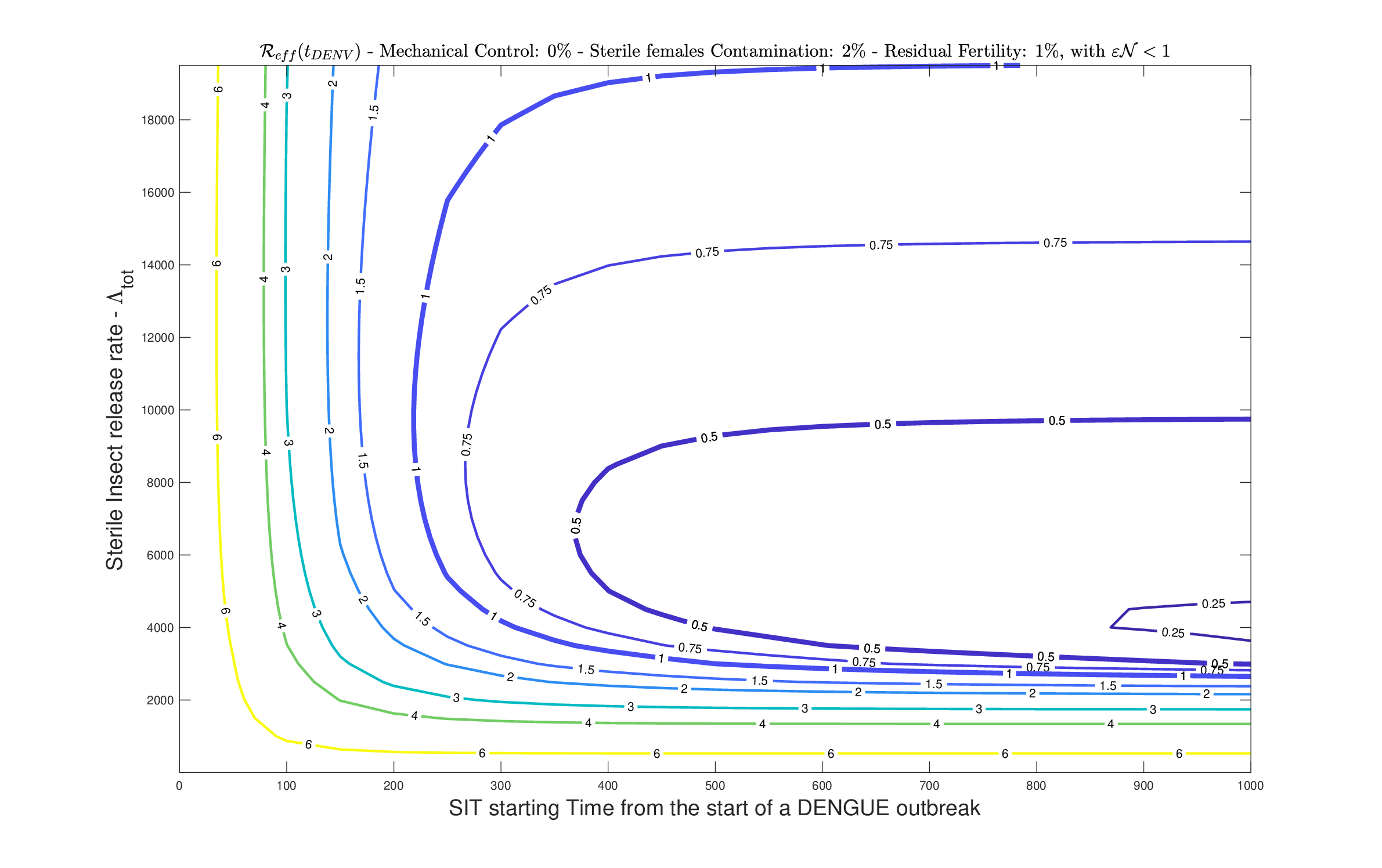} 
   \caption{$\R_{eff}(t_I)$ vs the starting time and the level of the control with $2\%$ of contamination by sterile females, $1\%$ of residual fertility, and without Mechanical control}
    \label{fig:2cde}
\end{figure}

\begin{figure}[h!]
\includegraphics[width=0.9\linewidth]{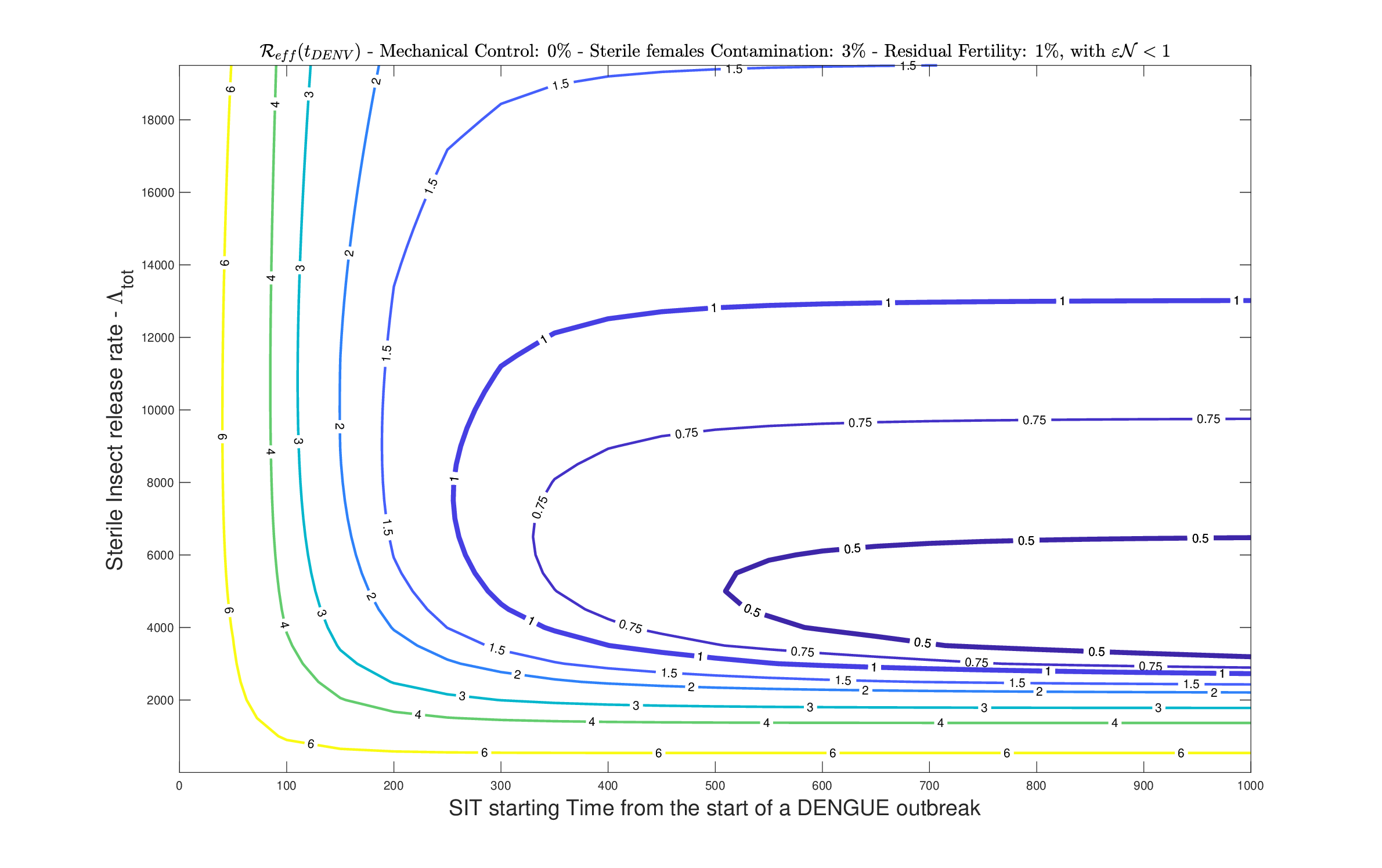} 
    \caption{$\R_{eff}(t_I)$ vs the starting time and the level of the control with $3\%$ of contamination by sterile females, $1\%$ of residual fertility, and without Mechanical control}
    \label{fig:4aa}
\end{figure}

\begin{figure}[h!]
\includegraphics[width=0.9\linewidth]{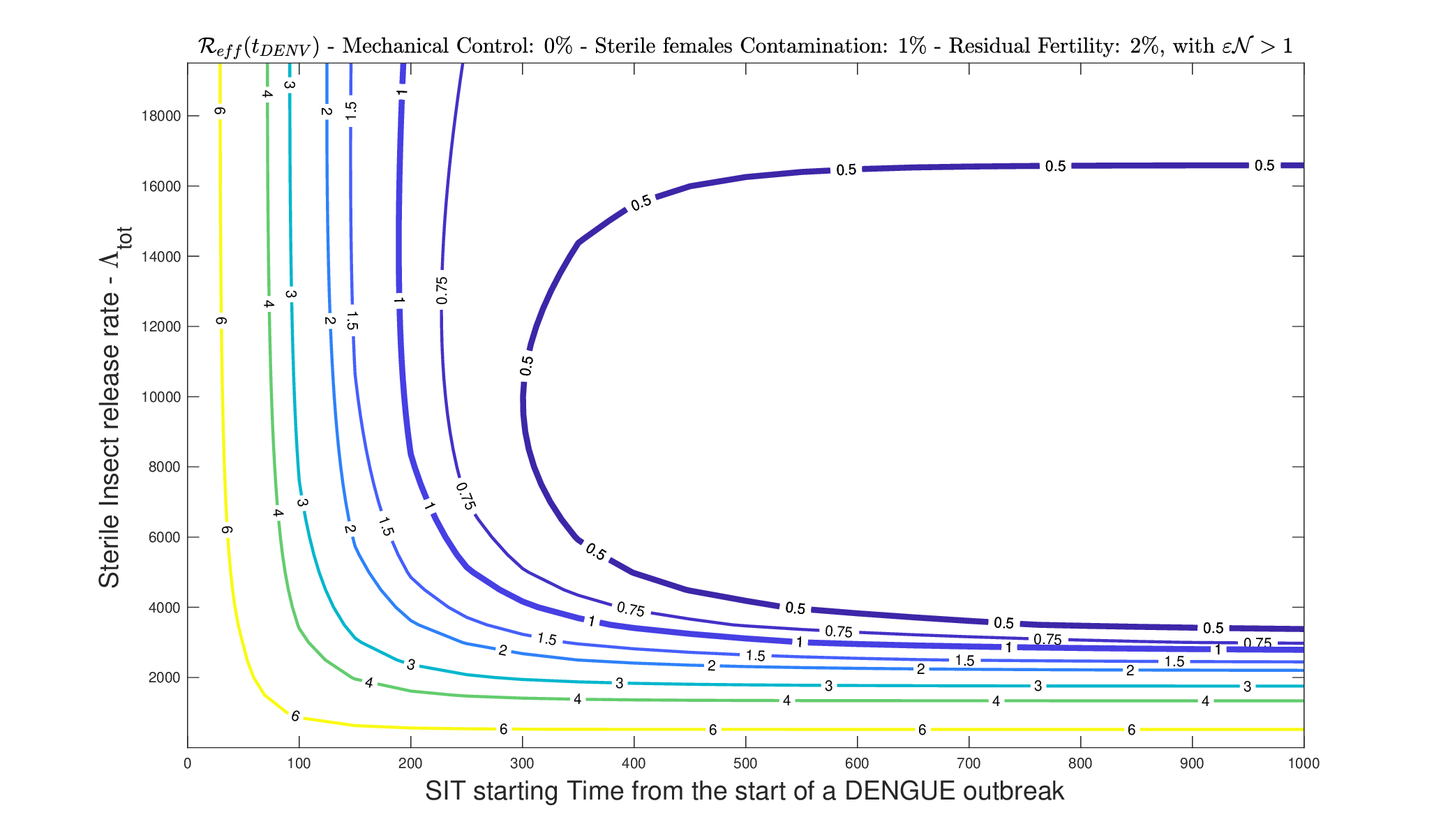} 
    \caption{$\R_{eff}(t_I)$ vs the starting time and the level of the control with $1\%$ of contamination by sterile females, $2\%$ of residual fertility, and without Mechanical control}
    \label{fig:2d}
\end{figure}

We provide simulations with several combination of values for $\epsilon_F$ from $0\%$ to $3\%$, and $\varepsilon$ from $0\%$ to $2\%$.

Since $\Lambda_{tot}$ varies from $0$ to $20 000$, then according to $\epsilon_F$, $\Lambda_F^{crit}$ varies from $0$ to $400$, when $\epsilon_F=0.01$, from $0$ to $800$, , when $\epsilon_F=0.02$, from $0$ to $1200$, when , when $\epsilon_F=0.03$. Thus, in the forthcoming simulations, for sufficiently large values of $\Lambda_{tot}$, we will have $\Lambda_F>\Lambda_F^{cri}$.

In Tables \ref{Table-threshold_1} and \ref{Table-threshold_2}, we illustrate some of the cases given in Tables \ref{table:recapitulatif} and \ref{table: recapitulatif2}. Clearly, when $\N \varepsilon >1$ (see Table \ref{Table-threshold_2}), we highlight the fact that it is more difficult to control the epidemiological risk, even with a release rate just above the critical threshold, and such that $\varepsilon_F \Lambda_{tot}<<\Lambda_F^{crit}$. In contrary, when $\N \varepsilon<1$, epidemiological control is easier to reach even with a substantial increase of the contamination by sterile females: see Table \ref{Table-threshold_1}. These results are also supported by the forthcoming simulations. 

\begin{table}[h!]
  \centering
  \begin{tabular}{|c|c|c|c|c|c|c|}
\hline
  $\epsilon_F$ & 0 & 0.01 & 0.02 & 0.03 & 0.05\\
  \hline \hline
    $(1-\epsilon_F)\Lambda_{tot}$ & 3700 & 3663 & 3626 & 3589 & 3515\\
  \hline
      $\epsilon_F\Lambda_{tot}$ & 0 & 37 & 74 & 111 & 185\\
    \hline
    $\R_{0,\N \varepsilon<1}^2$ & 3.51 & 3.314 & 3.143 & 2.99 & 2.72\\
    \hline 
    $\R_{0,SIT_c,W}^2$ & 0 & 0 & 0.422 & 0.527 & 0.701\\
     \hline
    $\R_{0,SIT_c,S}^2$ & 0 & 0.095 & 0.189 & 0.284  & 0.406\\
     \hline
    $\R_{0,SIT_c}^2$ & 0 & 0.095 & 0.61 & 0.81 & 1.17\\
    
\hline
\end{tabular}
\caption{Threshold values to lower the epidemiological risk for DENV when $\varepsilon=0.01$, such that $\N \varepsilon<1$, $\Lambda_M^{crit}=3653$, $\Lambda_F^{crit}=391$, and $\R_0^2>\R_{0,\N \varepsilon<1}^2$.}
\label{Table-threshold_1}
\end{table}
\begin{table}[h!]
  \centering
  \begin{tabular}{|c|c|c|c|c|c|}
\hline
  $\epsilon_F$ & 0 & 0.01 & 0.02 & 0.03\\
  \hline
  \hline
 $(1-\epsilon_F)\Lambda_{tot}$ & 3700 & 3663 & 3626 & 3515\\
  \hline
      $\epsilon_F\Lambda_{tot}$ & 0 & 37 & 74 & 111\\
  \hline
   $\R_{0,\N \varepsilon>1}^2$ & 116.7 & 105.6 & 94.58 & 83.53 \\
    \hline
    $\Lambda_M^{crit,\sharp}$ & 2869 & 2971 & 3074 & 3176 \\
      \hline
  $\Lambda_{M,\N \varepsilon>1}^{crit}$ & 3638 & 3718 & 3806 & 3905  \\
    \hline
    $\R_{0,SIT_c,W}^2$ & 0.925 & 0.0.969 & 1.01 & 1.06  \\
     \hline
    $\R_{0,SIT_c,S}^2$ & 0 & 0.095 & 0.19 & 0.28  \\
     \hline
    $\R_{0,SIT_c}^2$ & 0.925 & 1.064 & 1.20 & 1.34 \\
    
\hline
\end{tabular}
\caption{Threshold values to lower the epidemiological risk for DENV when $\varepsilon=0.02$, such that $\N \varepsilon>1$, $\Lambda_F^{crit}=391$, and $\R_0^2<\R_{0,\N\varepsilon>1}^2$.} \label{Table-threshold_2}.
\end{table}

In Figs. \ref{fig:2a}, and \ref{fig:2b}, page \pageref{fig:2a}, we consider the case where there is no contamination by sterile females, with $\varepsilon$ such that $\varepsilon \N <1$ and $\varepsilon \N >1$, that is where $\varepsilon=0$ and $\varepsilon=0.02$. Roughly speaking, it is easy to observe that residual fertility has less impact on the rate needed to decay  $\R_{eff}$ below $0.5$. When $\varepsilon \N>1$, it is not possible to lower the wild population under any given small threshold, to reduce the nuisance for instance, but it is still possible to reduce the epidemiological risk, at least when no female contamination occurs.

\begin{figure}[h!]
\includegraphics[width=0.9\linewidth]{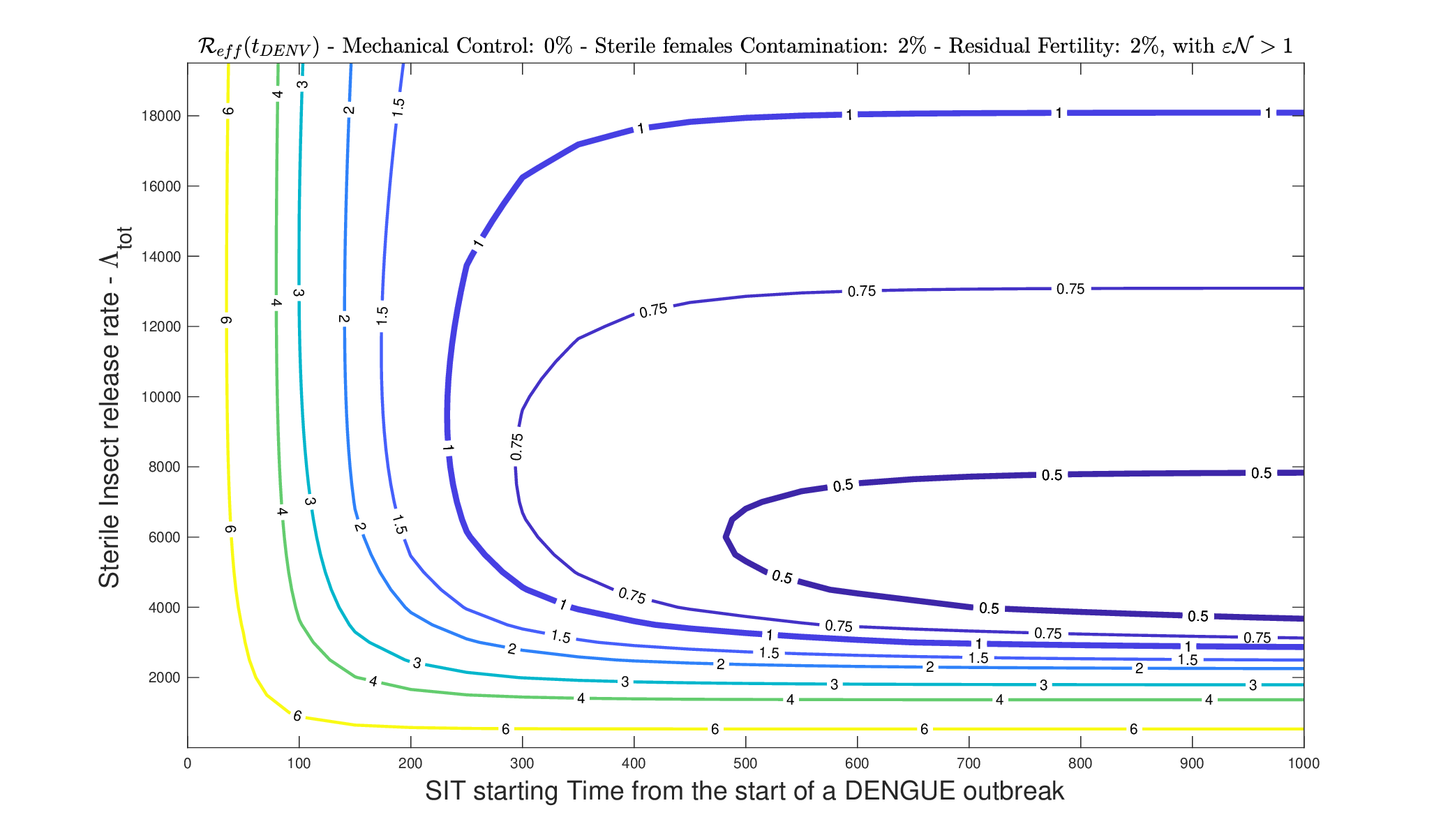} 
    \caption{$\R_{eff}(t_I)$ vs the starting time and the level of the control with $2\%$ of contamination by sterile females, $2\%$ of residual fertility, and without Mechanical control}
    \label{fig:3a}
\end{figure}

\begin{figure}[h!]
\includegraphics[width=0.9\linewidth]{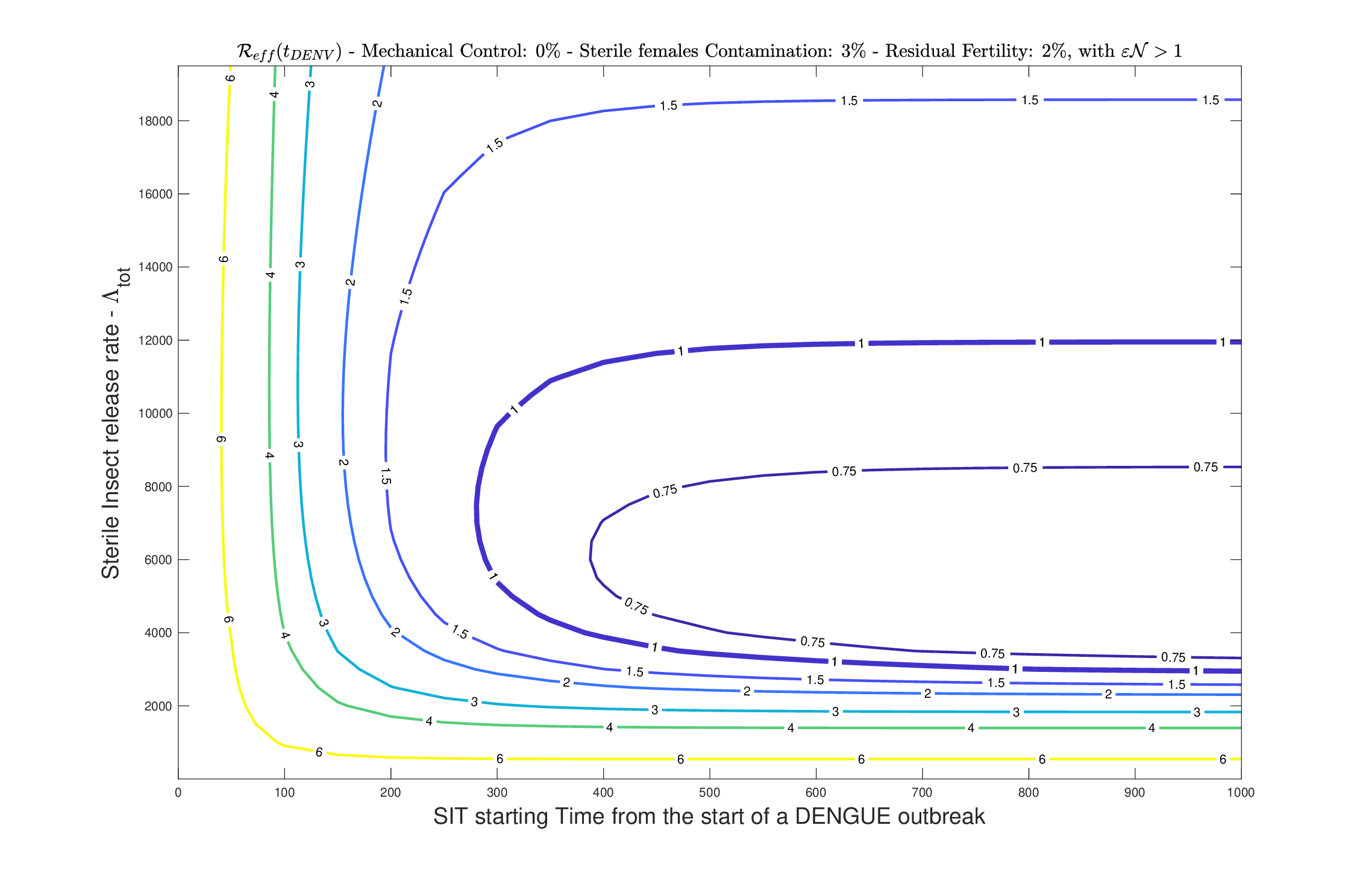} 
    \caption{$\R_{eff}(t_I)$ vs the starting time and the level of the control with $3\%$ of contamination by sterile females, $2\%$ of residual fertility, and without Mechanical control}
    \label{fig:4a}
\end{figure}

From Fig. \ref{fig:2c}, page \pageref{fig:2c}, to Fig. \ref{fig:4a}, page \pageref{fig:4a}, we consider contamination by sterile females with a residual fertility varying from $1\%$ to $2\%$ in order to consider both cases $\N \varepsilon<1$ and $\N \varepsilon>1$. It is interesting to notice that the shape of the level sets change according to $\epsilon_F$, such that when $\epsilon_F$ increases, the area where $\R_{eff}<0.5$ decays. In fact, when $\epsilon_F$ is large, say $2\%$ or $3\%$, then very massive releases are such that $\epsilon_F \Lambda_{tot}>\Lambda_F^{crit}$ which implies $\R_{0,TDFE}^2>1$ and $\R_{eff}>1$: see Figs. \ref{fig:4aa}, \ref{fig:3a}, and \ref{fig:4a}. This simulation clearly shows that increasing the release rate is not the right response, whatever if $\varepsilon \N$ is less or greater than $1$, when SIT is used to decay the epidemiological risk.  Clearly, as long as the female contamination is large, increasing the release rate will take the sterile females close to the release rate threshold, $\Lambda_F^{crit}$, such that $\R_{eff}>1$. Note also, that our simulations show In that an optimal release rate exists for a given, sufficiently large, SIT starting time.

\begin{figure}[h!]
\includegraphics[width=0.9\linewidth]{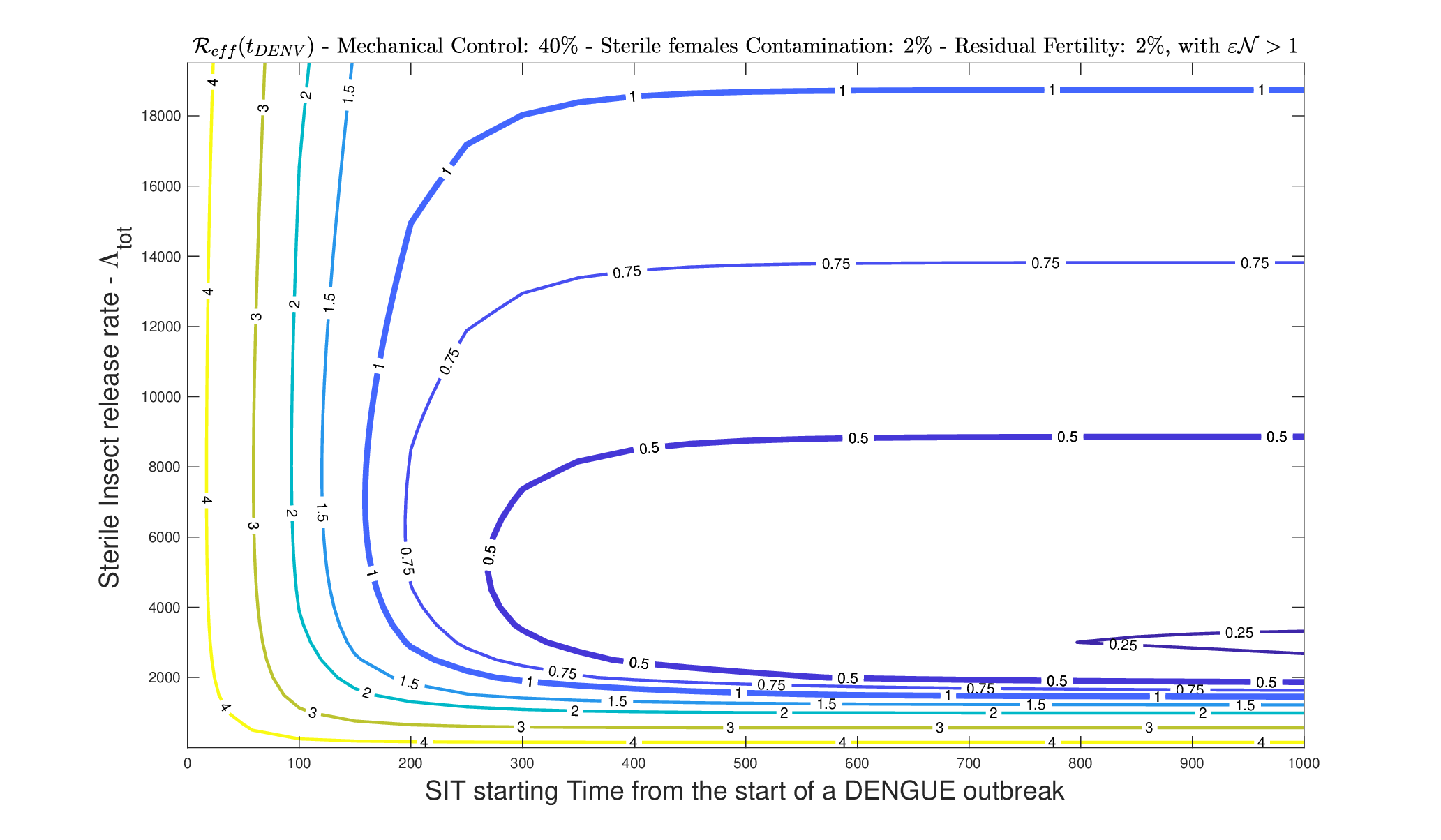} 
    \caption{$\R_{eff}(t_I)$ vs the starting time and the level of the control with $2\%$ of contamination by sterile females, $2\%$ of residual fertility, and $40\%$ of Mechanical control}
    \label{fig:3b}
\end{figure}

\begin{figure}[h!]
\includegraphics[width=0.9\linewidth]{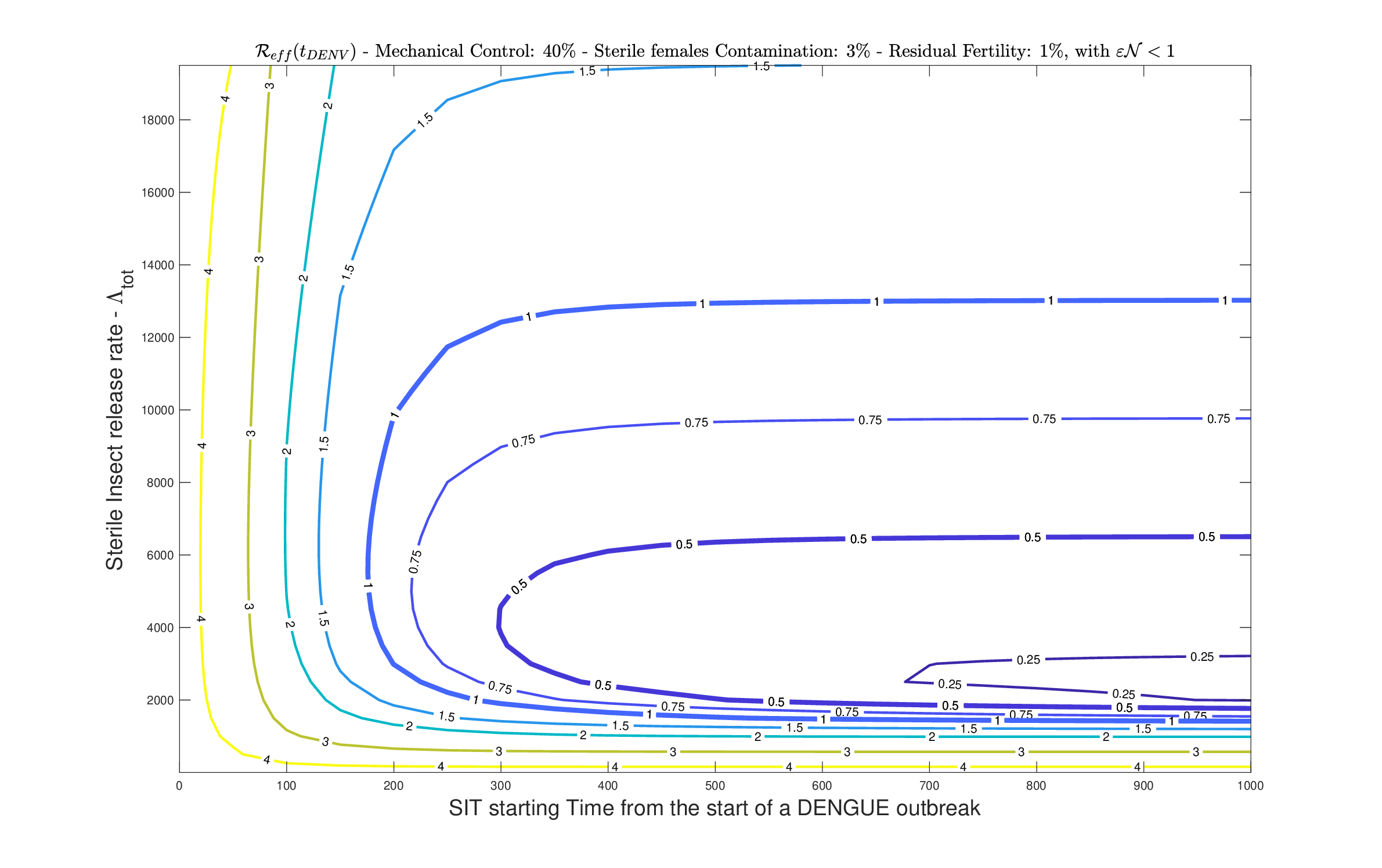} 
    \caption{$\R_{eff}(t_I)$ vs the starting time and the level of the control with $3\%$ of contamination by sterile females, $1\%$ of residual fertility, and $40\%$ of Mechanical control}
    \label{fig:4bb}
\end{figure}

\begin{figure}[h!]
\includegraphics[width=0.9\linewidth]{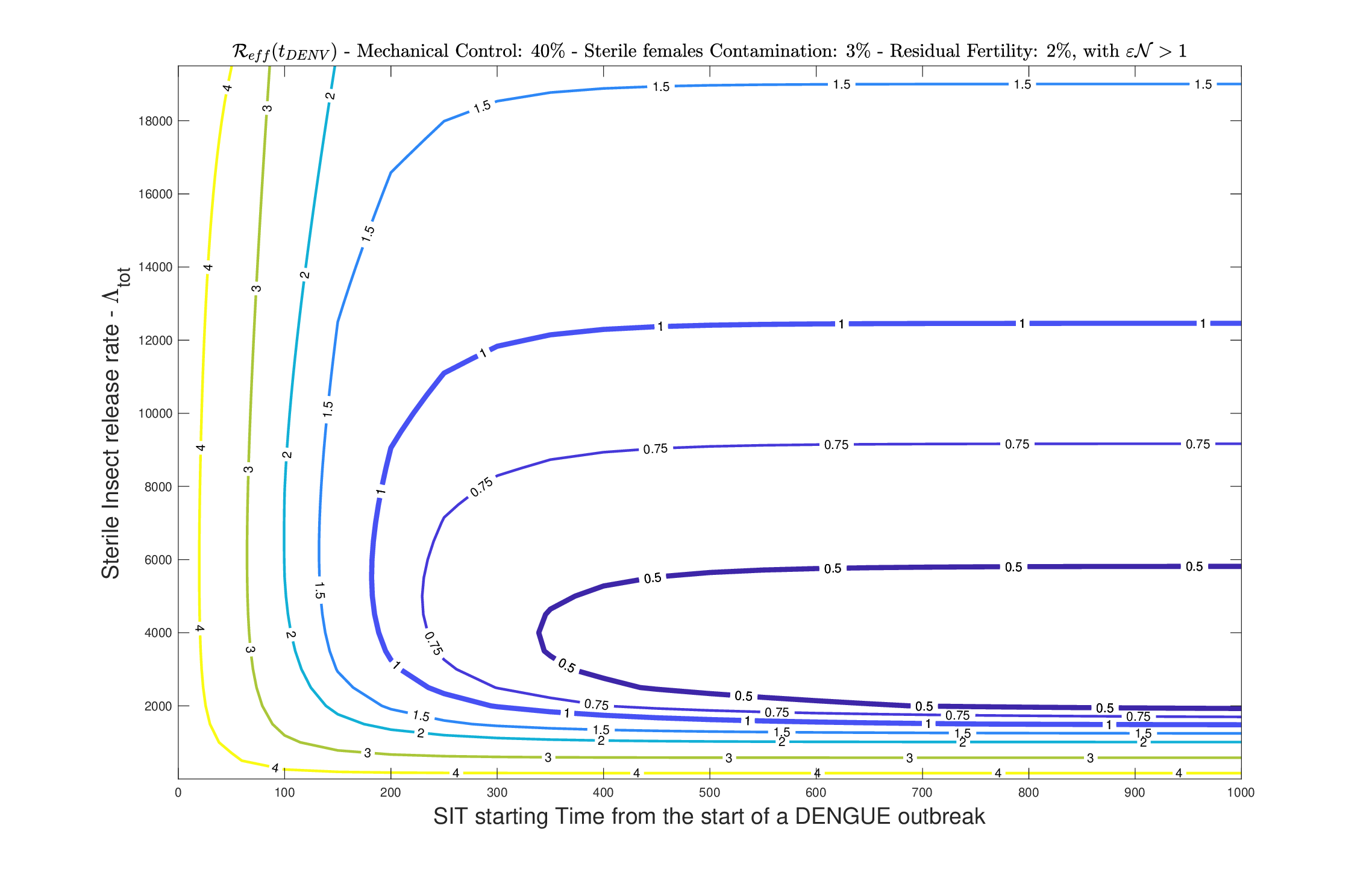} 
    \caption{$\R_{eff}(t_I)$ vs the starting time and the level of the control with $3\%$ of contamination by sterile females, $2\%$ of residual fertility, and $40\%$ of Mechanical control}
    \label{fig:4b}
\end{figure}

Mechanical control is clearly beneficial to reduce the time needed to decay $\R_{eff}$ below $0.5$ and also the (optimal) release rate: compare Figs. \ref{fig:3a} and \ref{fig:3b}, page \pageref{fig:3b}, where the time needed to reach $0.5$ for $\R_{eff}$ decay from $500$ days, for $\Lambda_{opt}\approx 6000$, to, only $300$ days with $\Lambda_{opt} \approx 4000$, to reduce $\R_{eff}$ before DENV starts to circulate. Compare also Figs. \ref{fig:4aa} and \ref{fig:4a} with Figs. \ref{fig:4bb} and \ref{fig:4b}.

In fact, when $\N \varepsilon>1$, serious problem occurs when contamination by sterile females increases, without mechanical control: see Fig. \ref{fig:4a}, page \pageref{fig:4a}. As seen, it is no more possible to decay $\R_{eff}$ below $0.5$ and, as explained before, very massive release can be such that $\R_{eff}>1$. In that case, SIT cannot be used to control the epidemiological risk, at least without mechanical control. In fig. \ref{fig:4b}, page \pageref{fig:4b}, mechanical control allows to lower the time needed to decay $\R_{eff}$ but does not really increase the maximal release rate such that $\R_{eff}<1$.

Altogether, our numerical simulations, that the first parameter to lower is $\varepsilon$, the residual fertility. However, even with a low residual fertility, say $1\%$, contamination by sterile females should be contained: compare Fig. \ref{fig:2cd}, page \pageref{fig:2cd}, with Fig. \ref{fig:2cde}, page \pageref{fig:2cde}.

\section{Conclusion \label{conclusion}}

Conducting SIT programs in the field is a very complex and difficult task. However, before reaching field releases and in order to be successful, several steps have to be checked in laboratory and in semi-field, before and during field releases. In fact, it is better to find and solve issues before starting field releases: to this aim control quality is an essential process within SIT programs. However, SIT programs against mosquitoes can fail, and this is in general due to a combination of several factors, among them residual male fertility and contamination by sterile females that seem not to be always studied as deep as they should be. Indeed, sometimes (numerical) upper bound values are given for these parameters but they do not rely on biological parameters related to the targeted vectors nor on epidemiological parameters when epidemiological control is the main objective. We aim to fill this gap.

Thus, using modelling and mathematical analysis, we provide threshold parameters for residual male fertility and contamination by sterile females. We also show that these thresholds impose constraints on SIT programs to be met. If not, then, the risk of SIT failure is high. 

Our results could be used and helpful for field experts to estimate the risk of SIT failures and, thus, to target the main parameters to improve before field releases and to follow carefully along the SIT process.


Theoretically, we show that while residual fertility can be an issue to control the wild population, i.e. to lower it under a given threshold, to reduce the nuisance, it is not when it comes to control the epidemiological risk. 
In other words, when $\varepsilon \N<1$, both nuisance reduction and epidemiological risk reduction are feasible as long as the sterile female contamination is low, that is $\epsilon \Lambda_{tot}<\Lambda_F^{crit}$. While, when $\varepsilon \N>1$, only epidemiological risk reduction is feasible but under rather severe constraints, that is $\epsilon \Lambda_{tot}<\Lambda_F^{crit}$ and $\R_0^2<\R_{0,\N\varepsilon>1}^2$, with releases that are sufficiently massive. 

In fact, once $\varepsilon \N<1$ is not met, we strongly encourage the SIT program to solve this issue before going further.

Finally, in several SIT reports/manuals  or SIT papers \cite{Iyaloo2020b}, a percentage is given for the maximal contamination by sterile females. We show that this percentage is useless since the maximal amount of sterile females allowed to be released will depend on the size of the total release. Indeed, you don't release the same amount of sterile females when you consider $1\%$ of $10 000$ or $1\%$ of $20 000$ sterile insects: for the first case, $\Lambda_F<\Lambda_F^{crit}$, while in the second case, $\Lambda_F>\Lambda^{crit}_F$, such that the dynamics of the whole system is completely different and so is the impact of SIT.

To conclude, our study shows that both contamination by sterile females, $\epsilon_F\Lambda_{tot}$, and residual male fertility, $\varepsilon$, matter in the efficiency of SIT. We provide upper bounds for these values that guarantee the efficiency of SIT, both for nuisance and epidemiological risk reduction. 

Of course, several improvements are possible, like considering impulsive releases, like in \cite{DumontYatat2022}. In addition other control quality tests could be taken into account in future SIT models in order to provide more realistic results, and eventually, when possible, to consider variable parameters, like in \cite{Duprez2022} to take into account temporal and spatial variation of the environmental parameters that can affect the dynamics of the vectors and thus its control. Last, migration could be also taken into account \cite{Bliman2022}.

\paragraph{Acknowledgments:} YD is (partially) supported by the DST/NRF SARChI Chair in Mathematical Models and Methods in Biosciences and Bioengineering at the University of Pretoria, South Africa (Grant 82770). YD acknowledges the support of the Conseil R\'egional de la R\'eunion (France), the Conseil D\'epartemental de la R\'eunion (France), the European Agricultural Fund for Rural Development (EAFRD) and the Centre de Coop\'eration Internationale en Recherche Agronomique pour le D\'eveloppement (CIRAD), France.

\bibliographystyle{plain}
\bibliography{main_DUMONT_YATAT_SIT_residual_fertility_sterileFemale}

\appendix

\section{A useful result on monotone systems}\label{AppendixA}
Let us consider an $n$ dimensional autonomous differential system:
\begin{equation}\label{ode}
\dfrac{dx}{dt}=f(x)
\end{equation}
where $f$ is a given vector function, i.e., $f = (f_{i})_{i=1,...,n}$, with $f_{i} : \mathbb{R}^{n} \mapsto \mathbb{R}$. System \eqref{ode} is called cooperative if for every $i, j \in \{1, 2, ..., n\}$ such that $i \ne j$, the function $f_{i}(x_{1}, ..., x_{n})$ is monotone increasing with respect to $x_{j}$. For cooperative system, the global asymptotic stability of an equilibrium can be studied by the following theorem, see also \cite{Anguelov2012TIS}:
\begin{theorem}\label{theo-monotone}
	Assume that system \eqref{ode} is a cooperative system. Let $a$, $b \in \Omega \subseteq \mathbb{R}^{n}$ such that $a < b$, $\left[a,b\right] \subseteq \Omega$ and $f(b) \leq 0 \leq f(a)$; where $\left[a, b\right] = \{x \in \mathbb{R}^{n}: a \leq x \leq b \}$. Then \eqref{ode} defines a (positive) dynamical system on $\left[a, b\right]$. Moreover, if $\left[a, b\right]$ contains a unique equilibrium $p$, then $p$ is globally asymptotically stable on $\left[a, b\right]$.
\end{theorem}

\section{Proof of Theorem \ref{theo-stabilite-equilibre-sit}}\label{AppendixB}
For reader convenience, we recall that $M_S^*=\dfrac{(1-\epsilon_F) \Lambda_{tot}}{\mu_{M_S}}$.
\begin{enumerate}
	\item Assume that $\mathcal{N}\varepsilon < 1$. By computing the eigenvalues of the Jacobian matrix of system \eqref{ODE-entomo} at the elimination equilibrium $E_{0}$ it is straightforward to obtain $E_{0}$ is locally asymptotically stable when $\mathcal{N}\varepsilon < 1$ while it is unstable when $\mathcal{N}\varepsilon >1$.
	\begin{itemize}
		\item Let us set $X = (A, M, F) \in \mathbb{R}_{+}^{4}$ and $f((1-\epsilon_F) \Lambda_{tot}, X)$ the right hand side of system \eqref{ODE-entomo}. For $(1-\epsilon_F) \Lambda_{tot} > 0$, we have that $f((1-\epsilon_F) \Lambda_{tot}, X) \leq f(0, X)$. Note that for $(1-\epsilon_F) \Lambda_{tot} = 0$, we recover \cite[system (1)]{Anguelov2020}.  If $(1-\epsilon_F) \Lambda_{tot} > \Lambda_M^{crit}$, then system \eqref{ODE-entomo} admits a unique equilibrium which is $E_{0}$. Using \cite[Theorem 3, point (1)]{Anguelov2020}, we deduce that $E_{0}$ is globally asymptotically stable in $\mathbb{R}_{+}^{3}$.
		\item The proof of points $(b)$ and $(c)$ is done in the same way like the proof of \cite[Theorem 3, points (2) \& (3)]{Anguelov2020}.
	\end{itemize}
	\item Assume that $\mathcal{N}\varepsilon > 1$. Then, the elimination equilibrium $E_{0}$ is unstable. Moreover, the inequality
	\begin{equation}\label{a1}
		\dfrac{\gamma+\mu_{1}+\mu_{2}A}{\gamma+\mu_{1}} > 4\mathcal{N}
	\end{equation}	
	holds for all sufficiently large $A$. Let $n > 0$ and let $A_{n}$ be so large that in addition to \eqref{a1} the following inequalities also hold:
	\begin{equation}\label{a2}
		\begin{array}{l}			
			A_{n} \geq n, \\
			M_{n} := \dfrac{2\left(1-r\right)\gamma}{\mu_{M}}A_{n} \geq n,\\
			F_{W,S_n} := \dfrac{\left(\gamma+\mu_{1}+\mu_{2}A_{n}\right)A_{n}}{2\phi} \geq n.\\					
		\end{array}
	\end{equation}
	Let $b_{n} = (A_{n}, M_{n}, F_{W,S_n})$ and $f$ be the right hand side of \eqref{ODE-entomo}. Then
	\begin{equation}\label{a3}
		f((1-\epsilon_F) \Lambda_{tot}, b_{n}) \leq f(0, b_{n}) =\begin{pmatrix}
			-\phi F_{W,S_n}\\
				-\dfrac{1}{2}\mu_MM_n\\
			r\gamma A_n\left(1-\dfrac{\gamma+\mu_1+\mu_2A_n}{4\mathcal{N}(\gamma+\mu_1)}\right)\\		
		\end{pmatrix} < 0_{\mathbb{R}^3}.
	\end{equation}
	Similarly, for an arbitrary $\delta > 0$, let  $a_{\delta} = (A_{\delta}, F_{W,S_\delta}, M_{\delta})$ with 
	\begin{equation}\label{a4}
		\begin{array}{l}
			M_{\delta} = \dfrac{(1-\epsilon_F) \Lambda_{tot}}{\mu_{M_S}(\alpha_{+}+\delta)} < M_{\dag}, \\
			A_{\delta} = \dfrac{\mu_{M}}{\left(1-r\right)\gamma}M_{\delta} < A_{\dag}, \\
			F_{W,S_\delta} = \dfrac{\left(\gamma+\mu_{1}+\mu_{2}A_{\delta}\right)}{\phi}A_{\delta} < F_{W,S_\dag}.
		\end{array}
	\end{equation}
	We also have that
	\begin{align}\label{a5}
		\dfrac{M_{\delta}+\varepsilon M_{S}^{\ast}}{M_{\delta}+ M_{S}^{\ast}}r\gamma A_{\delta} - \mu_{S}F_{W,S_\delta} &= \dfrac{\mu_S(\gamma+\mu_{A,1})A_{\delta}}{\phi}\left(\dfrac{1+\varepsilon\left(\alpha_{+}+\delta\right)}{1+\alpha_{+}+\delta}\mathcal{N}-1-\dfrac{\mathcal{Q}_S}{\alpha_{+}+\delta}\right)\\
		&= \dfrac{\mu_S(\gamma+\mu_{A,1})A_{\delta}}{\phi}\dfrac{\left(\mathcal{N}\varepsilon-1\right)\delta\alpha_{+}+\left(\mathcal{N}\varepsilon-1\right)\delta^{2}+\dfrac{\delta\mathcal{Q}_S}{\alpha_{+}}}{\left(1+\alpha_{+}+\delta\right)\left(\alpha_{+}+\delta\right)}\\
		&> 0.
	\end{align}
	Thus, it is straightforward to obtain that
	\begin{equation}\label{a6}
		f((1-\epsilon_F) \Lambda_{tot}, a_{\delta}) =\begin{pmatrix}
			0\\
			0\\
			\dfrac{M_{\delta}+\varepsilon M_{S}^{\ast}}{M_{\delta}+ M_{S}^{\ast}}r\gamma A_{\delta} - \mu_{S}F_{W,S_\delta}\\
		\end{pmatrix} > 0_{\mathbb{R}^3}.
	\end{equation}
	Applying Theorem \ref{theo-monotone} with $a = a_{\delta}$ and $b = b_{n}$, we obtain that for $n$  sufficiently large, system \eqref{ODE-entomo} defines a dynamical system on $\left[a_{\delta}, b_{n}\right]$ and that $E_{\dag}$ is globally asymptotically stable on $\left[a_{\delta}, b_{n}\right]$. Since $b_{n}$ can be selected to be larger than any point in $\mathbb{R}_{+}^{3}$ and $a_{\delta}$ can be selected to be lower than any point in $\mathbb{R}_{+}^{3}-\{0_{\mathbb{R}^3}\}$, this implies that $E_{\dag}$ is globally asymptotically stable in $\mathbb{R}_{+}^{3}-\{0_{\mathbb{R}^3}\}$.
	\item Assume that $\mathcal{N}\varepsilon = 1$.
	\begin{enumerate}
		\item If $(1-\epsilon_F) \Lambda_{tot} \geq \Lambda_{M, \sharp}^{crit}$, then $E_{0}=(0, 0, 0)$ is the only equilibrium of system \eqref{ODE-entomo}. Based on \eqref{a3} and Theorem \ref{theo-monotone}, we obtain that for $n$ sufficiently large, system \eqref{ODE-entomo} defines a dynamical system on $\left[E_{0}, b_{n}\right]$. Since $b_{n}$ can be selected to be larger than any point in $\mathbb{R}_{+}^{3}$, this implies that $E_{0}$ is globally asymptotically stable on $\mathbb{R}_{+}^{3}$. 
		\item $(1-\epsilon_F) \Lambda_{tot} \in \left(0, \Lambda_{M, \sharp}^{crit}\right)$, then we proceed as in point 2 by replacing $E_{\dag}$ by $E_{\sharp}$. Hence, we obtain that $E_{\sharp}$ is globally asymptotically stable. Then, the elimination equilibrium $E_{0}$ is unstable and the coexistence equilibrium $E_{\sharp}$ is globally asymptotically stable in $\mathbb{R}_{+}^{3}-\{0_{\mathbb{R}^3}\}$.\\
		This ends the proof.
	\end{enumerate}
\end{enumerate}

\section{Proofs of Propositions \ref{EE_existence}-\ref{EE-existence_general}: existence of endemic equilibria}\label{AppendixC}

First, it is interesting to check after an artificial endemic equilibrium, without wild insect, called WIFE, Wild Insect-Free boundary Equilibrium. This is a particular case, but it can exists. To find it it suffices to solve
$$
\left\{ %
\begin{array}{lcl}
0 & = & \mu_{h}N_{h}-B\beta_{mh}{\displaystyle \frac{S_{I}}{N_{h}}S_{h}-\mu_{h}S_{h},}\\
{\displaystyle 0} & = & B\beta_{mh}{\displaystyle \frac{S_{I}}{N_{h}}S_{h}-\nu_hI_{h}-\mu_{h}I_{h},}\\
0 & = & \nu_hI_{h}-\mu_{h}R_{h},
\end{array}\right.
$$
and
$$
\left\{ \begin{array}{ccl}
0 & = & \epsilon_{F}\Lambda_{tot}-B\beta_{hm}{\displaystyle \frac{I_{h}}{N_{h}}S_{S}-\mu_{S}S_{S},}\\
{\displaystyle 0} & = & B\beta_{hm}{\displaystyle \frac{I_{h}}{N_{h}}S_{S}-(\nu_{m}+\mu_{S})S_{E},}\\
{\displaystyle 0} & = & \nu_{m}S_{E}-\mu_{S}S_{I}.
\end{array}\right.
$$
Straightforward computations show that
$$
S_{I}^{\#}=\dfrac{\nu_{m}}{\mu_{I}\left(\nu_{m}+\mu_{S}\right)}\left(1-\dfrac{{\displaystyle \mu_{S}}}{\mu_{S}+B\beta_{hm}\dfrac{\mu_{h}}{\mu_{h}+\nu_h}\left(1-\dfrac{S_{h}^{\#}}{N_{h}}\right)}\right)\epsilon_{F}\Lambda_{tot},
$$
such that $0<\dfrac{S_h^*}{N_h}\leq 1$ is a positive root of the second order equation
$$
\left(\mu_{S}+B\beta_{hm}\dfrac{\mu_{h}}{\mu_{h}+\nu_h}\right)-\left(\mu_{S}+2B\beta_{hm}\dfrac{\mu_{h}}{\mu_{h}+\nu_h}+\mu_{S}\dfrac{\epsilon_{F}\Lambda_{tot}}{\Lambda_{F}^{crit}}\right)X+\left(B\beta_{hm}\dfrac{\mu_{h}}{\mu_{h}+\nu_h}+\mu_{S}\dfrac{\epsilon_{F}\Lambda_{tot}}{\Lambda_{F}^{crit}}\right)X^{2}=0.
$$
Assuming $\epsilon_{F}\Lambda_{tot}<\Lambda_{F}^{crit}$, we derive
$\dfrac{S_{h}^{\#}}{N_{h}}=1$, the TDFE equilibrium, and 
$$
\dfrac{S_{h}^{\#}}{N_{h}}=\dfrac{\mu_{S}+B\beta_{hm}\dfrac{\mu_{h}}{\mu_{h}+\nu_h}}{B\beta_{hm}\dfrac{\mu_{h}}{\mu_{h}+\nu_h}+\mu_{S}\dfrac{\epsilon_{F}\Lambda_{tot}}{\Lambda_{F}^{crit}}}>1,
$$
that is not a viable root.
When $\epsilon_{F}\Lambda_{tot}=\Lambda_{F}^{crit}$, we recover $\dfrac{S_{h}^{\#}}{N_{h}}=1$.
Then, assuming $\epsilon_{F}\Lambda_{tot}>\Lambda_{F}^{crit}$ or equivalently $\R_{0,TDFE}^2>1$, a boundary wild insects-free equilibrium $WIFE=(S_h^{\sharp},I_h^{\sharp},R_h^{\sharp},0,0,0,0,0,S_S^{\sharp},S_E^{\sharp},S_I^{\sharp})$ exists such that $$
\dfrac{S_{h}^{\#}}{N_{h}}=\dfrac{\mu_{S}+B\beta_{hm}\dfrac{\mu_{h}}{\mu_{h}+\nu_h}}{B\beta_{hm}\dfrac{\mu_{h}}{\mu_{h}+\nu_h}+\mu_{S}\dfrac{\epsilon_{F}\Lambda_{tot}}{\Lambda_{F}^{crit}}}<1,
$$
$$
    \begin{array}{ccl}
     I_h^\sharp &=& \dfrac{B\beta_{mh}}{\nu_h+\mu_h}\dfrac{S_I^\sharp}{N_h}S_h^\sharp,\\
     S_S^\sharp &=& \dfrac{\epsilon_F \Lambda_{tot}}{\mu_S+\dfrac{B\beta_{mh}}{N_h}I_h^\sharp},\\
     S_E^\sharp &=& \dfrac{\mu_S}{\nu_m}S_I^\sharp,\\
     R_h^\sharp &=& \dfrac{\nu_h}{\mu_h}I_h^\sharp.\\
    \end{array}
$$

The assumption $\mu_I=\mu_S$ is to simplify the forthcoming computations. 
In order to derive existence of a positive endemic equilibrium, such that $I_{h}>0$, $F_{W,I}>0$, and $S_{I}>0$, we solve
\begin{equation}
\label{solve_systH_EE}
\left\{ %
\begin{array}{lcl}
0 & = & \mu_{h}N_{h}-B\beta_{mh}{\displaystyle \frac{F_{W,I}+S_{I}}{N_{h}}S_{h}-\mu_{h}S_{h},}\\
{\displaystyle 0} & = & B\beta_{mh}{\displaystyle \frac{F_{W,I}+S_{I}}{N_{h}}S_{h}-\nu_hI_{h}-\mu_{h}I_{h},}\\
0 & = & \nu_hI_{h}-\mu_{h}R_{h},
\end{array}\right.
\end{equation}
\begin{equation}
\label{solve_systM_EE}
\left\{ %
\begin{array}{lcl}
0 & = & \phi(F_{W,S}+F_{W,E}+F_{W,I})-(\gamma+\mu_{A,1}+\mu_{A,2}A)A,\\
{\displaystyle 0} & = & (1-r)\gamma A-\mu_{M}M,\\
{\displaystyle 0} & = & \dfrac{M+\varepsilon M_{S}^{*}}{M+M_{S}^{*}}r\gamma A-B\beta_{hm}{\displaystyle \frac{I_{h}}{N_{h}}F_{W,S}-\mu_{S}F_{W,S},}\\
0 & = & B\beta_{hm}{\displaystyle \frac{I_{h}}{N_{h}}F_{W,S}-(\nu_{m}+\mu_{S})F_{W,E},}\\
0 & = & \nu_{m}F_{W,E}-\mu_{S}F_{W,I},\\
0 & = & \epsilon_F \Lambda_{tot}+\dfrac{(1-\varepsilon)M_{S}^{*}}{M+M_{S}^{*}}r\gamma A-B\beta_{hm}{\displaystyle \frac{I_{h}}{N_{h}}S_{S}-\mu_{S}S_{S},}\\
{\displaystyle 0} & = & B\beta_{hm}{\displaystyle \frac{I_{h}}{N_{h}}S_{S}-(\nu_{m}+\mu_{S})S_{E},}\\
{\displaystyle 0} & = & \nu_{m}S_{E}-\mu_{S}S_{I}.
\end{array}\right.
\end{equation}
Thanks to \eqref{solve_systM_EE}$_1$, and summing \eqref{solve_systM_EE}$_4$ and \eqref{solve_systM_EE}$_5$ such that 
\[
F_{W,E}+F_{W,I}=\dfrac{B\beta_{hm}}{\mu_{S}}\frac{I_{h}}{N_{h}}{\displaystyle F_{W,S},}
\]
we derive
\[
\phi\left(1+\dfrac{B\beta_{hm}}{\mu_{S}}\frac{I_{h}}{N_{h}}\right)F_{W,S}=(\gamma+\mu_{A,1}+\mu_{A,2}A)A,
\]
From \eqref{solve_systM_EE}$_3$, we have
\[
\dfrac{M+\varepsilon M_{S}^{*}}{M+M_{S}^{*}}r\gamma A=\left(B\beta_{hm}{\displaystyle \frac{I_{h}}{N_{h}}+\mu_{S}}\right)F_{W,S},
\]
that is, since $A>0$,
\[
\N\left(M+\varepsilon M_{S}^{*}\right)=\left(1+\dfrac{\mu_{A,2}}{\gamma+\mu_{A,1}}A\right)\left(M+M_{S}^{*}\right).
\]
Then, using \eqref{solve_systM_EE}$_2$,
\[
M=\dfrac{(1-r)\gamma}{\mu_{M}}A,
\]
we obtain
\[
\N\left(\dfrac{(1-r)\gamma}{\mu_{M}}A+\varepsilon M_{S}^{*}\right)=\left(1+\dfrac{\mu_{A,2}}{\gamma+\mu_{A,1}}A\right)\left(\dfrac{(1-r)\gamma}{\mu_{M}}A+M_{S}^{*}\right),
\]
that is equivalent to the following second order equation
\begin{equation}
\label{roots_A_EE}
\Q\left(\dfrac{(1-r)\gamma}{\mu_{M}}\right)^{2}A^{2}+\dfrac{(1-r)\gamma}{\mu_{M}}\left(\Q M_{S}^{*}-\N\right)A+\left(1-\N\varepsilon\right)M_{S}^{*}=0.
\end{equation}
We calculate
\[
\Delta=\left(\dfrac{(1-r)\gamma}{\mu_{M}}\left(\Q M_{S}^{*}-\N\right)\right)^{2}-4\dfrac{\mu_{A,2}}{\gamma+\mu_{A,1}}\dfrac{(1-r)\gamma}{\mu_{M}}\left(1-\N\varepsilon\right)M_{S}^{*},
\]
that is
\[
\Delta=\left(\dfrac{(1-r)\gamma}{\mu_{M}}\right)^{2}\left(\left(\Q M_{S}^{*}-\N\right)^{2}-4\Q\left(1-\N\varepsilon\right)M_{S}^{*}\right).
\]
When $\N\varepsilon \geq 1$, then $\Delta>0$, and we deduce the existence of one positive root
\[
A_{*}^{EE}=\dfrac{1}{2\Q\dfrac{(1-r)\gamma}{\mu_M}}\left(\N-\Q M_{S}^*+\sqrt{\left(\left(\Q M_{S}^*-\N\right)^{2}+4Q\left(\N \varepsilon-1 \right)M_{S}^*\right)}\right),\quad \mbox{when }\N\varepsilon>1,
\]
for all $M_S^*>0$, or 
\[
A_{*}^{EE}=\dfrac{1}{\Q\dfrac{(1-r)\gamma}{\mu_{M}}}\left(\N-\Q M_{S}^*\right),\quad \mbox{when }\N\varepsilon=1
\]
for all $\Q M_S^*<\N$, that is $(1-\epsilon_F) \Lambda_{tot}<\dfrac{\mu_{M_S}}{\Q}\N$.
 
When $\N\varepsilon<1$, we have now to study the sign of $\Delta$ according to $\Q M_{S}^{*}$, and solve
\[
\left(\Q M_{S}^{*}-\N\right)^{2}-4\left(1-\N\varepsilon\right)\Q M_{S}^{*}=\left(\Q M_{S}^{*}\right)^{2}+\left(\N\right)^{2}-2\left(\N+2\left(1-\N\varepsilon\right)\right)\Q M_{S}^{*}=0,
\]
for which with
\[
\Delta_{S}=4\left(\left(\left(\N+2\left(1-\N\varepsilon\right)\right)\right)^{2}-\left(\N\right)^{2}\right)=16\left(1-\N\varepsilon\right)\left(\N+\left(1-\N\varepsilon\right)\right)>0,
\]
we can deduce the following threshold
\[
\Q M_{S,1}^{*}=\left(\sqrt{\N+\left(1-\N\varepsilon\right)}-\sqrt{1-\N\varepsilon}\right)^{2}>0.
\]
\begin{remark}
Surprisingly, we derive a threshold almost similar to the threshold obtained in \eqref{aux5}. 
\end{remark}
Using the same reasoning than in section \ref{WimwSIT}, we deduce that, since $\Q M_{S}^*<\Q M_{S,1}^{*}$, then there exists two positive roots of \eqref{roots_A_EE}, that is
\[
A_{1}^{EE}=\dfrac{1}{2\Q\dfrac{(1-r)\gamma}{\mu_{M}}}\left(\N-\Q M_{S}^*-\sqrt{\left(\left(\Q M_{S}^*-\N\right)^{2}-4\Q\left(1-\N\varepsilon\right)M_{S}\right)}\right),
\]
\[
A_{2}^{EE}=\dfrac{1}{2\Q\dfrac{(1-r)\gamma}{\mu_{M}}}\left(\N-\Q M_{S}^*+\sqrt{\left(\left(\Q M_{S}^*-\N\right)^{2}-4\Q\left(1-\N\varepsilon\right)M_{S}^*\right)}\right).
\]
Then, we deduce that 
\[
M_{i}^{EE}=\dfrac{(1-r)\gamma}{\mu_{M}}A_{i}^{EE}\qquad i=1,2. 
\]
Thus for a given $A_{i}^{EE}$, we are able to estimate $\dfrac{(1-\varepsilon)M_{S}^{*}}{M+M_{S}^{*}}r\gamma A$
and $\dfrac{M+\varepsilon M_{S}^{*}}{M+M_{S}^{*}}r\gamma A$. In order
to deduce the other variables, some computations are needed. From \eqref{solve_systM_EE}$_7$ and \eqref{solve_systM_EE}$_8$, we have
\[
S_{I}=\dfrac{\nu_{m}}{\mu_{S}}S_{E}=\dfrac{\nu_{m}}{\mu_{S}(\nu_{m}+\mu_{S})}B\beta_{hm}{\displaystyle \frac{I_{h}}{N_{h}}S_{S}}=\dfrac{\dfrac{\nu_{m}}{\mu_{S}(\nu_{m}+\mu_{S})}B\beta_{hm}\frac{I_{h}}{N_{h}}}{\mu_{S}+B\beta_{hm}\frac{I_{h}}{N_{h}}}{\displaystyle \left(\epsilon_F \Lambda_{tot}+\dfrac{(1-\varepsilon)M_{S}^{*}}{M+M_{S}^{*}}r\gamma A\right),}
\]
Similarly, from \eqref{solve_systM_EE}$_3$, \eqref{solve_systM_EE}$_4$, and \eqref{solve_systM_EE}$_5$
\[
F_{W,I}=\dfrac{\nu_{m}}{\mu_{S}}F_{E,I}=\dfrac{\nu_{m}}{\mu_{S}(\nu_{m}+\mu_{S})}B\beta_{hm}\frac{I_{h}}{N_{h}}{\displaystyle F_{W,S}=\dfrac{\dfrac{\nu_{m}}{\mu_{S}(\nu_{m}+\mu_{S})}B\beta_{hm}\frac{I_{h}}{N_{h}}}{\mu_{S}+B\beta_{hm}\frac{I_{h}}{N_{h}}}}\dfrac{M+\varepsilon M_{S}^{*}}{M+M_{S}^{*}}r\gamma A
\]
Thus, from the two previous estimates, we deduce that
\[
\dfrac{S_{I}+F_{W,I}}{N_{h}}=\dfrac{\dfrac{\nu_{m}}{\mu_{S}(\nu_{m}+\mu_{S})}B\beta_{hm}\frac{I_{h}}{N_{h}}}{\mu_{S}+B\beta_{hm}\frac{I_{h}}{N_{h}}}{\displaystyle \dfrac{\left(\epsilon_F \Lambda_{tot}+r\gamma A\right)}{N_{h}}}.
\]
Then, from \eqref{solve_systH_EE}$_1$
\[
\mu_{h}N_{h}=B\beta_{mh}{\displaystyle \frac{F_{W,I}+S_{I}}{N_{h}}S_{h}+\mu_{h}S_{h}},
\]
and replacing $\dfrac{S_{I}+F_{W,I}}{N_{h}}$ leads to 
\[
\mu_{h}N_{h}={\displaystyle \left(B\beta_{mh}\dfrac{\nu_{m}}{\mu_{S}(\nu_{m}+\mu_{S})}\dfrac{B\beta_{hm}\frac{I_{h}}{N_{h}}}{\mu_{S}+B\beta_{hm}\frac{I_{h}}{N_{h}}}{\displaystyle \dfrac{\left(\epsilon_F \Lambda_{tot}+r\gamma A\right)}{N_{h}}}+\mu_{h}\right)S_{h},}
\]
\[
\mu_{h}N_{h}\left(\mu_{S}+B\beta_{hm}\frac{I_{h}}{N_{h}}\right)={\displaystyle \left(B\beta_{mh}\dfrac{\nu_{m}}{\mu_{S}(\nu_{m}+\mu_{S})}B\beta_{hm}\frac{I_{h}}{N_{h}}{\displaystyle \dfrac{\left(\epsilon_F \Lambda_{tot}+r\gamma A\right)}{N_{h}}}+\mu_{h}\left(\mu_{S}+B\beta_{hm}\frac{I_{h}}{N_{h}}\right)\right)S_{h},}
\]
from which we deduce
\[
S_{h}=\dfrac{\mu_{h}N_{h}\left(\mu_{S}+B\beta_{hm}\frac{I_{h}}{N_{h}}\right)}{{\displaystyle \left(B\beta_{mh}\dfrac{\nu_{m}}{\mu_{S}(\nu_{m}+\mu_{S})}B\beta_{hm}\frac{I_{h}}{N_{h}}{\displaystyle \dfrac{\left(\epsilon_F \Lambda_{tot}+r\gamma A\right)}{N_{h}}}+\mu_{h}\left(\mu_{S}+B\beta_{hm}\frac{I_{h}}{N_{h}}\right)\right)}}
\]
In particular, we can deduce
\[
\dfrac{S_{I}+F_{W,I}}{N_{h}}S_{h}={\displaystyle \dfrac{\left(\epsilon_F \Lambda_{tot}+r\gamma A\right)}{N_{h}}}\dfrac{\dfrac{\nu_{m}}{\mu_{S}(\nu_{m}+\mu_{S})}B\beta_{hm}\frac{I_{h}}{N_{h}}\mu_{h}N_{h}}{{\displaystyle \left(B\beta_{mh}\dfrac{\nu_{m}}{\mu_{S}(\nu_{m}+\mu_{S})}B\beta_{hm}\frac{I_{h}}{N_{h}}{\displaystyle \dfrac{\left(\epsilon_F \Lambda_{tot}+r\gamma A\right)}{N_{h}}}+\mu_{h}\left(\mu_{S}+B\beta_{hm}\frac{I_{h}}{N_{h}}\right)\right)}}
\]
and, using \eqref{solve_systH_EE}$_2$, i.e.
\[
B\beta_{mh}{\displaystyle \frac{F_{W,I}+S_{I}}{N_{h}}S_{h}=\left(\nu_h+\mu_{h}\right)I_{h},}
\]
we have
\[
B\beta_{mh}{\displaystyle {\displaystyle \dfrac{\left(\epsilon_F \Lambda_{tot}+r\gamma A\right)}{N_{h}}}\dfrac{\dfrac{\nu_{m}}{\mu_{S}(\nu_{m}+\mu_{S})}B\beta_{hm}\frac{I_{h}}{N_{h}}\mu_{h}}{{\displaystyle \left(B\beta_{mh}\dfrac{\nu_{m}}{\mu_{S}(\nu_{m}+\mu_{S})}B\beta_{hm}\frac{I_{h}}{N_{h}}{\displaystyle \dfrac{\left(\epsilon_F \Lambda_{tot}+r\gamma A\right)}{N_{h}}}+\mu_{h}\left(\mu_{S}+B\beta_{hm}\frac{I_{h}}{N_{h}}\right)\right)}}=\left(\nu_h+\mu_{h}\right)\dfrac{I_{h}}{N_{h}}.}
\]
Assuming $I_{h}>0$
\[
\begin{array}{l}
B\beta_{mh}{\displaystyle \dfrac{\left(\epsilon_F \Lambda_{tot}+r\gamma A\right)}{N_{h}}}\dfrac{\nu_{m}}{\mu_{S}(\nu_{m}+\mu_{S})}B\beta_{hm}\mu_{h}=\\
\left(\nu_h+\mu_{h}\right)\left(B\beta_{mh}\dfrac{\nu_{m}}{\mu_{S}(\nu_{m}+\mu_{S})}B\beta_{hm}\frac{I_{h}}{N_{h}}{\displaystyle \dfrac{\left(\epsilon_F \Lambda_{tot}+r\gamma A\right)}{N_{h}}}+\mu_{h}\left(\mu_{S}+B\beta_{hm}\frac{I_{h}}{N_{h}}\right)\right),
\end{array}
\]
we finally deduce
\[
I_{1,h}^{EE}=\mu_{h}\dfrac{B\beta_{mh}{\displaystyle \dfrac{\left(\epsilon_F \Lambda_{tot}+r\gamma A_{1}^{EE}\right)}{N_{h}}}\dfrac{\nu_{m}}{\mu_{S}(\nu_{m}+\mu_{S})}B\beta_{hm}-\left(\nu_h+\mu_{h}\right)\mu_{S}}{\left(\nu_h+\mu_{h}\right)B\beta_{hm}\left(\mu_{h}+B\beta_{mh}\dfrac{\nu_{m}}{\mu_{S}(\nu_{m}+\mu_{S})}{\displaystyle \dfrac{\left(\epsilon_F \Lambda_{tot}+r\gamma A_{1}^{EE}\right)}{N_{h}}}\right)}N_{h}>0,
\]
assuming
\[
{\displaystyle \dfrac{\left(\epsilon_F \Lambda_{tot}+r\gamma A_{1}^{EE}\right)}{N_{h}}}\dfrac{\nu_{m}B\beta_{hm}B\beta_{mh}}{\mu_{S}(\nu_{m}+\mu_{S})\left(\nu_h+\mu_{h}\right)\mu_{S}}>1,
\]
that is 
\begin{equation}
\epsilon_F \Lambda_{tot}+ r\gamma A_{1}^{EE}>\dfrac{F_{W,S}^{*}}{\R_{0}^{2}}.
\label{cond_EE}
\end{equation}
From the previous formulae, we deduce $S_{h,1}^{EE}$, $R_{h,1}^{EE}$, $F_{W,S,1}^{EE}$,
$F_{W,E,1}^{EE}$, $F_{W,I,1}^{EE}$, $S_{I,1}^{EE}$, $S_{E,1}^{EE}$, and finally $S_{S,1}^{EE}$. We proceed similarly to get the second endemic equilibrium $EE_{SIT,2}$ or $EE_{SIT,*}$, under the same condition \eqref{cond_EE} because $A_1^{EE}<A_1^{EE}$.

We now assume that $\mu_S<\mu_I$.  To derive the equilibria, such that $I_{h}>0$, $A>0$
and $S_{I}>0$, we have to solve
\begin{equation}
\label{solve_systH_EE2}
\left\{ %
\begin{array}{lcl}
0 & = & \mu_{h}N_{h}-B\beta_{mh}{\displaystyle \frac{F_{W,I}+S_{I}}{N_{h}}S_{h}-\mu_{h}S_{h},}\\
{\displaystyle 0} & = & B\beta_{mh}{\displaystyle \frac{F_{W,I}+S_{I}}{N_{h}}S_{h}-\nu_hI_{h}-\mu_{h}I_{h},}\\
0 & = & \nu_hI_{h}-\mu_{h}R_{h},
\end{array}\right.
\end{equation}
\begin{equation}
\label{solve_systM_EE2}
\left\{ %
\begin{array}{lcl}
0 & = & \phi(F_{W,S}+F_{W,E}+F_{W,I})-(\gamma+\mu_{A,1}+\mu_{A,2}A)A,\\
{\displaystyle 0} & = & (1-r)\gamma A-\mu_{M}M,\\
{\displaystyle 0} & = & \dfrac{M+\varepsilon M_{S}^{*}}{M+M_{S}^{*}}r\gamma A-B\beta_{hm}{\displaystyle \frac{I_{h}}{N_{h}}F_{W,S}-\mu_{S}F_{W,S},}\\
0 & = & B\beta_{hm}{\displaystyle \frac{I_{h}}{N_{h}}F_{W,S}-(\nu_{m}+\mu_{S})F_{W,E},}\\
0 & = & \nu_{m}F_{W,E}-\mu_{I}F_{W,I},\\
0 & = & \epsilon_F \Lambda_{tot}+\dfrac{(1-\varepsilon)M_{S}^{*}}{M+M_{S}^{*}}r\gamma A-B\beta_{hm}{\displaystyle \frac{I_{h}}{N_{h}}S_{S}-\mu_{S}S_{S},}\\
{\displaystyle 0} & = & B\beta_{hm}{\displaystyle \frac{I_{h}}{N_{h}}S_{S}-(\nu_{m}+\mu_{S})S_{E},}\\
{\displaystyle 0} & = & \nu_{m}S_{E}-\mu_{I}S_{I}.
\end{array}\right.
\end{equation}
Let us consider the auxiliary variable $X=\dfrac{M+\varepsilon M_S^*}{M+M_S^*}$. It follows from system \eqref{solve_systH_EE2}-\eqref{solve_systM_EE2} that:
$$
\begin{array}{ccl}
    F_{W,S} & = & \dfrac{r\gamma AX}{\mu_S+B\beta_{hm}\dfrac{I_h}{N_h}},  \\
     F_{W,E} & = & \dfrac{B\beta_{hm}}{\nu_m+\mu_S}\dfrac{I_h}{N_h}F_{W,S},\\  
      & = & \dfrac{B\beta_{hm}}{\nu_m+\mu_S}\dfrac{I_h}{N_h}\dfrac{r\gamma AX}{\mu_S+B\beta_{hm}\dfrac{I_h}{N_h}},  \\
      F_{W,I} & = & \dfrac{\nu_m}{\mu_I}F_{W,E},\\  
      & = & \dfrac{\nu_m}{\mu_I}\dfrac{B\beta_{hm}}{\nu_m+\mu_S}\dfrac{I_h}{N_h}\dfrac{r\gamma AX}{\mu_S+B\beta_{hm}\dfrac{I_h}{N_h}},  \\
\end{array}
$$
$$
\begin{array}{ccl}
      S_S &=&\dfrac{\epsilon_F \Lambda_{tot}+\dfrac{(1-\varepsilon)M_S^*}{M+M_S^*}r\gamma A}{\mu_S+B\beta_{hm}\dfrac{I_h}{N_h}},\\
      S_E &=& \dfrac{B\beta_{hm}}{\nu_m+\mu_S}\dfrac{I_h}{N_h}S_S,\\
      &=& \dfrac{B\beta_{hm}}{\nu_m+\mu_S}\dfrac{I_h}{N_h}\dfrac{\epsilon_F \Lambda_{tot}+\dfrac{(1-\varepsilon)M_S^*}{M+M_S^*}r\gamma A}{\mu_S+B\beta_{hm}\dfrac{I_h}{N_h}},\\
      S_I &=& \dfrac{\nu_m}{\mu_I}S_E,\\
      &=& \dfrac{\nu_m}{\mu_I}\dfrac{B\beta_{hm}}{\nu_m+\mu_S}\dfrac{I_h}{N_h}\dfrac{\epsilon_F \Lambda_{tot}+\dfrac{(1-\varepsilon)M_S^*}{M+M_S^*}r\gamma A}{\mu_S+B\beta_{hm}\dfrac{I_h}{N_h}}.\\
\end{array}
$$
Therefore, equation \eqref{solve_systM_EE2}$_1$ assumes the form $A=0$  or  
\begin{equation}\label{appenddix-aux1}
    r\gamma\phi X\left(1+\dfrac{B\beta_{hm}}{\nu_m+\mu_S}\dfrac{I_h}{N_h}\left(1+\dfrac{\nu_m}{\mu_I}\right)\right)-(\gamma+\mu_{A,1})\left(\mu_S+B\beta_{hm}\dfrac{I_h}{N_h}\right)-\mu_{A,2}\left(\mu_S+B\beta_{hm}\dfrac{I_h}{N_h}\right)A=0.
\end{equation}
\comment{
{\color{blue} 
$\bullet$ We first consider the case $A=0$. Reasoning as above, we obtain either the trivial disease-free equilibrium, TDFE or the wild insects-free endemic equilibrium $WIFEE=(S_h^{\sharp},I_h^{\sharp},R_h^{\sharp},0,0,0,0,0,S_S^{\sharp},S_E^{\sharp},S_I^{\sharp})$ whenever $\R_{0,TDFE}^2>\dfrac{\mu_S}{\mu_I}$, with
$$
    \begin{array}{ccl}
     S_I^\sharp  &= & \dfrac{\mu_hN_h}{B\beta_{mh}}\dfrac{1}{1+\dfrac{B\beta_{hm}\mu_h}{\mu_S(\nu_h+\mu_h)}}\left(\dfrac{\mu_I}{\mu_S}\R_{0,TDFE}^2-1\right), \\
    \end{array}
$$
$$S_S^\sharp = \dfrac{\epsilon_F \Lambda_{tot}}{\mu_S}\dfrac{1+\dfrac{\mu_S(\nu_m+\mu_S)N_h\mu_h}{\epsilon_F \Lambda_{tot}\nu_mB\beta_{mh}}}{1+\dfrac{B\beta_{hm}}{\mu_S}\dfrac{\mu_h}{\nu_h+\mu_h}}$$ and
$$
\begin{array}{ccl}
     S_E^\sharp &=& \dfrac{\mu_S}{\nu_m}S_I^\sharp,\\
         I_h^\sharp &=& \dfrac{B\beta_{mh}}{\nu_h+\mu_h}\dfrac{S_I^\sharp}{N_h}\dfrac{\mu_hN_h}{\mu_h+B\beta_{mh}\dfrac{S_I^\sharp}{N_h}},\\
         S_h^\sharp &=& \dfrac{\mu_hN_h}{\mu_h+B\beta_{mh}\dfrac{S_I^\sharp}{N_h}},\\
         R_h^\sharp &=& \dfrac{\nu_h}{\mu_h}I_h^\sharp.\\
\end{array}
$$

\noindent
$\bullet$ We now consider the case $A>0$.} 
}

From \eqref{solve_systH_EE2}$_2$, we derive $$I_h=\dfrac{B\beta_{mh}}{\nu_h+\mu_h}\dfrac{F_{W,I}+S_I}{N_h}.$$
However,

$$
\begin{array}{ccl}
     F_{W,I}+S_I &=& \dfrac{\nu_m}{\mu_I}\dfrac{B\beta_{hm}}{\nu_m+\mu_S}\dfrac{I_h}{N_h}\dfrac{r\gamma AX}{\mu_S+B\beta_{hm}\dfrac{I_h}{N_h}}+ \dfrac{\nu_m}{\mu_I}\dfrac{B\beta_{hm}}{\nu_m+\mu_S}\dfrac{I_h}{N_h}\dfrac{\Lambda_F+\dfrac{(1-\varepsilon)M_S^*}{M+M_S^*}r\gamma A}{\mu_S+B\beta_{hm}\dfrac{I_h}{N_h}},\\
     &=& \dfrac{\nu_m}{\mu_I}\dfrac{B\beta_{hm}}{\nu_m+\mu_S}\dfrac{I_h}{N_h}\dfrac{\Lambda_F+r\gamma A}{\mu_S+B\beta_{hm}\dfrac{I_h}{N_h}}.\\
\end{array}
$$
Therefore, for $I_h>0$, we have 

\begin{equation}\label{cas-general-I}
    \mu_S+B\beta_{hm}\dfrac{I_h}{N_h}=\dfrac{\nu_m}{\mu_I}\dfrac{B\beta_{mh}}{\nu_h+\mu_h}\dfrac{B\beta_{hm}}{\nu_m+\mu_S}\dfrac{1}{N_h^2}(\Lambda_F+r\gamma A)=\alpha (\Lambda_F+r\gamma A),
\end{equation}
where for simplicity, we set $$\alpha=\dfrac{\nu_m}{\mu_I}\dfrac{B\beta_{mh}}{\nu_h+\mu_h}\dfrac{B\beta_{hm}}{\nu_m+\mu_S}\dfrac{1}{N_h^2}.$$
Hence \eqref{appenddix-aux1} assumes the form
\begin{equation}\label{appenddix-aux2}r\gamma\phi \left(1+\dfrac{1+\dfrac{\nu_m}{\mu_I}}{\nu_m+\mu_S}(\alpha(\Lambda_F+r\gamma A)-\mu_S)\right)X-(\gamma+\mu_{A,1}+\mu_{A,2}A)(\Lambda_F+r\gamma A)\alpha=0
\end{equation}

or equivalently

\begin{equation}\label{appenddix-aux3}
    a_3A^3+a_2A^2+a_1A+a_0=0,
\end{equation}
where
$$
\begin{array}{ccl}
    a_3 & = & - \left( 1-r \right) {\gamma}^{2}\mu_{{A,2}}\alpha\,r<0,  \\
     a_2 & = & \dfrac {{r}^{2}{\gamma}^{3}\phi\, \left( 1-r \right)  \left( 1+{\dfrac {\nu_{{m
}}}{\mu_{{I}}}} \right) \alpha}{\nu_{{m}}+\mu_{{S}}}-\mu_{{M}}M_{{S}}^*\mu_{{A,2}}\alpha\,r\gamma- \left( 1-r \right) \gamma
\, \left(  \left( \gamma+\mu_{{A,1}} \right) r\gamma\,\alpha+\mu_{{A,2}}
\alpha\,\Lambda_{{F}} \right) 
,\\
& = & (1-r)\alpha\gamma\left(r\gamma(\gamma+\mu_{A,1})\left(\N\dfrac{1+\dfrac{\nu_m}{\mu_I}}{1+\dfrac{\nu_m}{\mu_S}}-\Q M_S^*-1\right)-\mu_{A,2}\Lambda_F\right),  
\end{array}
$$

$$
\begin{array}{ccl}
   a_1    & = & \dfrac{{r}^{2}{\gamma}^{2}\phi\,\mu_{{M}}\varepsilon\,M_{{S}}^* \left( 1+{\dfrac {
\nu_{{m}}}{\mu_{{I}}}} \right) \alpha}{\nu_{{m}}+\mu_{{S}}}+r{\gamma}^{2}\phi\, \left( 1-r \right)  \left( 1+
 \dfrac{\left( 1+{\dfrac {\nu_{{m}}}{\mu_{{I}}}} \right)  \left( \alpha\,
\Lambda_{{F}}-\mu_{{S}} \right)}{  \nu_{{m}}+\mu_{{S}} }\right)\\
&&-\mu_{{M}}M_{{S}}^* \left(  \left( \gamma+\mu_{{A,1}} \right) 
r\gamma\,\alpha+\mu_{{A,2}}\alpha\,\Lambda_{{F}} \right) - \left( 1-r
 \right) \gamma\, \left( \gamma+\mu_{{A,1}} \right) \Lambda_{{F}}\alpha,\\

 & = & \mu_MM_S^*r\gamma\alpha(\gamma+\mu_{A,1})\left(\N\varepsilon\dfrac{1+\dfrac{\nu_m}{\mu_I}}{1+\dfrac{\nu_m}{\mu_S}}-1\right)-(1-r)\gamma(\gamma+\mu_{A,1})\Lambda_F\alpha(\Q M_S^*+1)\\
  & & +r\gamma^2\phi(1-r)\left(1+\dfrac{1+\dfrac{\nu_m}{\mu_I}}{1+\dfrac{\nu_m}{\mu_S}}\left(\dfrac{\alpha\Lambda_F}{\mu_S}-1\right)\right),\\
 
 & = & (1-r)\gamma(\gamma+\mu_{A,1})\left(\dfrac{\alpha\mu_M M_S^*r}{1-r}\left(\N\varepsilon\dfrac{1+\dfrac{\nu_m}{\mu_I}}{1+\dfrac{\nu_m}{\mu_S}}-1\right)+\N\mu_S\left(1-\dfrac{1+\dfrac{\nu_m}{\mu_I}}{1+\dfrac{\nu_m}{\mu_S}}\right)+\alpha\Lambda_F\left(\N\dfrac{1+\dfrac{\nu_m}{\mu_I}}{1+\dfrac{\nu_m}{\mu_S}}-\Q M_S^*-1\right)\right),\\
\end{array}
$$

$$
\begin{array}{ccl}
  a_0 & = & \mu_MM_S^*\left(r\gamma\,\phi\,\varepsilon\,\left( 1+ \dfrac{\left( 1+{\dfrac {
\nu_{{m}}}{\mu_{{I}}}} \right)  \left( \alpha\,\Lambda_{{F}}-\mu_{{S}}
 \right) }{\nu_{{m}}+\mu_{{S}}}  \right) - \left( 
\gamma+\mu_{{A,1}} \right) \alpha\,\Lambda_{{F}}\right),\\
& = & \mu_MM_S^*\left( 
\gamma+\mu_{{A,1}} \right)\left(\mu_S\N\varepsilon\left(1+\dfrac{1+\dfrac{\nu_m}{\mu_I}}{1+\dfrac{\nu_m}{\mu_S}}\left(\dfrac{\alpha\Lambda_F}{\mu_S}-1\right)\right)-\alpha\Lambda_F\right),\\
&=& \mu_MM_S^*\left( 
\gamma+\mu_{{A,1}} \right)\mu_S\left(\N\varepsilon\left(1-\dfrac{1+\dfrac{\nu_m}{\mu_I}}{1+\dfrac{\nu_m}{\mu_S}}+\dfrac{1+\dfrac{\nu_m}{\mu_I}}{1+\dfrac{\nu_m}{\mu_S}}\dfrac{\alpha\Lambda_F}{\mu_S}\right)-\dfrac{\alpha\Lambda_F}{\mu_S}\right).\\
\end{array}
$$
Recall that since $\mu_S<\mu_I$, it follows that $$\dfrac{1+\dfrac{\nu_m}{\mu_I}}{1+\dfrac{\nu_m}{\mu_S}}<1.$$
To discuss the number of real positive solutions of equation \eqref{appenddix-aux3} we use the Descartes' rule of sign, see for instance Table \ref{table: descarte-sign-reule}.

\begin{table}[H]
    \centering
    \begin{tabular}{|c|c|c|c|c|}
    \hline
       $a_3$  & $a_2$ & $a_1$ & $a_0$ & Number of positive real solutions \\
       \hline
         - & - &- &- & 0 \\
         \hline
         - & - &- & + & 1 \\
         \hline
         - & - & + &- & 2\, \mbox{or}\,0 \\
         \hline
         - & + &- &- & 2\, \mbox{or}\,0  \\
         \hline
         - & + & + &- & 2\, \mbox{or}\,0  \\
         \hline
         - & + &- &+ & 3\, \mbox{or}\,1  \\
         \hline
         - & - &+ &+ & 1 \\
         \hline
         - & + & + & + & 1 \\
         \hline
    \end{tabular}
    \caption{Number of positive solutions of equation \eqref{appenddix-aux3} with the Descartes' rule of sign.}
    \label{table: descarte-sign-reule}
\end{table}

As explained in the numerical part, we consider a total release rate of sterile insects, $\Lambda_{tot}$, and a parameter $\epsilon_F$, the percentage of sterile females released, such that
$$
\Lambda_{M}=\left(1-\epsilon_F \right) \Lambda_{tot},\qquad \mbox{and} \qquad \Lambda_F=\epsilon_F \Lambda_{tot}.
$$
Doing like that, $M_S^*=(1-\epsilon_F)\dfrac{\Lambda_{tot}}{\mu_{M_{S}}}$ and we deduce all parameters thanks to $\Lambda_{tot}$, that is
$$
a_0=\mu_MM_S^*\left( 
\gamma+\mu_{{A,1}} \right)\mu_S \left(\N\varepsilon\left(1-\dfrac{1+\dfrac{\nu_m}{\mu_I}}{1+\dfrac{\nu_m}{\mu_S}}\right)-\dfrac{\alpha \epsilon_F\Lambda_{tot}}{\mu_S}\left(1-\N \varepsilon \left(\dfrac{1+\dfrac{\nu_m}{\mu_I}}{1+\dfrac{\nu_m}{\mu_S}} \right) \right) \right),
$$
$$
\begin{array}{ll}
a_1=&(1-r)\gamma(\gamma+\mu_{A,1})\left(\N\mu_{S}\left(1-\dfrac{1+\dfrac{\nu_{m}}{\mu_{I}}}{1+\dfrac{\nu_{m}}{\mu_{S}}}\right)-\left(\dfrac{\alpha\mu_{M}(1-\epsilon_F)r}{\left(1-r\right)\mu_{M_{S}}}\left(1-\N\varepsilon\dfrac{1+\dfrac{\nu_{m}}{\mu_{I}}}{1+\dfrac{\nu_{m}}{\mu_{S}}}\right)+\alpha\epsilon_F\left(1-\N\dfrac{1+\dfrac{\nu_{m}}{\mu_{I}}}{1+\dfrac{\nu_{m}}{\mu_{S}}}\right)\right)\Lambda_{tot} \right.\\
& \left.-\dfrac{\Q(1-\epsilon_F)}{\mu_{M_{S}}}\alpha\epsilon_F\left(\Lambda_{tot}\right)^{2}\right),
\end{array}
$$
and 
$$
a_2=(1-r)\alpha \gamma \left( r\gamma(\gamma+\mu_{A,1})\left(\N\dfrac{1+\dfrac{\nu_{m}}{\mu_{I}}}{1+\dfrac{\nu_{m}}{\mu_{S}}}-1\right)-\mu_{A,2}\left((1-\epsilon_F)\dfrac{r}{1-r}\dfrac{\mu_{M}}{\mu_{M_{S}}}+\epsilon_F\right)\Lambda_{tot} \right).
$$
From $a_0$, it is easy to deduce the following discussion thanks to $\N \varepsilon$ and for a given $\epsilon_F$:
\begin{itemize}
    \item $a_0>0$ if $\N \varepsilon \geq \dfrac{1+\dfrac{\nu_m}{\mu_S}}{1+\dfrac{\nu_m}{\mu_I}}>1$ or if $\N \varepsilon \leq \dfrac{1+\dfrac{\nu_m}{\mu_S}}{1+\dfrac{\nu_m}{\mu_I}}$ and
    $$
    \Lambda_{tot} \leq \Lambda_{F,EE}^{crit,1}=\dfrac{ \mu_S}{\epsilon_F\alpha}\dfrac{\N\varepsilon\left(1-\dfrac{1+\dfrac{\nu_m}{\mu_I}}{1+\dfrac{\nu_m}{\mu_S}}\right)}{1-\N \varepsilon \left(\dfrac{1+\dfrac{\nu_m}{\mu_I}}{1+\dfrac{\nu_m}{\mu_S}} \right)}.
    $$
    Otherwise, when $\Lambda_{tot} > \Lambda_{F,EE}^{crit,1}$, $a_0<0$
    \item It is interesting to notice that $a_2<0$, whatever $\Lambda_{tot}\geq 0$ if $\N \leq \dfrac{1+\dfrac{\nu_m}{\mu_S}}{1+\dfrac{\nu_m}{\mu_I}}$. In addition $\N \leq \dfrac{1+\dfrac{\nu_m}{\mu_S}}{1+\dfrac{\nu_m}{\mu_I}}$ implies $\N \varepsilon \leq \dfrac{1+\dfrac{\nu_m}{\mu_S}}{1+\dfrac{\nu_m}{\mu_I}}$ because $0\leq \varepsilon <1$. When $\N > \dfrac{1+\dfrac{\nu_m}{\mu_S}}{1+\dfrac{\nu_m}{\mu_I}}$, then $a_2>0$ if
    $$
    \Lambda_{tot} < \Lambda_{tot}^{crit,2}=\dfrac{r\gamma(\gamma+\mu_{A,1})\left(\N\dfrac{1+\dfrac{\nu_{m}}{\mu_{I}}}{1+\dfrac{\nu_{m}}{\mu_{S}}}-1\right)}{\mu_{A,2}\left((1-\epsilon_F)\dfrac{r}{1-r}\dfrac{\mu_{M}}{\mu_{M_{S}}}+\epsilon_F\right)}.
    $$
    It is negative, otherwise.
    \item Straightforward computations show that $a_1>0$ if 
    $$
    \Lambda_{tot}< \Lambda_{tot}^{crit,3}=\dfrac{1}{2\dfrac{\Q(1-\epsilon_F)}{\mu_{M_{S}}}\alpha\epsilon_F}\left[\sqrt{\Delta}+\left(\dfrac{\alpha\mu_{M}(1-\epsilon_F)r}{\left(1-r\right)\mu_{M_{S}}}\left(1-\N\varepsilon\dfrac{1+\dfrac{\nu_{m}}{\mu_{I}}}{1+\dfrac{\nu_{m}}{\mu_{S}}}\right)+\alpha\epsilon_F\left(1-\N\dfrac{1+\dfrac{\nu_{m}}{\mu_{I}}}{1+\dfrac{\nu_{m}}{\mu_{S}}}\right)\right)\right],
    $$
    where 
    $$
    \Delta=\left(\left(\dfrac{\alpha\mu_{M}(1-\epsilon_F)r}{\left(1-r\right)\mu_{M_{S}}}\left(1-\N\varepsilon\dfrac{1+\dfrac{\nu_{m}}{\mu_{I}}}{1+\dfrac{\nu_{m}}{\mu_{S}}}\right)+\alpha\epsilon_F\left(1-\N\dfrac{1+\dfrac{\nu_{m}}{\mu_{I}}}{1+\dfrac{\nu_{m}}{\mu_{S}}}\right)\right)\right)^{2}+4\dfrac{\Q(1-\epsilon_F)}{\mu_{M_{S}}}\alpha\epsilon_F\N\mu_{S}\left(1-\dfrac{1+\dfrac{\nu_{m}}{\mu_{I}}}{1+\dfrac{\nu_{m}}{\mu_{S}}}\right)>0.
    $$
    \item Thus, we derive 
    \begin{itemize}
        \item Assume $\N \varepsilon \geq \dfrac{1+\dfrac{\nu_m}{\mu_S}}{1+\dfrac{\nu_m}{\mu_I}}$. If $\Lambda_{tot}< \Lambda_{tot}^{crit,3}$, then $a_0>0$ and $a_1>0$,
        \item Assume $\N \varepsilon \geq \dfrac{1+\dfrac{\nu_m}{\mu_S}}{1+\dfrac{\nu_m}{\mu_I}}$. If $\Lambda_{tot}> \max \{\Lambda_{tot}^{crit,2}, \Lambda_{tot}^{crit,3} \}$, then $a_0>0$, $a_1<0$, and $a_2<0$ because $\N \varepsilon \geq \dfrac{1+\dfrac{\nu_m}{\mu_S}}{1+\dfrac{\nu_m}{\mu_I}}$ implies $\N \geq \dfrac{1+\dfrac{\nu_m}{\mu_S}}{1+\dfrac{\nu_m}{\mu_I}}$,
        \item Assume $\N \varepsilon \leq  \dfrac{1+\dfrac{\nu_m}{\mu_S}}{1+\dfrac{\nu_m}{\mu_I}}$ and $\Lambda_{tot}<\min \{\Lambda_{tot}^{crit,3},\Lambda_{tot}^{crit,1} \}$, then $a_0>0$ and $a_1>0$,
    \end{itemize}
    such that, thanks to Table \ref{table: descarte-sign-reule} page \pageref{table: descarte-sign-reule}, we deduce that only one positive equilibrium exists.
    \begin{itemize}
        \item Assume $\N \varepsilon \geq \dfrac{1+\dfrac{\nu_m}{\mu_S}}{1+\dfrac{\nu_m}{\mu_I}}$. If $\Lambda_{tot}^{crit,3} <\Lambda_{tot} <\Lambda_{tot}^{crit,2}$, then $a_0>0$, $a_1<0$, and $a_2>0$,
    \end{itemize}
    such that, thanks to Table \ref{table: descarte-sign-reule}, we deduce that there exists $1$ or $3$ positive equilibria.
    \item Assume $\N \varepsilon \leq  \dfrac{1+\dfrac{\nu_m}{\mu_S}}{1+\dfrac{\nu_m}{\mu_I}}$. 
    \begin{itemize}
        \item If $\Lambda_{tot}^{crit,1}<\Lambda_{tot}< \Lambda_{tot}^{crit,3}$, then $a_0<0$ and $a_1>0$,
        \item If $\Lambda_{tot}>\max \{\Lambda_{tot}^{crit,3},\Lambda_{tot}^{crit,1} \}$, then $a_0<0$ and $a_1<0$. If $\N \geq \dfrac{1+\dfrac{\nu_m}{\mu_S}}{1+\dfrac{\nu_m}{\mu_I}}$ and $\Lambda_{tot} < \Lambda_{tot}^{crit,2}$, then $a_2>0$,
    \end{itemize}
     such that, thanks to Table \ref{table: descarte-sign-reule}, whatever the sign of $a_2$, we deduce that no or $2$ positive equilibria.
     \item Assume $\N \varepsilon \leq  \dfrac{1+\dfrac{\nu_m}{\mu_S}}{1+\dfrac{\nu_m}{\mu_I}}$. If $\Lambda_{tot}>\max \{\Lambda_{tot}^{crit,3},\Lambda_{tot}^{crit,1} \}$, then $a_0<0$ and $a_1<0$.
     \begin{itemize}
         \item  If $\N \leq \dfrac{1+\dfrac{\nu_m}{\mu_S}}{1+\dfrac{\nu_m}{\mu_I}}$, then $a_2<0$,
        \item If $\N \geq \dfrac{1+\dfrac{\nu_m}{\mu_S}}{1+\dfrac{\nu_m}{\mu_I}}$, and $\Lambda_{tot}> \Lambda_{tot}^{crit,2}$, then $a_2<0$, 
     \end{itemize}
     such that, according to Table \ref{table: descarte-sign-reule} page \pageref{table: descarte-sign-reule}, there is no positive equilibrium.
     \end{itemize}

\section{Proof of Theorem  \ref{Persistence}}\label{AppendixE}
Assume that $\N\varepsilon>1$ and $\R_{0,TDFE}^2>1$. Let us consider the following sets
$$
\begin{array}{ccl}
    x(t) & := & (S_h,I_h,R_h,A,M,F_{W,S},F_{W,E},F_{W,I},S_S,S_E,S_I)(t),\\
    \Gamma & := & \{ x\in \mathbb{R}^{11}_+: I_h>0, R_h>0,A>0,M>0,F_{W,S}>0,F_{W,E}>0,F_{W,I}>0, S_E>0,S_I>0 \}, \\
   \partial\Gamma & := & \{ x\in \mathbb{R}^{11}_+: I_h\times 
 R_h\times A\times M\times F_{W,S}\times F_{W,E}\times F_{W,I}\times S_E\times S_I=0 \}. 
\end{array}
$$
Direct computations, see e.g. \cite{Traore2021}, lead that the sets $\Gamma$ and $\partial\Gamma$ are positively invariant with respect to system \eqref{Human-ode-good}-\eqref{Mosquitoes-ode-SIT-constant}. All solutions are bounded and system \eqref{Human-ode-good}-\eqref{Mosquitoes-ode-SIT-constant} is a point dissipative system. We denote $\varphi_t(x_0)$ the flow corresponding to system \eqref{Human-ode-good}-\eqref{Mosquitoes-ode-SIT-constant}, such that the solution of system \eqref{Human-ode-good}-\eqref{Mosquitoes-ode-SIT-constant} starting at $x_0$ at $t>0$ is $x(t,x_0)=\varphi_t(x_0)$. Let $M_\partial=\{x\in\partial\Gamma: \varphi_t(x)\in\partial\Gamma ~\mbox{for}~ t\geq0 \}$. Then we have $M_\partial=\partial\Gamma$. The trivial disease-free equilibrium $TDFE$, the wild insects-free equilibrium $WIFE$ and the disease-free equilibrium $DFE$ are in $M_\partial$. Let $W^s(TDFE)$, $W^s(WIFE)$ and $W^s(DFE)$ be the stable manifold of $TDFE$, $WIFE$ and $DFE$, respectively. In the sequel, we prove that $W^s(TDFE)\cap\Gamma=\emptyset$, $W^s(WIFE)\cap\Gamma=\emptyset$ and $W^s(DFE)\cap\Gamma=\emptyset$ hold when $\N\varepsilon>1$ and $\R_{0,TDFE}^2>1$. We first show that $W^s(TDFE)\cap\Gamma=\emptyset$. Since $\N\varepsilon>1$ and $\R_{0,TDFE}^2>1$, by continuity, there exists $\epsilon_0$ such that for all $\epsilon\in[0,\epsilon_0]$, we have
$$\dfrac{r\gamma\phi}{\left(\mu_S+\dfrac{B\beta_{hm}}{N_h}\epsilon\right)(\gamma+\mu_{A,1}+\mu_{A,2}\epsilon)}\varepsilon>1$$

and

$$
\dfrac{\nu_m}{\nu_m+\mu_S}\dfrac{B\beta_{mh}}{\mu_I}\dfrac{B\beta_{hm}}{\nu_h+\mu_h}\dfrac{1}{N_h^2}\left(\dfrac{\epsilon_F\Lambda_{tot}}{\mu_S}-\epsilon\right)\left(N_h-\epsilon\right)>1.$$

We claim that there exists $\eta_0>0$, such that for all $x_0\in\Gamma$, $\lim\sup\limits_{t\rightarrow+\infty}\|\varphi_t(x_0)-TDFE\|>\eta_0$. Indeed, suppose that this is not true. Hence, there exists $T>0$ such that for $t>T$, we have:
$$N_h-\epsilon\leq S_h\leq N_h+\epsilon, \quad \dfrac{\epsilon_F\Lambda_{tot}}{\mu_S}-\epsilon \leq S_S \leq \dfrac{\epsilon_F\Lambda_{tot}}{\mu_S}+\epsilon,\quad I_h\leq \epsilon, \quad A\leq \epsilon. 
$$
From system \eqref{Human-ode-good}-\eqref{Mosquitoes-ode-SIT-constant}, it follows that

\begin{equation}\label{sous-system1}    
\left\{
\begin{array}{ccl}
     \dfrac{dA}{dt} &\geq& \phi F_{W,S} -(\gamma+\mu_{A,1}+\mu_{A,2}\epsilon)A,  \\
     \dfrac{dF_{W,S}}{dt} &\geq&  \varepsilon r\gamma A- B\dfrac{\beta_{hm}\epsilon}{N_h}F_{W,S}-\mu_SF_{W,S},\\
\end{array}
\right.
\end{equation}

and

\begin{equation}\label{sous-system2} 
\left\{
\begin{array}{lcl}
{\displaystyle \frac{dI_{h}}{dt}} & \geq & B\beta_{mh}{\displaystyle \frac{S_{I}}{N_{h}}(N_h-\epsilon)-\nu_hI_{h}-\mu_{h}I_{h},}\\
\dfrac{dR_{h}}{dt} & = & \nu_hI_{h}-\mu_{h}R_{h},\\
{\displaystyle \frac{dS_{E}}{dt}} & \geq & B\beta_{hm}{\displaystyle \frac{I_{h}}{N_{h}}\left(\dfrac{\epsilon_F\Lambda_{tot}}{\mu_S}-\epsilon\right)-(\nu_{m}+\mu_{S})S_{E},}\\
{\displaystyle \frac{dS_{I}}{dt}} & = & \nu_{m}S_{E}-\mu_{I}S_{I}.
\end{array}
\right.
\end{equation}

Let us consider the matrices

$$J_1=\begin{pmatrix}
    -(\gamma+\mu_{A,1}+\mu_{A,2}\epsilon) & \phi \\
    \varepsilon r\gamma & -B\dfrac{\beta_{hm}\epsilon}{N_h}F_{W,S}-\mu_S
\end{pmatrix}$$
and 
$$
J_2=\begin{pmatrix}
    -(\nu_h+\mu_h) & 0 &  0 &  B\beta_{mh}\dfrac{N_h-\epsilon}{N_h} \\
    \nu_h &-\mu_h & 0 & 0\\
\displaystyle \frac{B\beta_{hm}}{N_{h}}\left(\dfrac{\epsilon_F\Lambda_{tot}}{\mu_S}-\epsilon\right) & 0 & -(\nu_m+\mu_S) & 0 \\
0 & 0 & \nu_m & -\mu_I\\
\end{pmatrix}.
$$
Let $s(J)$ be the stability modulus of the matrix $J$. It therefore follows that $s(J_1)>0$ and $s(J_2)>0$. Hence, the positive solutions $A$, $F_{W,S}$, $I_h$, $R_h$, $S_E$ and $S_I$ of systems \eqref{Human-ode-good}-\eqref{Mosquitoes-ode-SIT-constant} are unbounded which is a contradiction. Thus, $W^s(TDFE)\cap\Gamma=\emptyset$. Exactly the same computations show also that $W^s(WIFE)\cap\Gamma=\emptyset$. To show that $W^s(DFE)\cap\Gamma=\emptyset$, we first recall that following \eqref{R0-SIT-d}, $\R_{0,TDFE}^2>1 \Rightarrow \R_{0,SIT_c}^2>1$. Hence by continuity there exists $\epsilon_0$ such that for all $\epsilon\in[0,\epsilon_0]$, we also have

$$
\dfrac{\nu_m}{\nu_m+\mu_S}\dfrac{B\beta_{mh}}{\mu_I}\dfrac{B\beta_{hm}}{\nu_h+\mu_h}\dfrac{1}{N_h^2}\left(S_{S,DFE}-\epsilon\right)\left(N_h-\epsilon\right)>1.$$
As before, we claim that there exists $\eta_0>0$, such that for all $x_0\in\Gamma$, $\lim\sup\limits_{t\rightarrow+\infty}\|\varphi_t(x_0)-TDFE\|>\eta_0$. Indeed, suppose that this is not true. Hence, there exists $T>0$ such that for $t>T$, we have:
$$N_h-\epsilon\leq S_h\leq N_h+\epsilon, \quad S_{S,DFE}-\epsilon \leq S_S \leq S_{S,DFE}+\epsilon,\quad I_h\leq \epsilon, \quad A\leq \epsilon. 
$$
Hence, \eqref{sous-system2}$_3$ assumes the form $${\displaystyle \frac{dS_{E}}{dt}}  \geq  B\beta_{hm}{\displaystyle \frac{I_{h}}{N_{h}}\left(S_{S,DFE}-\epsilon\right)-(\nu_{m}+\mu_{S})S_{E},}$$
and matrix $J_2$ now becomes
$$
J_2=\begin{pmatrix}
    -(\nu_h+\mu_h) & 0 &  0 &  B\beta_{mh}\dfrac{N_h-\epsilon}{N_h} \\
    \nu_h &-\mu_h & 0 & 0\\
\displaystyle \frac{B\beta_{hm}}{N_{h}}\left(S_{S,DFE}-\epsilon\right) & 0 & -(\nu_m+\mu_S) & 0 \\
0 & 0 & \nu_m & -\mu_I\\
\end{pmatrix}.
$$
As previously, we have that $s(J_1)>0$ and $s(J_2)>0$. Hence, the positive solutions $A$, $F_{W,S}$, $I_h$, $R_h$, $S_E$ and $S_I$ of system \eqref{Human-ode-good}-\eqref{Mosquitoes-ode-SIT-constant} are unbounded which a contradiction. Thus, $W^s(DFE)\cap\Gamma=\emptyset$. Therefore, we have 
$W^s(TDFE)=\{TDFE\}$, $W^s(WIFE)=\{x\in\mathbb{R}_+^{11}: A=M=F_{W,S}=F_{W,E}=F_{W,I}=0, I_h>0, R_h>0, S_E>0, S_I>0\}$ and $W^s(DFE)=\{x\in\mathbb{R}_+^{11}: A>0, M>0, F_{W,S}>0,  F_{W,E}=F_{W,I}=I_h==S_E=S_I=0\}$ such that $M_\partial=W^s(TDFE)\cup W^s(WIFE)\cup W^s(DFE)$. In addition, each equilibrium is isolated and acyclic in $M_\partial$. Based on Theorem \cite[Theorem 4.6]{Thieme1992}, we found that system \eqref{Human-ode-good}-\eqref{Mosquitoes-ode-SIT-constant} is uniformly persistent with respect to $(\Gamma,\partial\Gamma)$ whenever $\N\varepsilon>1$ and $\R_{0,TDFE}^2>1$. Moreover, using the invariance of $\Gamma$, the dissipativity of system \eqref{Human-ode-good}-\eqref{Mosquitoes-ode-SIT-constant} and its uniform persistence, we can deduce, following \cite[Theorem D.3]{smith-waltman-1995}, the existence of a least one positive coexistence equilibrium.



    

\end{document}